\documentclass[oneside,english]{amsart}
\usepackage[T1]{fontenc}
\usepackage[latin9]{inputenc}
\usepackage{geometry}
\usepackage{array}
\usepackage{longtable}
\usepackage{float}
\usepackage{amstext}
\usepackage{amsthm}
\usepackage{amssymb}
\usepackage{stmaryrd}
\usepackage{setspace}
\usepackage{wasysym}

\makeatletter

\providecommand{\tabularnewline}{\\}

\numberwithin{equation}{section}
\numberwithin{figure}{section}
\theoremstyle{plain}
\newtheorem{thm}{\protect\theoremname}
\theoremstyle{remark}
\newtheorem{rem}[thm]{\protect\remarkname}
\theoremstyle{plain}
\newtheorem{lem}[thm]{\protect\lemmaname}
\theoremstyle{plain}
\newtheorem{prop}[thm]{\protect\propositionname}
\theoremstyle{plain}
\newtheorem{cor}[thm]{\protect\corollaryname}

\makeatother

\usepackage{babel}
\providecommand{\corollaryname}{Corollary}
\providecommand{\lemmaname}{Lemma}
\providecommand{\propositionname}{Proposition}
\providecommand{\remarkname}{Remark}
\providecommand{\theoremname}{Theorem}

\begin{document}
\title{Reachable elements in basic classical Lie superalgebras}
\author{Leyu Han}
\begin{abstract}
\noindent Let $\mathfrak{g}=\mathfrak{g}_{\bar{0}}\oplus\mathfrak{g}_{\bar{1}}$
be a basic classical Lie superalgebra over $\mathbb{C}$, $e\in\mathfrak{g}_{\bar{0}}$
a nilpotent element and $\mathfrak{g}^{e}$ the centralizer of $e$
in $\mathfrak{g}$. We study various properties of nilpotent elements
in $\mathfrak{g}$, which have previously only been considered in
the case of Lie algebras. In particular, we prove that $e$ is reachable
if and only if $e$ satisfies the Panyushev property for $\mathfrak{g}=\mathfrak{sl}(m|n)$,
$m\neq n$ or $\mathfrak{psl}(n|n)$ and $\mathfrak{osp}(m|2n)$.
For exceptional Lie superalgebras $\mathfrak{g}=D(2,1;\alpha)$, $G(3)$,
$F(4)$, we give the classification of $e$ which are reachable, strongly
reachable or satisfy the Panyushev property. In addition, we give
bases for $\mathfrak{g}^{e}$ and its centre $\mathfrak{z}(\mathfrak{g}^{e})$
for $\mathfrak{g}=\mathfrak{psl}(n|n)$, which completes results of
Han on the relationship between $\dim\mathfrak{g}^{e}$, $\dim\mathfrak{z}(\mathfrak{g}^{e})$
and the labelled Dynkin diagrams for all basic classical Lie superalgebras.{\normalsize{} }{\normalsize\par}
\end{abstract}

\maketitle

\section{Introduction\label{sec:Introduction}}

\noindent Let $\mathfrak{g}=\mathfrak{g}_{\bar{0}}\oplus\mathfrak{g}_{\bar{1}}$
be a basic classical Lie superalgebra over $\mathbb{C}$ and let $e\in\mathfrak{g}_{\bar{0}}$.
Write $\mathfrak{g}^{e}=\{x\in\mathfrak{g}:[e,x]=0\}$ for the centralizer
of an element $e$ in $\mathfrak{g}$. The main objects of our study
here are $\mathfrak{g}^{e}$ and its derived algebra $[\mathfrak{g}^{e},\mathfrak{g}^{e}]$.
The element $e$ is called \textit{reachable} if $e\in[\mathfrak{g}^{e},\mathfrak{g}^{e}]$.
The element $e$ is called \textit{strongly reachable} if $[\mathfrak{g}^{e},\mathfrak{g}^{e}]=\mathfrak{g}^{e}$.
Using the Jordan decomposition, it can be shown that only a nilpotent
element can satisfies any of the above conditions, see \cite[Section 1]{Yakimova}.
Throughout this paper, let $e\in\mathfrak{g}_{\bar{0}}$ be nilpotent
and let $\{e,h,f\}$ be an $\mathfrak{sl}(2)$-triple in $\mathfrak{g}_{\bar{0}}$.
We say that $e$ satisfies the \textit{Panyushev property} if, in
the ad$h$ grading $\mathfrak{g}^{e}=\bigoplus_{j\geq0}\mathfrak{g}^{e}(j)$,
the subalgebra $\mathfrak{g}^{e}(\geq1)=\bigoplus_{j\geq1}\mathfrak{g}^{e}(j)$
is generated by $\mathfrak{g}^{e}(1)$. Note that the Panyushev property
is equivalent to $\mathfrak{g}^{e}$ being generated by $\mathfrak{g}^{e}(0)$
and $\mathfrak{g}^{e}(1)$.

Research on reachable elements in the case of Lie algebras has been
developed over the years since Elashvili and Grelaud \cite{Elashvili=000026Gr=0000E9laud}
classified all reachable elements in simple Lie algebras, such elements
are called compact in \cite{Elashvili=000026Gr=0000E9laud}. Many
mathematicians studied the properties of reachable nilpotent elements
in different types of Lie algebras $\mathfrak{g}$ after that. In
\cite{Panyushev}, Panyushev showed that in case of $\mathfrak{g}=\mathfrak{sl}(n)$,
a nilpotent element $e$ is reachable is equivalent to $e$ satisfies
the Panyushev property. Yakimova \cite{Yakimova} extended Panyushev's
results to cases in which $\mathfrak{g}=\mathfrak{so}(n)$ and $\mathfrak{sp}(n)$
and further determined the condition for $e$ such that it is strongly
reachable. Properties of nilpotent orbits in exceptional Lie algebras
over $\mathbb{C}$ were investigated by de Graaf in \cite{de Graaf}
and were extended to positive characteristic by Premet and Stewart
\cite{Premet=000026Stewart}. 

It is natural to ask what are the conditions for nilpotent orbits
in Lie superalgebras to be reachable, strongly reachable or satisfy
the Panyushev property. This was a question posed by V.G. Kac. The
purpose of this paper is to classify all reachable nilpotent elements
in $\mathfrak{g}_{\bar{0}}$ for all types of basic classical Lie
superalgebras and study the relationship between the property of being
reachable, strongly reachable and the Panyushev property. In particular,
we prove that $e$ is reachable if and only if $e$ satisfies the
Panyushev property for $\mathfrak{g}=\mathfrak{sl}(m|n)$, $m\neq n$,
$\mathfrak{psl}(n|n)$ and $\mathfrak{osp}(m|2n)$. For exceptional
Lie superalgebras $\mathfrak{g}=D(2,1;\alpha)$, $G(3)$, $F(4)$,
we apply results on the structure of $\mathfrak{g}^{e}$ in \cite{han-exc}
to give the classification of $e$ which are reachable, strongly reachable
or satisfy the Panyushev property. In contrary to Lie algebra case,
there are two cases for which $e$ are reachable but $e$ do not satisfy
the Panyushev property. 

The main results of our paper are as follows:

For Lie superalgebras $\mathfrak{g}=\mathfrak{sl}(m|n)$, $m\neq n$,
$\mathfrak{psl}(n|n)$ or $\mathfrak{osp}(m|2n)$, the nilpotent orbits
are parametrized by partitions. Our first result can be viewed as
an analogue of the result of Elashvili and Grelaud \cite[Theorem 1]{Elashvili=000026Gr=0000E9laud}
which classifies all reachable elements in $\mathfrak{g}$ in terms
of partitions.
\begin{thm}
\label{thm:Re iff partition}Let $\mathfrak{g}=\mathfrak{g}_{\bar{0}}\oplus\mathfrak{g}_{\bar{1}}=\mathfrak{sl}(m|n)$,
$m\neq n$, $\mathfrak{psl}(n|n)$ or $\mathfrak{osp}(m|2n)$ and
$e\in\mathfrak{g}_{\bar{0}}$ be nilpotent. Let $\lambda=(\lambda_{1},\dots,\lambda_{r+s})$
be a partition of $(m|n)$ (resp. $(m|2n)$) with respect to $e$
for $\mathfrak{g}=\mathfrak{sl}(m|n)$, $m\neq n$ or $\mathfrak{psl}(n|n)$
(resp. $\mathfrak{g}=\mathfrak{osp}(m|2n)$) such that $\lambda_{1}\geq\dots\geq\lambda_{r+s}$.
Then $e$ is reachable if and only if $\lambda_{i}-\lambda_{i+1}\in\{0,1\}$
for all $1\leq i<r+s$ and $\lambda_{r+s}=1$.
\end{thm}

Our next theorem is an analogue of results of Panyushev \cite[Theorem 4.5]{Panyushev}
and Yakimova \cite[Lemma 10]{Yakimova} for Lie superalgebras.
\begin{thm}
\label{thm:Re iff Pan}Let $\mathfrak{g}=\mathfrak{g}_{\bar{0}}\oplus\mathfrak{g}_{\bar{1}}=\mathfrak{sl}(m|n)$,
$m\neq n$, $\mathfrak{psl}(n|n)$ or $\mathfrak{osp}(m|2n)$ and
$e\in\mathfrak{g}_{\bar{0}}$ be nilpotent. Then $e$ is reachable
if and only if $\mathfrak{g}^{e}(\geq1)=\bigoplus_{j\geq1}\mathfrak{g}^{e}(j)$
is generated by $\mathfrak{g}^{e}(1)$, i.e. $[\mathfrak{g}^{e}(1),\mathfrak{g}^{e}(j)]=\mathfrak{g}^{e}(j+1)$
for each $j\geq1$.
\end{thm}

For exceptional Lie superalgebras $\mathfrak{g}$, the author described
the structure of the centralizer $\mathfrak{g}^{e}$ in \cite{han-exc},
in particular, its dimension and a basis. We use those results to
prove the following theorem, for which we note that some notation
for nilpotent orbits is introduced in Section \ref{sec:exceptional}.
\begin{thm}
\label{thm:Re iff Pan exceptional}Let $\mathfrak{g}=\mathfrak{g}_{\bar{0}}\oplus\mathfrak{g}_{\bar{1}}=D(2,1;\alpha)$,
$G(3)$ or $F(4)$ and $e\in\mathfrak{g}_{\bar{0}}$ be nilpotent.
Then $e$ is reachable if and only if $e$ satisfies the Panyushev
property except for $\mathfrak{g}=G(3)$ and $e=x_{1}$ or $e=E+x_{1}$,
in which cases $e$ are reachable but the Panyushev property do not
hold. The nilpotent orbits that are reachable or strongly reachable
are listed in Tables \ref{tab:Reachable D(2,1;)}, \ref{tab:Reachable G(3)}
and \ref{tab:Reachable F(4)}.
\end{thm}

In order to study the structure of $[\mathfrak{g}^{e},\mathfrak{g}^{e}]$
for the proof of Theorems \ref{thm:Re iff partition} and \ref{thm:Re iff Pan},
we need to investigate the structure for $\mathfrak{g}^{e}$. A basis
for $\mathfrak{g}^{e}$ is given in \cite[Sections 3--4]{Han-sl-osp}
for $\mathfrak{g}=\mathfrak{sl}(m|n)$, $m\neq n$ or $\mathfrak{osp}(m|2n)$
and in \cite[Sections 4--6]{han-exc} for $\mathfrak{g}=D(2,1;\alpha)$,
$G(3)$ or $F(4)$. Hence, it remains to deal with the case of $\mathfrak{g}=\mathfrak{psl}(n|n)$.
In Section \ref{sec:centralizer =000026 centre psl}, we give bases
for $\mathfrak{g}^{e}$ and its centre $\mathfrak{z}(\mathfrak{g}^{e})$
for $\mathfrak{g}=\mathfrak{psl}(n|n)$. This also leads to the following
theorems that completes the results in \cite{Han-sl-osp} and \cite{han-exc}
for all basic classical Lie superalgebras. 

Combining with \cite{Han-sl-osp} and \cite{han-exc}, we state final
results relating $\dim\mathfrak{g}^{e}$, $\dim\mathfrak{z}(\mathfrak{g}^{e})$
and the labelled Dynkin diagrams in Theorems \ref{thm:LLD=000026centre}--\ref{thm:g^e,z(g^e)}.
We fix some notation as below. Let $G$ be the reductive algebraic
group over $\mathbb{C}$ given as in Table \ref{tab:Linear-algebraic-groups}.
We have that $\mathrm{Lie}(G)=\mathfrak{g}_{\bar{0}}$ and there is
a representation $\rho:G\rightarrow\mathrm{GL}(\mathfrak{g}_{\bar{1}})$
such that $d_{\rho}:\mathrm{Lie}(G)\rightarrow\mathfrak{gl}(\mathfrak{g}_{\bar{1}})$
determines the adjoint action of $\mathfrak{g}_{\bar{0}}$ on $\mathfrak{g}_{\bar{1}}$.
Fix $\Delta$ to be a labelled Dynkin diagram with respect to $e$.
Let $n_{2}(\Delta)$ be the number of nodes which have labels equal
to $2$ in $\Delta$. For $\mathfrak{g}=\mathfrak{sl}(m|n)$, $m\neq n$,
$\mathfrak{psl}(n|n)$ or $\mathfrak{osp}(m|2n)$, suppose $P$ is
the Dynkin pyramid of shape $\lambda$ which will be defined in Subsection
\ref{subsec:Dynkin-pyramids}. Let $r_{i}$ (resp. $s_{i}$) be the
number of boxes on the $i$th column with parity $\bar{0}$ (resp.
$\bar{1}$) in $P$.
\begin{thm}
\label{thm:LLD=000026centre}Let $\mathfrak{g}=\mathfrak{g}_{\bar{0}}\oplus\mathfrak{g}_{\bar{1}}$
be a basic classical Lie superalgebra and $e\in\mathfrak{g}_{\bar{0}}$
be nilpotent. Let $G$ be the reductive algebraic group defined as
in Table \ref{tab:Linear-algebraic-groups}. Denote by $a_{1},\dots,a_{l}$
the labels in $\Delta$. Then

\noindent (1) Assume $\Delta$ has no label equal to $1$, then 
\[
\dim\left(\mathfrak{z}\left(\mathfrak{g}^{e}\right)\right)^{G^{e}}=\dim\mathfrak{z}\left(\mathfrak{g}^{e}\right)=n_{2}\left(\Delta\right)=\dim\mathfrak{z}\left(\mathfrak{g}^{h}\right)
\]
 except for the case $\mathfrak{g}=\mathfrak{psl}(n|n)$ and there
exists some $i$ such that $r_{i}\neq s_{i}$, in which case we have
$\dim\left(\mathfrak{z}\left(\mathfrak{g}^{e}\right)\right)^{G^{e}}=\dim\mathfrak{z}\left(\mathfrak{g}^{e}\right)=n_{2}\left(\Delta\right)=\dim\mathfrak{z}\left(\mathfrak{g}^{h}\right)+1$.

\noindent (2) $\dim\left(\mathfrak{z}\left(\mathfrak{g}^{e}\right)\right)^{G^{e}}=\frac{1}{2}\left\lceil \sum_{i=1}^{l}a_{i}\right\rceil +\varepsilon$
where the value of $\varepsilon$ is equal to $0$ except when $\mathfrak{g}=D(2,1;\alpha)$,
$e=E_{1}+E_{2}+E_{3}$ or $\mathfrak{g}=F(4)$, $e=E+e_{(7)}$, in
which cases $\varepsilon=-1$.
\end{thm}

Recall that the sub-labelled Dynkin diagram $\Delta_{0}$ in Theorem
\ref{thm:g^e,z(g^e)} is the 2-free core of $\Delta$ obtained by
removing all nodes with labels equal to $2$ from $\Delta$. 
\begin{thm}
\label{thm:g^e,z(g^e)}Let $\mathfrak{g}=\mathfrak{g}_{\bar{0}}\oplus\mathfrak{g}_{\bar{1}}$
be a basic classical Lie superalgebra and $e\in\mathfrak{g}_{\bar{0}}$
be nilpotent. Let $G$ be the reductive algebraic group defined as
in Table \ref{tab:Linear-algebraic-groups}. Let $\Delta_{0}$ be
the 2-free core of $\Delta$. Denote by $\mathfrak{g}_{0}$ be the
subalgebra of $\mathfrak{g}$ generated by the root vectors corresponding
to the simple roots in $\Delta_{0}$. Let $G_{0}$ be a closed connected
subgroup of $G$ such that $\mathrm{Lie}(G_{0})=(\mathfrak{g}_{0})_{\bar{0}}$.
Denote by $e_{0}$ a nilpotent element in $(\mathfrak{g}_{0})_{\bar{0}}$
whose corresponding labelled Dynkin diagram is $\Delta_{0}$. Then
we have 

(1) $\dim\mathfrak{g}^{e}-\dim\mathfrak{g}_{0}^{e_{0}}=n_{2}(\Delta)$;

(2) $\dim\left(\mathfrak{z}(\mathfrak{g}^{e})\right)^{G^{e}}-\dim\left(\mathfrak{z}(\mathfrak{g}_{0}^{e_{0}})\right)^{G_{0}^{e_{0}}}=n_{2}(\Delta)$
except for some cases when \textup{$\mathfrak{g}=\mathfrak{sl}(m|n)$,
$m\neq n$, $\mathfrak{psl}(n|n)$ and $\mathfrak{osp}(m|2n)$. }

The explicit explanation for exceptions in case of $\mathfrak{psl}(n|n)$
can be found in Subsection \ref{subsec:Proof-of-theorems} and results
in case of $\mathfrak{sl}(m|n)$, \textup{$m\neq n$} and $\mathfrak{osp}(m|2n)$
can be found in \cite{Han-sl-osp}.
\end{thm}

\begin{doublespace}
\noindent %
\begin{longtable}[c]{|c|c|}
\caption{\label{tab:Linear-algebraic-groups}Groups $G$}
\tabularnewline
\endfirsthead
\hline 
Lie superalgebras $\mathfrak{g}$ & Groups $G$\tabularnewline
\hline 
\hline 
$\mathfrak{sl}(m|n),m\neq n$ & $\left\{ (A,B)\in\mathrm{GL}_{m}(\mathbb{C})\times\mathrm{GL}_{n}(\mathbb{C}):\mathrm{det}(A)=\mathrm{det}(B)\right\} $\tabularnewline
\hline 
$\mathfrak{psl}(n|n)$ & $\left\{ (A,B)\in\mathrm{GL}_{n}(\mathbb{C})\times\mathrm{GL}_{n}(\mathbb{C}):\mathrm{det}(A)=\mathrm{det}(B)\}/\left\{ aI_{n|n}:a\in\mathbb{C}^{\times}\right\} \right\} $\tabularnewline
\hline 
$\mathfrak{osp}(m|2n)$ & $\mathrm{O}_{m}(\mathbb{C})\times\mathrm{Sp}_{2n}(\mathbb{C})$\tabularnewline
\hline 
$D(2,1;\alpha)$ & $\mathrm{SL}_{2}(\mathbb{C})\times\mathrm{SL}_{2}(\mathbb{C})\times\mathrm{SL}_{2}(\mathbb{C})$\tabularnewline
\hline 
$G(3)$ & $\mathrm{SL}_{2}(\mathbb{C})\times G_{2}$\tabularnewline
\hline 
$F(4)$ & $\mathrm{SL}_{2}(\mathbb{C})\times\mathrm{Spin}_{7}(\mathbb{C})$\tabularnewline
\hline 
\end{longtable}
\end{doublespace}

This paper is organized as follows. We first recall some fundamental
concepts of Lie superalgebras such as basic classical Lie superalgebras,
root space decomposition and labelled Dynkin diagrams in Section \ref{sec:Preliminaries}.
In Section \ref{sec:Lie-superalgebras-psl}, we recall the structure
of the Lie superalgebra $\mathfrak{psl}(n|n)$ and construct the corresponding
labelled Dynkin diagrams. In Section \ref{sec:centralizer =000026 centre psl},
we determine bases of the centralizer $\mathfrak{g}^{e}$ of a nilpotent
element $e\in\mathfrak{g}_{\bar{0}}$ and centre of centralizer $\mathfrak{z}(\mathfrak{g}^{e})$
for $\mathfrak{g}=\mathfrak{psl}(n|n)$ and proves Theorems \ref{thm:LLD=000026centre}
and \ref{thm:g^e,z(g^e)}. In Sections \ref{sec:sl(m|n)} and \ref{sec:osp(m|2n)},
we classify all reachable elements in $\mathfrak{sl}(m|n)$, $m\neq n$,
$\mathfrak{psl}(n|n)$ or $\mathfrak{osp}(m|2n)$ which proves Theorem
\ref{thm:Re iff partition}. Then we demonstrate that $e\in\mathfrak{g}_{\bar{0}}$
is reachable is equivalent to $e\in\mathfrak{g}_{\bar{0}}$ satisfies
the Panyushev property for the above Lie superalgebras. In Section
\ref{sec:exceptional}, we compute the lists of nilpotent elements
$e\in\mathfrak{g}_{\bar{0}}$ that are reachable, strongly reachable
or satisfy the Panyushev property for $\mathfrak{g}$ an exceptional
Lie superalgebra.

\noindent \textbf{Acknowledgements.} The author acknowledges financial
support from the Engineering and Physical Sciences Research Council
(EP/W522478/1). We would like to thank Simon Goodwin for some
helpful comments.

\section{Preliminaries\label{sec:Preliminaries}}

\subsection{Basic classical Lie superalgebras}

Let $\mathfrak{g}=\mathfrak{g}_{\bar{0}}\oplus\mathfrak{g}_{\bar{1}}$
be a finite-dimensional simple Lie superalgebra over $\mathbb{C}$.
The classification of finite-dimensional complex simple Lie superalgebras
was given by Kac in \cite{Kac}. Among those simple Lie superalgebras,
basic classical Lie superalgebras are of major research interest.
Recall that $\mathfrak{g}$ is called a \textit{basic classical Lie
superalgebra} if $\mathfrak{g}_{\bar{0}}$ is a reductive Lie algebra
and there exists an even non-degenerate invariant supersymmetric bilinear
form $(\cdotp,\cdotp)$ on $\mathfrak{g}$. The basic classical Lie
superalgebras that are not Lie algebras are listed in Table \ref{tab:Basic-classical},
they are $\mathfrak{sl}(m|n)$, $m\neq n$, $\mathfrak{psl}(n|n)=\mathfrak{sl}(m|n)/\mathbb{C}I_{n|n}$,
$n>1$, $\mathfrak{osp}(m|2n)$ and three exceptional types $D(2,1;\alpha)$,
$G(3)$, $F(4)$. 

\begin{table}[H]
\begin{centering}
\begin{tabular}{|c|c|}
\hline 
$\mathfrak{g}$ & $\mathfrak{g}_{\bar{0}}$\tabularnewline
\hline 
\hline 
$A(m,n)=\mathfrak{sl}(m|n)$, $m,n\geq1$, $m\neq n$ & $\mathfrak{sl}(m)\oplus\mathfrak{sl}(n)\oplus\mathbb{C}$\tabularnewline
\hline 
$A(n,n)=\mathfrak{psl}(n|n)$, $n>1$ & $\mathfrak{sl}(n)\oplus\mathfrak{sl}(n)$\tabularnewline
\hline 
$B(m,n)=\mathfrak{osp}(m|2n)$, $m$ is odd, $m,n\geq1$ & $\mathfrak{o}(m)\oplus\mathfrak{sp}(2n)$\tabularnewline
\hline 
$C(n)=\mathfrak{osp}(2|2n)$, $n\geq1$ & $\mathfrak{o}(2)\oplus\mathfrak{sp}(2n)$\tabularnewline
\hline 
$D(m,n)=\mathfrak{osp}(m|2n)$, $m\geq4$ is even, $n\geq1$ & $\mathfrak{o}(m)\oplus\mathfrak{sp}(2n)$\tabularnewline
\hline 
$D(2,1;\alpha)$, $\alpha\neq0,-1$ & $\mathfrak{sl}(2)\oplus\mathfrak{sl}(2)\oplus\mathfrak{sl}(2)$\tabularnewline
\hline 
$G(3)$ & $\mathfrak{sl}(2)\oplus G_{2}$\tabularnewline
\hline 
$F(4)$ & $\mathfrak{sl}(2)\oplus\mathfrak{so}(7)$\tabularnewline
\hline 
\end{tabular}
\par\end{centering}
\caption{\label{tab:Basic-classical}Basic classical Lie superalgebras}
\end{table}

More details about structure of the above Lie superalgebras can be
found for example in \cite[Sections 3--4]{Han-sl-osp} for types $A(m,n)$,
$A(n,n)$, $B(m,n)$, $C(n)$, $D(m,n)$ and \cite[Sections 4--6]{han-exc}
for $D(2,1;\alpha)$, $G(3)$, $F(4)$.

\subsection{Decompositions of $\mathfrak{g}$}

Fix a Cartan subalgebra $\mathfrak{h}\subseteq\mathfrak{g}_{\bar{0}}$
of $\mathfrak{g}$. The following definitions can be found in \cite[Chapter 1]{Cheng2012}.
There exists a root space decomposition $\mathfrak{g}=\mathfrak{h}\oplus\bigoplus_{\alpha\in\mathfrak{h}^{*}}\mathfrak{g}_{\alpha}$
where $\mathfrak{g_{\alpha}}:=\{x\in\mathfrak{g}:[h,x]=\alpha(h)x\ \text{for\ all}\ h\in\mathfrak{\mathfrak{h}}\}$
and $\mathfrak{h}=\mathfrak{g}_{0}$. The set $\Phi=\{\alpha\in\mathfrak{h^{*}:\alpha\neq}0,\mathfrak{g}_{\alpha}\neq0\}$
is called the root system of $\mathfrak{g}$. The sets of even roots
and odd roots are defined to be $\Phi_{\bar{0}}=:\{\alpha\in\Phi:\mathfrak{g}_{\alpha}\subseteq\mathfrak{g}_{\bar{0}}\}$
and $\Phi_{\bar{1}}=:\{\alpha\in\Phi:\mathfrak{g}_{\alpha}\subseteq\mathfrak{g}_{\bar{1}}\}$
respectively. A system of positive roots $\Phi^{+}$ is a subset of
$\Phi$ such that for all roots $\alpha\in\Phi$ exactly one of $\alpha$,$-\alpha$
contained in $\Phi^{+}$ and for any two distinct roots $\alpha,\beta\in\Phi^{+}$,
if $\alpha+\beta\in\Phi$, then $\alpha+\beta\in\Phi^{+}$. The system
of simple roots $\Pi=\{\alpha_{i}:i=1,...,l\}$ determined by $\Phi^{+}$
is a subset of $\Phi^{+}$ such that any $\alpha_{i}\in\Pi$ cannot
be written as sums of positive roots. Note that systems of positive
roots do exist, but they are not all conjugate by Weyl group of $G$.

Fix the system of positive roots $\Phi^{+}$. Let $\mathfrak{b}\subseteq\mathfrak{g}$
be the Borel subalgebra such that $\mathfrak{b}=\mathfrak{h}\oplus\bigoplus_{\alpha\in\Phi^{+}}\mathfrak{g}_{\alpha}$.
Since there are in general many inequivalent conjugacy classes of
Borel subalgebras, there exist more than one system of simple roots
that are determined by $\mathfrak{b}$ up to a transformation of the
Weyl group of $\mathfrak{g}$ in some cases, see for example \cite[Section 2.3]{Frappat}.

\subsection{Labelled Dynkin diagram}

In this subsection, we recall the concept of the labelled Dynkin diagrams
as defined in \cite[Subsection 2.3]{Han-sl-osp}. For $\mathfrak{g}=\mathfrak{g}_{\bar{0}}\oplus\mathfrak{g}_{\bar{1}}$
a basic classical Lie superalgebra, there is a non-degenerate even
invariant supersymmetric bilinear form $(\cdot,\cdot)$ on $\mathfrak{g}$.
Since $(\cdot,\cdot)$ restricts to a non-degenerate symmetric bilinear
form on $\mathfrak{h}$, there exists an isomorphism from $\mathfrak{h}$
to $\mathfrak{h}^{*}$ which provides a symmetric bilinear form on
$\mathfrak{h}^{*}$. The \textit{Dynkin diagram} of $\mathfrak{g}$
with respect to a system of simple roots $\Pi$ is the graph where 

(i) Vertices are labelled by $\Pi$. We say an odd root $\alpha\in\Phi$
is \textit{isotropic} if $(\alpha,\alpha)=0$ and is \textit{non-isotropic}
if $(\alpha,\alpha)\neq0$. Each even root is represented by a white
node $\ocircle$ and each odd isotropic (resp. non-isotropic root)
root is represented by a grey node $\varotimes$ (resp. a black node
$\newmoon$);

(ii) The edges between the vertices correspond to simple roots $\alpha_{i}$
and $\alpha_{j}$ is labelled by $\mu_{i,j}$ where
\begin{equation}
\mu_{i,j}=\begin{cases}
\ensuremath{\vert(\alpha_{i},\alpha_{j})\ensuremath{\vert}} & \text{if }(\alpha_{i},\alpha_{i})=(\alpha_{j},\alpha_{j})=0,\\
\frac{2\ensuremath{\vert}(\alpha_{i},\alpha_{j})\ensuremath{\vert}}{min\{\vert(\alpha_{i},\alpha_{i})\vert,\ensuremath{\vert}(\alpha_{j},\alpha_{j})\ensuremath{\vert}\}} & \text{if }(\alpha_{i},\alpha_{i})(\alpha_{j},\alpha_{j})\neq0,\\
\frac{2\ensuremath{\vert}(\alpha_{i},\alpha_{j})\ensuremath{\vert}}{min_{(\alpha_{k},\alpha_{k})\neq0,\alpha_{k}\in\Phi}\ensuremath{\vert}(\alpha_{k},\alpha_{k})\ensuremath{\vert}} & \text{if }(\alpha_{i},\alpha_{i})\neq0,\ (\alpha_{j},\alpha_{j})=0.
\end{cases}\label{eq:lines-=0003BC}
\end{equation}
 If the value of $\mu_{i,j}$ is a natural number, then we draw $\mu_{i,j}$
edges rather than label the edge.

(iii) When $\mu_{i,j}>1$, we put an arrow pointing from the vertex
labelled by $\alpha_{i}$ to the vertex labelled by $\alpha_{j}$
if $(\alpha_{i},\alpha_{i})(\alpha_{j},\alpha_{j})\neq0$ and $(\alpha_{i},\alpha_{i})>(\alpha_{j},\alpha_{j})$
or if $(\alpha_{i},\alpha_{i})=0,(\alpha_{j},\alpha_{j})\neq0$ and
$\vert(\alpha_{j},\alpha_{j})\vert<2$, or pointing from the vertex
labelled by $\alpha_{j}$ to the vertex labelled by $\alpha_{i}$
if $(\alpha_{i},\alpha_{i})=0,(\alpha_{j},\alpha_{j})\neq0$ and $\vert(\alpha_{j},\alpha_{j})\vert>2$. 

Let $e\in\mathfrak{g}_{\bar{0}}$ be nilpotent. By the Jacobson--Morozov
Theorem \cite[ III, 4.3]{Springer}, there exists an $\mathfrak{sl}(2)$-triple
$\{e,h,f\}$ in $\mathfrak{g}_{\bar{0}}$. Since this $\mathfrak{sl}(2)$-triple
determines a grading on $\mathfrak{g}$ according to the eigenvalues
of ad$h,$ we can decompose $\mathfrak{g}$ into its ad$h$-eigenspaces
such that $\mathfrak{g}=\bigoplus_{j\in\mathbb{Z}}\mathfrak{g}(j)$
and $\mathfrak{g}(j)=\{x\in\mathfrak{g}:[h,x]=jx\}$. 

With the above notation, the \textit{labelled Dynkin diagram} with
respect to $e$ is constructed following the steps below. Note that
in contrast to the case of Lie algebras, in general $e$ determines
more than one labelled Dynkin diagram.

(i) Choose a system $\Phi^{+}(0)$ of positive roots for $\mathfrak{g}(0)$
to get a Borel subalgebra $\mathfrak{b}(0)$ of $\mathfrak{g}(0)$.
Let $\mathfrak{b}$ be the Borel subalgebra such that $\mathfrak{b}=\mathfrak{b}(0)\oplus\mathfrak{g}(j>0)$,
then we obtain the corresponding system of positive roots $\Phi^{+}$
and a system of simple roots $\Pi=\{\alpha_{1},...,\alpha_{l}\}$.

(ii) Construct the the Dynkin diagram of $\mathfrak{g}$ and label
each node $\alpha_{i}$ with $\alpha_{i}(h)$. 
\begin{rem}
For each $\alpha_{i},i=1,...,l$, the corresponding root space $\mathfrak{g}_{\alpha_{i}}\subseteq\mathfrak{g}(j)$
for some $j\geq0$. Hence, we have $\alpha_{i}(h)\geq0$ for $i=1,...,l$.
In fact, $\alpha_{i}(h)\in\{0,1,2\}$ based on our calculation in
\cite{Han-sl-osp} and \cite{han-exc}.
\end{rem}

\section{Lie superalgebras of type $A(m,n)$ and labelled Dynkin diagrams\label{sec:Lie-superalgebras-psl}}

\subsection{Structure of Lie superalgebras of type $A(m,n)$}

We begin with recalling the definition of $\mathfrak{gl}(m|n)$, $\mathfrak{sl}(m|n)$
and $\mathfrak{psl}(n|n)$. Let $V=V_{\bar{0}}\oplus V_{\bar{1}}$
be a finite-dimensional vector superspace such that $\dim V_{\bar{0}}=m$
and $\dim V_{1}=n$. Denote by $M_{a,b}$ the set of $a\times b$
matrices. As a vector space $\mathfrak{gl}(m|n)=\mathfrak{gl}(V)=\mathfrak{g}_{\bar{0}}\oplus\mathfrak{g}_{\bar{1}}$
is $M_{m+n,m+n}$ where 
\[
\mathfrak{g}_{\bar{0}}=\left\{ \begin{pmatrix}A & 0\\
0 & D
\end{pmatrix}:A\in M_{m,m},D\in M_{n,n}\right\} ,\text{ and }\mathfrak{g}_{\bar{1}}=\left\{ \begin{pmatrix}0 & B\\
C & 0
\end{pmatrix}:B\in M_{m,n},C\in M_{n,m}\right\} .
\]
 For an element $x\in\mathfrak{g}_{\bar{i}}$, we define the \textit{parity}
of $x$ by $\bar{x}=\bar{i}$. The Lie bracket is defined by 
\[
[x,y]=xy-(-1)^{\bar{x}\bar{y}}yx\text{ for all homogeneous }x,y\in\mathfrak{gl}(m|n).
\]
The Lie superalgebra $\mathfrak{sl}(m|n)$ is defined to be 
\[
\mathfrak{sl}(m|n)=\left\{ x=\begin{pmatrix}A & B\\
C & D
\end{pmatrix}\in\mathfrak{gl}(m|n):\mathrm{str}(x)=\mathrm{tr}(A)-\mathrm{tr}(D)=0\right\} .
\]
 We know that $\mathfrak{sl}(m|n)$ is simple when $m\neq n$, $m,n>0$.
When $m=n$, the identity matrix $I_{n|n}\in\mathfrak{sl}(m|n)$ is
central and we denote $\mathfrak{psl}(n|n)=\mathfrak{sl}(n|n)/\mathbb{C}I_{n|n}$.
Note that $\mathfrak{sl}(1|1)/\mathbb{C}I_{n|n}$ is abelian thus
is not simple. For $n\geq2$, we have that $\mathfrak{psl}(n|n)$
is simple, see for example \cite[Subsection 1.1.2]{Cheng2012}.

\subsection{Dynkin pyramids and labelled Dynkin diagrams\label{subsec:Dynkin-pyramids}}

Let $\mathfrak{g}=\mathfrak{g}_{\bar{0}}\oplus\mathfrak{g}_{\bar{1}}=\mathfrak{gl}(m|n),\mathfrak{sl}(m|n)$
or $\mathfrak{psl}(n|n)$. Recall that the nilpotent orbits in $\mathfrak{g}_{\bar{0}}$
are parametrized by the partitions of $(m|n)$. Let $\lambda$ be
a partition of $(m|n)$ such that 
\begin{equation}
\lambda=(p\mid q)=(p_{1},\dots,p_{r}\mid q_{1},\dots,q_{s})\label{eq:partition(p,q)}
\end{equation}
 where $p$ and $q$ are partitions of $m$ and $n$ respectively
and $p_{1}\geq\dots\geq p_{r},q_{1}\geq\dots\geq q_{s}$. In this
paper, we also use the following notation for $\lambda$. Write 
\begin{equation}
\lambda=(\lambda_{1},\dots,\lambda_{r+s})\label{eq:partition 1,2,3...}
\end{equation}
 by rearranging the order of numbers in $(p|q)$ in (\ref{eq:partition(p,q)})
such that $\lambda_{1}\geq\dots\geq\lambda_{r+s}$. For $i=1,\dots,r+s$,
we define $\bar{i}\in\{\bar{0},\bar{1}\}$ such that for $c\in\mathbb{Z}$,
we have that $|\{i:\lambda_{i}=c,\bar{i}=\bar{0}\}|=|\{j:p_{j}=c\}|$,
$|\{i:\lambda_{i}=c,\bar{i}=\bar{1}\}|=|\{j:q_{j}=c\}|$ and if $\lambda_{i}=\lambda_{j}$,
$\bar{i}=\bar{0},\bar{j}=\bar{1}$, then $i<j$. In particular, $\sum_{\ensuremath{\bar{i}}=\bar{0}}\lambda_{i}=m$,
$\sum_{\bar{i}=\bar{1}}\lambda_{i}=n$.

In order to determine the labelled Dynkin diagrams for $\mathfrak{g}$,
we recall the definition of the \textit{Dynkin pyramid of shape} $\lambda$,
following \cite{Elashvili=000026Kac} and \cite{Hoyt2012}. This is
a diagram consisting of a finite collection of boxes of size $2\times2$
drawn in the upper half of the $xy$-plane. The coordinates of a box
is represented by the coordinate of its midpoint. Write $\lambda=(\lambda_{1},\dots,\lambda_{r+s})$
as defined in (\ref{eq:partition 1,2,3...}), we number rows of the
Dynkin pyramid from $1,\dots,r+s$ from bottom to top such that the
$i$th row has length $\lambda_{i}$ and boxes in the row $i$ has
parity $\bar{i}$. We fix a numbering $1,\dots,m+n$ of the boxes
of the Dynkin pyramid from top to bottom and from left to right. We
also denote by $\mathrm{row}(i)$ (resp. $\mathrm{col}(i)$) the $x$-coordinate
(resp. $y$-coordinate) of the $i$th box and the parity of $\mathrm{row}(i)$
is defined to be $|\mathrm{row}(i)|$.

According to \cite[Section 7]{Hoyt2012}, the Dynkin pyramid determines
a nilpotent element $e\in\mathfrak{g}_{\bar{0}}$ and a semisimple
element $h\in\mathfrak{g}_{\bar{0}}$ such that $\{e,h\}$ can be
extended to an $\mathfrak{sl}(2)$-triple in $\mathfrak{g}_{\bar{0}}$.
Write $e_{ij}$ for the $ij$-matrix unit, then the nilpotent element
\[
e=\sum_{\substack{\mathrm{row}(i)=\mathrm{row}(j)\\
\mathrm{col}(j)=\mathrm{col}(i)-2
}
}e_{ij}
\]
 for all $i,j\in\{1,\dots,m+n\}$ and the semisimple element $h=\sum_{i=1}^{m+n}-\mathrm{col}(i)e_{ii}$. 

Let $\mathfrak{h}$ be the Cartan subalgebra of $\mathfrak{g}$ consisting
of all diagonal matrices in $\mathfrak{g}$. We define $\{\varepsilon_{1},\dots,\varepsilon_{m+n}\}\subseteq\mathfrak{h}^{*}$
such that for $a=\mathrm{diag}(a_{1},\dots,a_{m+n})\in\mathfrak{h}$,
we have $\varepsilon_{i}(a)=a_{i}$ for $i=1,\dots,m+n$ and the parity
of $\varepsilon_{i}$ is equal to $\ensuremath{\left|\mathrm{row}(i)\right|}$.
According to \cite[Section 3]{Dimitrov}, the root system of $\mathfrak{g}$
with respect to $\mathfrak{h}$ is $\Phi=\Phi_{\bar{0}}\cup\Phi_{\bar{1}}$
where 
\[
\Phi_{\bar{0}}=\{\varepsilon_{i}-\varepsilon_{j}:i\neq j,\ensuremath{\bar{i}}=\bar{j}\},\Phi_{\bar{1}}=\{\varepsilon_{i}-\varepsilon_{j}:i\neq j,\ensuremath{\bar{i}}\neq\ensuremath{\bar{j}}\}
\]
 and $(\varepsilon_{i},\varepsilon_{j})=(-1)^{\ensuremath{\bar{i}}}\delta_{ij}$.
For $i\neq j,\bar{i}\neq\bar{j}$, all odd roots are isotropic because
$(\varepsilon_{i}-\varepsilon_{j},\varepsilon_{i}-\varepsilon_{j})=0$.

Fix the Dynkin pyramid $P$ of shape $\lambda$. Next we use $P$
to determine the labelled Dynkin diagram $\Delta$ for $e$. For the
box labelled by $i$ such that $i=1,\dots,m+n-1$, we associate a
white node $\ocircle$ (resp. a grey node $\varotimes$) for root
$\alpha_{i}=\varepsilon_{i}-\varepsilon_{i+1}$ if $\ensuremath{\left|\mathrm{row}(i+1)\right|}=\ensuremath{\left|\mathrm{row}(i)\right|}$
(resp. $\ensuremath{\left|\mathrm{row}(i+1)\right|}\neq\ensuremath{\left|\mathrm{row}(i)\right|}$).
The $i$th node is labelled by the value $\mathrm{col}(i+1)-\mathrm{col}(i)$. 

\section{The centralizer of a nilpotent element in $\mathfrak{psl}(n|n)$
and its centre\label{sec:centralizer =000026 centre psl}}

\subsection{\label{subsec:Centralizer-of-nilpoent}Centralizer of nilpotent elements
$e$ in $\mathfrak{psl}(n|n)_{\bar{0}}$}

Let $\mathfrak{g}=\mathfrak{g}_{\bar{0}}\oplus\mathfrak{g}_{\bar{1}}=\mathfrak{psl}(n|n)$,
$n>1$ and let $e\in\mathfrak{g}_{\bar{0}}$ be nilpotent. Define $\phi:\mathfrak{\mathfrak{sl}}(n|n)\rightarrow\mathfrak{g}$
to be the quotient map. Let $\underline{e}\in\mathfrak{sl}(n|n)_{\bar{0}}$
be the lift of $e$ in $\mathfrak{sl}(n|n)$, i.e. $\phi\left(\underline{e}\right)=e$.
In order to compute the dimension of $\mathfrak{g}^{e}$, we first
recall bases for $\mathfrak{gl}(n|n)^{\underline{e}}$ and $\mathfrak{sl}(n|n)^{\underline{e}}$
based on \cite[Section 3.2]{Hoyt2012} and \cite[Section 3.2]{Han-sl-osp}. 

We know that $\mathfrak{gl}(n|n)=\mathrm{End}(V_{\bar{0}}\oplus V_{\bar{1}})$
where $\mathfrak{gl}(n|n)_{\bar{0}}=\mathrm{End}(V_{\bar{0}})\oplus\mathrm{End}(V_{\bar{1}})$
and $\mathfrak{gl}(n|n)_{\bar{1}}=\mathrm{Hom}(V_{\bar{0}},V_{\bar{1}})\oplus\mathrm{Hom}(V_{\bar{1}},V_{\bar{0}})$.
There exist $v_{1},\dots,v_{r+s}\in V$ such that $\{\underline{e}^{j}v_{i}:0\leq j\leq\lambda_{i}-1,\bar{i}=\bar{0}\}$
is a basis for $V_{\bar{0}}$ , $\{\underline{e}^{j}v_{i}:0\leq j\leq\lambda_{i}-1,\bar{i}=\bar{1}\}$
is a basis for $V_{\bar{1}}$ and $\underline{e}^{\lambda_{i}}v_{i}=0$
for $1\leq i\leq r+s$. Each element $\xi\in\mathfrak{gl}(n|n)^{\underline{e}}$
is determined by $\xi(v_{i})$. According to arguments in \cite[Section 3.2]{Hoyt2012}
and \cite[Section 3.2]{Han-sl-osp}, a basis for $\mathfrak{gl}(n|n)^{\underline{e}}$
is of the form
\begin{equation}
\left\{ \xi_{i}^{j,k}:1\leq i,j\leq r+s\text{ and }\max\{\lambda_{j}-\lambda_{i},0\}\leq k\leq\lambda_{j}-1\right\} \label{eq:basis gl}
\end{equation}
 where $\xi_{i}^{j,k}$ sends $v_{i}$ to $\underline{e}^{k}v_{j}$
and the other basis elements to zero. Hence, the dimension of $\mathfrak{gl}(n|n)^{\underline{e}}$
is determined by the number of indices $i,j,k$ and we have 
\begin{equation}
\text{\ensuremath{\dim}}\mathfrak{gl}(n|n)_{\bar{0}}^{\underline{e}}=\left(n+2\sum_{i=1}^{r}(i-1)p_{i}\right)+\left(n+2\sum_{j=1}^{s}(j-1)q_{j}\right),\label{eq:dim(gl0^e)}
\end{equation}
\begin{equation}
\dim\mathfrak{gl}(n|n)_{\bar{1}}^{\underline{e}}=2\sum_{i,j=1}^{r,s}\min(p_{i},q_{j})\label{eq:dim(gl1^e)}
\end{equation}
where $p_{i}$ and $q_{j}$ are defined in (\ref{eq:partition(p,q)}).

Next we consider a basis for $\mathfrak{sl}(n|n)^{\underline{e}}$.
Note that the author proved that $\dim\mathfrak{sl}(n|n)^{\underline{e}}=\dim\mathfrak{gl}(n|n)^{\underline{e}}-1$
in \cite[Theorem 10]{Han-sl-osp}. Since the diagonal matrices $\xi_{i}^{i,0}$
are not in $\mathfrak{sl}(n|n)^{\underline{e}}$, a basis for $\mathfrak{sl}(n|n)^{\underline{e}}$
can be written as 
\[
\left\{ \xi_{i}^{j,k}:1\leq i\neq j\leq r+s\text{ and }\max\{\lambda_{j}-\lambda_{i},0\}\leq k\leq\lambda_{j}-1\right\} 
\]
\[
\cup\left\{ \xi_{i}^{i,k}:0<k\leq\lambda_{i}-1\right\} \cup\left\{ \lambda_{i+1}\xi_{i}^{i,0}-(-1)^{\bar{i}}\lambda_{i}\xi_{i+1}^{i+1,0}:1\leq i<r+s\right\} .
\]

Now we are ready to calculate the dimension of $\mathfrak{g}^{e}$
and further give a basis for $\mathfrak{g}^{e}$.
\begin{lem}
\label{lem:g^e}Let $e\in\mathfrak{g}_{\bar{0}}$ be nilpotent, we
have that $\dim\mathfrak{g}^{e}=\dim\mathfrak{\mathfrak{sl}}(n|n)^{\underline{e}}-1$.
\end{lem}

\begin{proof}
Suppose $\phi(\underline{x})\in\mathfrak{g}^{e}$, then $[\underline{x},\underline{e}]=0$
or $aI_{n\mid n}$ for some $a\in\mathbb{C}^{\times}$. Assume $[\underline{x},\underline{e}]=aI_{n\mid n}\in\mathfrak{\mathfrak{sl}}(n|n)_{\bar{0}}$
for some $a\in\mathbb{C}^{\times}$, we can write $\underline{x}=\underline{x}_{\bar{0}}+\underline{x}_{\bar{1}}$
such that $\underline{x}_{\bar{0}}\in\mathfrak{\mathfrak{sl}}(n|n)_{\bar{0}}$
and $\underline{x}_{\bar{1}}\in\mathfrak{\mathfrak{sl}}(n|n)_{\bar{1}}$.
Then we have that $[\underline{x},\underline{e}]=[\underline{x}_{\bar{0}},\underline{e}]+[\underline{x}_{\bar{1}},\underline{e}]$,
thus $[\underline{x}_{\bar{0}},\underline{e}]$ must be equal to $aI_{n\vert n}$.
However, this is impossible because $\underline{e}\in\mathfrak{sl}(n)\oplus\mathfrak{sl}(n)$
so that $[\underline{x}_{\bar{0}},\underline{e}]\in\mathfrak{sl}(n)\oplus\mathfrak{sl}(n)$.
Hence, we deduce that $\phi^{e}:\mathfrak{\mathfrak{sl}}(n|n)^{\underline{e}}\rightarrow\mathfrak{g}^{e}$
is surjective and $\ker\phi^{e}=\langle I_{n\vert n}\rangle$. Therefore,
we have that $\dim\mathfrak{g}^{e}=\dim\mathfrak{\mathfrak{sl}}(n|n)^{\underline{e}}-1$.
\end{proof}
Since the identity matrix $I_{n|n}\in\mathfrak{\mathfrak{sl}}(n|n)^{\underline{e}}$
and $\phi$ sends $I_{n|n}$ to zero, we can remove one element from
$\lbrace\lambda_{i+1}\xi_{i}^{i,0}-(-1)^{\bar{i}}\lambda_{i}\xi_{i+1}^{i+1,0}:1\leq i<r+s\rbrace$
in the basis for $\mathfrak{sl}(n|n)^{\underline{e}}$ to obtain a
basis for $\mathfrak{g}^{e}$. Hence, we write a basis for $\mathfrak{g}^{e}$ as
\[
\left\{ \xi_{i}^{j,k}:1\leq i\neq j\leq r+s\text{ and }\max\{\lambda_{j}-\lambda_{i},0\}\leq k\leq\lambda_{j}-1\right\} 
\]
\[
\cup\left\{ \xi_{i}^{i,k}:1\leq i\leq r+s,0<k\leq\lambda_{i}-1\right\} \cup\left\{ \lambda_{i+1}\xi_{i}^{i,0}-(-1)^{\bar{i}}\lambda_{i}\xi_{i+1}^{i+1,0}:1<i<r+s\right\} .
\]

In order to prove the theorems in Section \ref{subsec:Proof-of-theorems},
we use an alternative formula for $\dim\mathfrak{g}^{e}$ based on
\cite[Proposition 9]{Han-sl-osp}.
\begin{prop}
\label{prop:dimg^e~c_i}Let $\lambda$ be a partition of $(n|n)$
defined as in (\ref{eq:partition 1,2,3...}). Denote by $P$ the Dynkin
pyramid of shape $\lambda$ and $e\in\mathfrak{g}_{\bar{0}}$ the
nilpotent element determined by $P$. Let $c_{i}$ be the number of
boxes in the $i$th column of $P$. Then $\dim\mathfrak{g}^{e}=\sum_{i\in\mathbb{Z}}c_{i}^{2}+\sum_{i\in\mathbb{Z}}c_{i}c_{i+1}-2$.
\end{prop}

\begin{proof}
The proof of \cite[Proposition 9]{Han-sl-osp} implies that $\dim\mathfrak{gl}(n|n)^{\underline{e}}=\sum_{i\in\mathbb{Z}}c_{i}^{2}+\sum_{i\in\mathbb{Z}}c_{i}c_{i+1}$.
We have already shown that $\dim\mathfrak{g}^{e}=\dim\mathfrak{\mathfrak{sl}}(n|n)^{\underline{e}}-1=\dim\mathfrak{gl}(n|n)^{\underline{e}}-2$,
as required.
\end{proof}

\subsection{Centre of centralizer of nilpotent elements $e$ in $\mathfrak{psl}(n|n)_{\bar{0}}$}

The following theorem provides a basis for $\mathfrak{z}(\mathfrak{g}^{e})$.
\begin{thm}
\label{thm:z(g^e)}Let $\mathfrak{g}=\mathfrak{g}_{\bar{0}}\oplus\mathfrak{g}_{\bar{1}}=\mathfrak{psl}(n|n)$
and let $\lambda$ be a partition of $(n|n)$ denoted as (\ref{eq:partition 1,2,3...}).
Let $e\in\mathfrak{g}_{\bar{0}}$ be a nilpotent element with Jordan
type $\lambda$. We have that $\mathfrak{z}(\mathfrak{g}^{e})=\mathrm{Span}\{e,\dots,e^{\lambda_{1}-1}\}$.
\end{thm}

\begin{proof}
For the partition $\lambda=\left(\lambda_{1},\dots,\lambda_{r+s}\right)$,
we know that elements $\lambda_{h}\xi_{t}^{t,0}-(-1)^{\bar{t}}\lambda_{t}\xi_{h}^{h,0}$
with $t\neq h$ lie in the Cartan subalgebra $\mathfrak{h}\subseteq\mathfrak{g}$. 

Suppose $x\in\mathfrak{z}(\mathfrak{g}^{e})$ is of the form 
\[
x=\sum_{\substack{1\leq i\neq j\leq r+s\\
\max\{\lambda_{j}-\lambda_{i},0\}\leq k\leq\lambda_{j}-1
}
}a_{i}^{j,k}\xi_{i}^{j,k}+\sum_{\substack{0<k\leq\lambda_{i}-1\\
1\leq i\leq r+s
}
}b_{i}^{i,k}\xi_{i}^{i,k}+\sum_{1<i<r+s}c_{i}^{i,0}\left(\lambda_{i+1}\xi_{i}^{i,0}-(-1)^{\bar{i}}\lambda_{i}\xi_{i+1}^{i+1,0}\right)
\]
 where $a_{i}^{j,k},b_{i}^{i,k},c_{i}^{i,0}\in\mathbb{C}$ are coefficients.
If $r+s\geq3$, then $x$ commutes with elements $\lambda_{h}\xi_{t}^{t,0}-(-1)^{\bar{t}}\lambda_{t}\xi_{h}^{h,0}$
for all $t\neq h$. It is clear that elements $\xi_{i}^{i,k}$, $\lambda_{i+1}\xi_{i}^{i,0}-(-1)^{\bar{i}}\lambda_{i}\xi_{i+1}^{i+1,0}\in\mathfrak{g}_{\bar{0}}$
commutes with $\lambda_{h}\xi_{t}^{t,0}-(-1)^{\bar{t}}\lambda_{t}\xi_{h}^{h,0}$
for all $t\neq h$. Hence, we have 
\begin{align}
[x,\lambda_{h}\xi_{t}^{t,0}-(-1)^{\bar{t}}\lambda_{t}\xi_{h}^{h,0}] & =\lambda_{h}\sum_{t,j,k}a_{t}^{j,k}\xi_{t}^{j,k}-\lambda_{h}\sum_{i,t,k}a_{i}^{t,k}\xi_{i}^{t,k}\nonumber \\
 & -(-1)^{\bar{t}}\lambda_{t}\sum_{h,j,k}a_{h}^{j,k}\xi_{h}^{j,k}+(-1)^{\bar{t}}\lambda_{t}\sum_{i,h,k}a_{i}^{h,k}\xi_{i}^{h,k}.\label{eq:centre psl-1}
\end{align}
We know that $\mathfrak{z}\left(\mathfrak{sl}\left(n|n\right)\right)$
can be expressed as $\left\langle \sum_{1\leq i\leq r+s}\xi_{i}^{i,0}\right\rangle $.
Clearly summands in (\ref{eq:centre psl-1}) can not be equal to $\mathfrak{z}\left(\mathfrak{sl}\left(n|n\right)\right)$,
thus $[x,\lambda_{h}\xi_{t}^{t,0}-(-1)^{\bar{t}}\lambda_{t}\xi_{h}^{h,0}]$
must be zero. This implies that $a_{i}^{j,k}=0$ for all $1\leq i\neq j\leq r+s$.
Hence, an element $x\in\mathfrak{z}(\mathfrak{g}^{e})$ is of the
form 
\[
x=\sum_{\substack{0<k\leq\lambda_{i}-1\\
1\leq i\leq r+s
}
}b_{i}^{i,k}\xi_{i}^{i,k}+\sum_{1<i<r+s}c_{i}^{i,0}\left(\lambda_{i+1}\xi_{i}^{i,0}-(-1)^{\bar{i}}\lambda_{i}\xi_{i+1}^{i+1,0}\right).
\]

Next we consider an element $\xi_{t}^{h,0}\in\mathfrak{g}^{e}$ for
some $t<h$. By computing $[x,\xi_{t}^{h,0}]$ we have 
\begin{align*}
[x,\xi_{t}^{h,0}] & =\left(\sum_{h,k}b_{h}^{h,k}-\sum_{t,k}b_{t}^{t,k}\right)\xi_{t}^{h,k}+\xi_{t}^{h,0}(\lambda_{h+1}c_{h}^{h,0}-\lambda_{t+1}c_{t}^{t,0}\\
 & -(-1)^{\overline{h-1}}\lambda_{h-1}c_{h-1}^{h-1,0}+(-1)^{\overline{t-1}}\lambda_{t-1}c_{t-1}^{t-1,0}).
\end{align*}
Similarly this can not be equal to $\mathfrak{z}\left(\mathfrak{sl}\left(n|n\right)\right)$,
thus $[x,\xi_{t}^{h,0}]$ must be equal to zero. Hence, we have that
$c_{i}^{i,0}=0$ for $1<i<r+s$ and $b_{h}^{h,k}=b_{t}^{t,k}$ for
all $1\leq t\neq h\leq r+s$ and $0<k\leq\lambda_{1}-1$. We know
that $e^{k}=\sum_{i=1}^{r+s}\xi_{i}^{i,k}$. Therefore, we deduce
that $x\in\mathfrak{z}(\mathfrak{g}^{e})$ if and only if  $x\in\langle e,\dots,e^{\lambda_{1}-1}\rangle$.

If $r+s=2$, we have $\lambda=(\lambda_{1}|\lambda_{2})=(n|n)$ for
$n>1$. In this case, an element $y\in\mathfrak{g}^{e}$ is of the
form 
\[
y=\sum_{\substack{1\leq i\neq j\leq2\\
0\leq k\leq n-1
}
}a_{i}^{j,k}\xi_{i}^{j,k}+\sum_{\substack{0<k\leq n-1\\
1\leq i\leq2
}
}b_{i}^{i,k}\xi_{i}^{i,k}
\]
 where $a_{i}^{j,k},b_{i}^{i,k}\in\mathbb{C}$ are coefficients. Now
suppose $y\in\mathfrak{z}(\mathfrak{g}^{e})$, we consider an element
$\xi_{1}^{1,1}\in\mathfrak{g}^{e}$. By calculating 
\[
[y,\xi_{1}^{1,1}]=\sum_{k=0}^{n-1}a_{2}^{1,k}\xi_{2}^{1,k+1}-\sum_{k=0}^{n-1}a_{1}^{2,k}\xi_{1}^{2,k+1}.
\]
This cannot be equal to $\mathfrak{z}\left(\mathfrak{sl}\left(n|n\right)\right)$
implies that $a_{1}^{2,k}=a_{2}^{1,k}=0$ for $0\leq k\leq n-2$.
Then we calculate 
\[
[y,\xi_{2}^{1,0}]=a_{1}^{2,n-1}(\xi_{1}^{1,n-1}+\xi_{2}^{2,n-1})+(\sum_{k=1}^{n-1}b_{1}^{1,k}-\sum_{k=1}^{n-1}b_{2}^{2,k})\xi_{2}^{1,k}.
\]
 Since the above summands cannot be written as $\left\langle \xi_{1}^{1,0}+\xi_{2}^{2,0}\right\rangle =\mathfrak{z}\left(\mathfrak{sl}\left(n|n\right)\right)$,
we deduce that $[y,\xi_{2}^{1,0}]=0$. This implies that $a_{1}^{2,n-1}=0$
and $b_{1}^{1,k}=b_{2}^{2,k}$ for $0<k\leq n-1$. Similarly, we have
that $a_{2}^{1,n-1}=0$ by considering $\xi_{1}^{2,0}\in\mathfrak{g}^{e}$.
Therefore, we obtain that $\mathfrak{z}(\mathfrak{g}^{e})=\mathrm{Span}\{e,\dots,e^{\lambda_{1}-1}\}$.
\end{proof}
\begin{thm}
\label{thm:group action}Let $\mathfrak{g}=\mathfrak{g}_{\bar{0}}\oplus\mathfrak{g}_{\bar{1}}=\mathfrak{psl}(n|n)$
and $e\in\mathfrak{g}_{\bar{0}}$ be nilpotent. Let $G$ be the reductive
algebraic group such that $G=\{(A,B)\in\mathrm{GL}_{n}(\mathbb{C})\times\mathrm{GL}_{n}(\mathbb{C}):\mathrm{det}(A)=\mathrm{det}(B)\}/\{aI_{n|n}:a\in\mathbb{C}^{\times}\}$.
Then $\left(\mathfrak{z}(\mathfrak{g}^{e})\right)^{G^{e}}=\mathfrak{z}(\mathfrak{g}^{e})$. 
\end{thm}

\begin{proof}
Let $\bar{\mathfrak{g}}=\bar{\mathfrak{g}}_{\bar{0}}\oplus\bar{\mathfrak{g}}_{\bar{1}}=\mathfrak{gl}(n|n)/\mathbb{C}I_{n|n}$
and $\bar{G}=\{(A,B)\in\mathrm{GL}_{n}(\mathbb{C})\times\mathrm{GL}_{n}(\mathbb{C})\}/\{aI_{n|n}:a\in\mathbb{C}^{\times}\}$. 

We know that $\bar{G}=GZ(\bar{G})=\{gg':g\in G,g'\in Z(\bar{G})\}$
and thus $\bar{G}^{e}=G^{e}Z(\bar{G})$. Since $Z(\bar{G})$ centralizes
$\mathfrak{g}_{\bar{0}}$ and we have already shown that $\mathfrak{z}(\mathfrak{g}^{e})\subseteq\mathfrak{g}_{\bar{0}}$
in Theorem \ref{thm:z(g^e)}, we deduce that $\left(\mathfrak{z}(\mathfrak{g}^{e})\right)^{G^{e}}=\left(\mathfrak{z}(\mathfrak{g}^{e})\right)^{\bar{G}^{e}}$.
By definition $\left(\mathfrak{z}(\mathfrak{g}^{e})\right)^{\bar{G}^{e}}=\{x\in\mathfrak{z}(\mathfrak{g}^{e}):g\cdotp x=x\text{ for all }g\in\bar{G}^{e}\}$,
this implies that $x\in\left(\mathfrak{z}(\mathfrak{g}^{e})\right)^{\bar{G}^{e}}$
if and only if $\left(\bar{G}^{e}\right)^{x}=\bar{G}^{e}$. By \cite[Theorem 13.2]{Humphreys1975},
we deduce that $\mathrm{Lie}\left(\left(\bar{G}^{e}\right)^{x}\right)=\left(\bar{\mathfrak{g}}_{\bar{0}}^{e}\right)^{x}$
for any $x\in\mathfrak{z}(\mathfrak{g}^{e})$. We also know that $\mathrm{Lie}\left(\left(\bar{G}^{e}\right)^{x}\right)=\bar{\mathfrak{g}}_{\bar{0}}^{e}$
if and only if $\left(\bar{G}^{e}\right)^{x}=\bar{G}^{e}$. Moreover,
by definition $\left(\mathfrak{z}(\mathfrak{g}^{e})\right)^{\bar{\mathfrak{g}}_{\bar{0}}^{e}}=\{x\in\mathfrak{z}(\mathfrak{g}^{e}):[x,y]=0\text{ for all }y\in\bar{\mathfrak{g}}_{\bar{0}}^{e}\}$
and thus $x\in\left(\mathfrak{z}(\mathfrak{g}^{e})\right)^{\bar{\mathfrak{g}}_{\bar{0}}^{e}}$
if and only if $\left(\bar{\mathfrak{g}}_{\bar{0}}^{e}\right)^{x}=\bar{\mathfrak{g}}_{\bar{0}}^{e}$.
Hence, we obtain that $x\in\left(\mathfrak{z}(\mathfrak{g}^{e})\right)^{\bar{G}^{e}}$
if and only if $x\in\left(\mathfrak{z}(\mathfrak{g}^{e})\right)^{\bar{\mathfrak{g}}_{\bar{0}}^{e}}$
and thus $\left(\mathfrak{z}(\mathfrak{g}^{e})\right)^{\bar{G}^{e}}=\left(\mathfrak{z}(\mathfrak{g}^{e})\right)^{\bar{\mathfrak{g}}_{\bar{0}}^{e}}$.
Moreover, we know that $\left(\mathfrak{z}(\mathfrak{g}^{e})\right)^{\bar{\mathfrak{g}}_{\bar{0}}^{e}}=\left(\mathfrak{z}(\mathfrak{g}^{e})\right)^{\mathfrak{g}_{\bar{0}}^{e}}=\mathfrak{z}(\mathfrak{g}^{e})$
because $\bar{\mathfrak{g}}_{\bar{0}}^{e}=\mathfrak{g}_{\bar{0}}^{e}+\mathfrak{z}(\bar{\mathfrak{g}}_{\bar{0}})$
and $\mathfrak{z}(\bar{\mathfrak{g}}_{\bar{0}})$ centralizes all
of $\mathfrak{g}_{\bar{0}}^{e}$. Therefore, we deduce that $\left(\mathfrak{z}(\mathfrak{g}^{e})\right)^{G^{e}}=\mathfrak{z}(\mathfrak{g}^{e})$.
\end{proof}

\subsection{\label{subsec:Proof-of-theorems}Relation between centralizers and
labelled Dynkin diagrams}

In this subsection, we prove Theorems \ref{thm:LLD=000026centre}
and \ref{thm:g^e,z(g^e)} for $\mathfrak{g}=\mathfrak{g}_{\bar{0}}\oplus\mathfrak{g}_{\bar{1}}=\mathfrak{psl}(n|n)$,
$n>1$. Let $\lambda$ be a partition of $(n|n)$ as in (\ref{eq:partition 1,2,3...}).
Construct a Dynkin pyramid $P$ of shape $\lambda$ following Section
\ref{subsec:Dynkin-pyramids}. Denote by $r_{i}$ (resp. $s_{i}$)
the number of boxes with parity $\bar{0}$ (resp. $\bar{1}$) in the
$i$th column of $P$. Let $e\in\mathfrak{g}_{\bar{0}}$ be the nilpotent
element determined by $P$. 

\noindent \textit{Proof of Theorem \ref{thm:LLD=000026centre}.(2)
for $\mathfrak{g}=\mathfrak{psl}(n|n)$.} Based on the way that the
labelled Dynkin diagram $\Delta$ is constructed, any $\Delta$ with
respect to $e\in\mathfrak{g}_{\bar{0}}$ is the same as the labelled
Dynkin diagram for $\underline{e}\in\mathfrak{sl}(n|n)$. According
to \cite[Theorem 2]{Han-sl-osp}, we have that $\dim\mathfrak{z}\left(\mathfrak{sl}\left(n|n\right)^{\underline{e}}\right)=\frac{1}{2}\sum_{i=1}^{l}a_{i}$.
By Theorems \ref{thm:z(g^e)}--\ref{thm:group action} and \cite[Theorem 13]{Han-sl-osp},
we know that $\dim\left(\mathfrak{z}(\mathfrak{g}^{e})\right)^{G^{e}}=\dim\mathfrak{z}\left(\mathfrak{sl}\left(n|n\right)^{\underline{e}}\right)$.
Therefore, we deduce that $\dim\left(\mathfrak{z}(\mathfrak{g}^{e})\right)^{G^{e}}=\frac{1}{2}\sum_{i=1}^{l}a_{i}$.
\qed

Next we look at the case when $\Delta$ has no label equal to $1$. 

\noindent \textit{Proof of Theorem \ref{thm:LLD=000026centre}.(1)
for $\mathfrak{g}=\mathfrak{psl}(n|n)$.} We know that the labels
in the labelled Dynkin diagram $\Delta$ are the horizontal difference
between consecutive boxes in the Dynkin pyramid. Thus $\Delta$ has
no label equal to $1$ if and only if $\lambda_{i}$ are all even
or $\lambda_{i}$ are all odd for $i=1,\dots,r+s$. Based on the way
that $\Delta$ is constructed, we have that $n_{2}\left(\Delta\right)=\lambda_{1}-1$.
Hence, we have that $\dim\left(\mathfrak{z}(\mathfrak{g}^{e})\right)^{G^{e}}=n_{2}\left(\Delta\right)$
by Theorems \ref{thm:z(g^e)} and \ref{thm:group action}. 

Next we look at $\mathfrak{g}^{h}$. Recall that for any element $x\in\mathfrak{g}$,
we denote by $\underline{x}$ the lift of $x$ in $\mathfrak{\mathfrak{sl}}(n|n)$
such that $\phi(\underline{x})=x$. Using a similar argument as in
the proof of Lemma \ref{lem:g^e}, we obtain that $\mathfrak{\mathfrak{sl}}(n|n)^{\underline{h}}\rightarrow\mathfrak{g}^{h}$
is surjective. Thus an element in $\mathfrak{g}^{h}$ is of the form
\[
\begin{pmatrix}x_{-\lambda_{1}+1} & \dotsc & 0\\
\vdots & \ddots & \vdots\\
0 & \dotsb & x_{\lambda_{1}-1}
\end{pmatrix}+\mathfrak{z}\left(\mathfrak{\mathfrak{sl}}(n|n)\right)
\]
 with respect to the basis determined by the Dynkin pyramid such that
$x_{i}\in\mathfrak{gl}(r_{i}|s_{i})$ for $i=-\lambda_{1}+1,-\lambda_{1}+3,\dots,\lambda_{1}-3,\lambda_{1}-1$
and $\sum_{i}\mathrm{str}(x_{i})=0$. 

To determine the structure of $\mathfrak{z}(\mathfrak{g}^{h})$, we
first show that $\mathfrak{z}(\mathfrak{\mathfrak{sl}}(n|n)^{\underline{h}})\rightarrow\mathfrak{z}(\mathfrak{g}^{h})$
is surjective. Suppose there exist $x\in\mathfrak{z}(\mathfrak{g}^{h})$
such that $\underline{x}\notin\mathfrak{z}(\mathfrak{\mathfrak{sl}}(n|n)^{\underline{h}})$.
Then $[\underline{x},\underline{y}]\in\mathfrak{z}(\mathfrak{\mathfrak{sl}}(n|n))$
for all $\underline{y}\in\mathfrak{\mathfrak{sl}}(n|n)^{\underline{h}}$
and $[\underline{x},\underline{y}]\neq0$ for some $\underline{y}$.
Denote by $\mathfrak{h}$ (resp. $\underline{\mathfrak{h}}$) the
Cartan subalgebra of $\mathfrak{g}$ (resp. $\mathfrak{\mathfrak{sl}}(n|n)$).
Let $\underline{d}=\mathrm{diag}(d_{1},\dots,d_{2n})\in\underline{\mathfrak{h}}$,
$x=\left[x_{ij}\right]\in\mathfrak{g}$. We have that the commutator
between $d$ and $x$ equals to $\left[(d_{i}-d_{j})x_{ij}\right]$.
Note that $\left[(d_{i}-d_{j})x_{ij}\right]=0$ if and only if $x_{ij}=0$
when $d_{i}\neq d_{j}$. Pick $d\in\mathfrak{h}\subseteq\mathfrak{g}^{h}$
such that $d_{i}\neq d_{j}$ for all $i\neq j$. Then $\mathfrak{g}^{d}=\mathfrak{h}$
and thus $\mathfrak{z}(\mathfrak{g}^{h})\subseteq\mathfrak{h}$. For
$x\in\mathfrak{h}$, $y=y_{\bar{0}}+y_{\bar{1}}\in\mathfrak{g}^{h}$,
it is impossible for $[\underline{x},\underline{y}]\in\mathfrak{z}(\mathfrak{\mathfrak{sl}}(n|n))$
and $[\underline{x},\underline{y}]\neq0$. This is because $[\underline{x},\underline{y}]\in\mathfrak{z}(\mathfrak{\mathfrak{sl}}(n|n))$
and $[\underline{x},\underline{y}]\neq0$ implies $[\underline{x},\underline{y_{\bar{0}}}]\in\mathfrak{z}(\mathfrak{\mathfrak{sl}}(n|n))$
and $[\underline{x},\underline{y_{\bar{0}}}]\neq0$. Note that $[\underline{x},\underline{y_{\bar{0}}}]$
lies in the derived subalgebra $[\mathfrak{sl}(n)\oplus\mathfrak{sl}(n),\mathfrak{sl}(n)\oplus\mathfrak{sl}(n)]$
but $\mathfrak{z}(\mathfrak{\mathfrak{sl}}(n|n))\notin[\mathfrak{sl}(n)\oplus\mathfrak{sl}(n),\mathfrak{sl}(n)\oplus\mathfrak{sl}(n)]$.
Therefore, we deduce that $\mathfrak{z}(\mathfrak{\mathfrak{sl}}(n|n)^{\underline{h}})\rightarrow\mathfrak{z}(\mathfrak{g}^{h})$
is surjective. 

According to \cite[Section 3.4]{Han-sl-osp}, an element in $\mathfrak{z}(\mathfrak{\mathfrak{sl}}(n|n)^{\underline{h}})$
is of the form
\[
\begin{pmatrix}a_{-\lambda_{1}+1}I_{r_{-\lambda_{1}+1}|s_{-\lambda_{1}+1}} & \dotsb & 0\\
\vdots & \ddots & \vdots\\
0 & \dotsb & a_{\lambda_{1}-1}I_{r_{\lambda_{1}-1}|s_{\lambda_{1}-1}}
\end{pmatrix}
\]
 Since $\mathfrak{z}(\mathfrak{\mathfrak{sl}}(n|n)^{\underline{h}})\rightarrow\mathfrak{z}(\mathfrak{g}^{h})$
is surjective, we have that an element in $\mathfrak{z}(\mathfrak{g}^{h})$
is of the form 
\[
\begin{pmatrix}a_{-\lambda_{1}+1}I_{r_{-\lambda_{1}+1}|s_{-\lambda_{1}+1}} & \dotsb & 0\\
\vdots & \ddots & \vdots\\
0 & \dotsb & a_{\lambda_{1}-1}I_{r_{\lambda_{1}-1}|s_{\lambda_{1}-1}}
\end{pmatrix}+\mathfrak{z}\left(\mathfrak{\mathfrak{sl}}(n|n)\right)
\]
 for some $a_{i}\in\mathbb{Z}$ such that $\sum_{i}a_{i}(r_{i}-s_{i})=0$.
Therefore, we have that $\dim\mathfrak{z}(\mathfrak{g}^{h})=\lambda_{1}-1$
if $r_{i}=s_{i}$ for all $i$ and $\dim\mathfrak{z}(\mathfrak{g}^{h})=\lambda_{1}-2$
if $r_{i}\neq s_{i}$ for some $i$.\qed

We next state the $\mathfrak{psl}(n|n)$ version of Theorem \ref{thm:g^e,z(g^e)}
and prove it. Recall that the sub-labelled Dynkin diagram $\Delta_{0}$
is the 2-free core of $\Delta$ obtained by removing all nodes with
labels equal to $2$ from $\Delta$. Let $k>0$ be minimal such that
the $k$th column in $P$ contains no boxes. Then we define $\tau$
to be the number of $i$ such that $r_{i}=s_{i}\neq0$ for all $|i|>k$.
We also define 
\[
\sigma=\begin{cases}
0 & \text{if }\sum_{|i|<k}r_{i}\neq\sum_{|i|<k}s_{i}\\
1 & \text{if }\sum_{|i|<k}r_{i}=\sum_{|i|<k}s_{i}
\end{cases}.
\]

\begin{thm}
Let $\mathfrak{g}=\mathfrak{g}_{\bar{0}}\oplus\mathfrak{g}_{\bar{1}}=\mathfrak{psl}(n|n)$,
$n>1$ and $e\in\mathfrak{g}_{\bar{0}}$ be nilpotent. Let $G=\{(A,B)\in\mathrm{GL}_{n}(\mathbb{C})\times\mathrm{GL}_{n}(\mathbb{C}):\mathrm{det}(A)=\mathrm{det}(B)\}/\{aI_{n|n}:a\in\mathbb{C}^{\times}\}\}$,
which is a reductive algebraic group over $\mathbb{C}$ such that
$\mathrm{Lie}(G)=\mathfrak{g}_{\bar{0}}$. Let $\Delta_{0}$ be the
2-free core of $\Delta$. Denote by $\mathfrak{g}_{0}$ be the subalgebra
of $\mathfrak{g}$ generated by the root vectors corresponding to
the simple roots in $\Delta_{0}$. Let $G_{0}$ be a closed connected
subgroup of $G$ such that $\mathrm{Lie}(G_{0})=(\mathfrak{g}_{0})_{\bar{0}}$.
Denote by $e_{0}$ a nilpotent element in $(\mathfrak{g}_{0})_{\bar{0}}$
whose corresponding labelled Dynkin diagram is $\Delta_{0}$. Then
we have 

(1) $\dim\mathfrak{g}^{e}-\dim\mathfrak{g}_{0}^{e_{0}}=n_{2}\left(\Delta\right)$;

(2) (i)When $\Delta$ has no label equal to $1$, 
\[
\dim\left(\mathfrak{z}\left(\mathfrak{g}^{e}\right)\right)^{G^{e}}-\dim\left(\mathfrak{z}\left(\mathfrak{g}_{0}^{e_{0}}\right)\right)^{G_{0}^{e_{0}}}=\begin{cases}
n_{2}\left(\Delta\right)-\tau+1=0 & \text{if }r_{i}=s_{i}\text{ for all }i\\
n_{2}\left(\Delta\right)-\tau & \text{otherwise}
\end{cases};
\]

(ii) When $\Delta$ has some label equal to $1$, 
\[
\dim\left(\mathfrak{z}\left(\mathfrak{g}^{e}\right)\right)^{G^{e}}-\dim\left(\mathfrak{z}\left(\mathfrak{g}_{0}^{e_{0}}\right)\right)^{G_{0}^{e_{0}}}=\begin{cases}
n_{2}\left(\Delta\right)-\tau & \text{if }r_{i}=s_{i}\text{ for all }i\text{ and }\sigma=1\\
n_{2}\left(\Delta\right)-\sigma-\tau & \text{otherwise}
\end{cases}.
\]
\end{thm}

\begin{proof}
\noindent We divide our analysis into two cases. 

\textbf{Case $1$:} When $\Delta$ has no label equal to $1$. We
have shown that $n_{2}\left(\Delta\right)=\lambda_{1}-1$ in the proof
of Theorem \ref{thm:LLD=000026centre}. In this case, since $\Delta_{0}$
has all labels equal to $0$, we have that $e_{0}=0$ and thus $\mathfrak{g}_{0}^{e_{0}}=\mathfrak{g}_{0}$.
Based on \cite[Section 3.4]{Han-sl-osp}, we know that $\mathfrak{g}_{0}=\left(\bigoplus_{i\in\mathbb{Z}}\mathfrak{sl}(r_{i}|s_{i})\right)/\mathfrak{z}(\mathfrak{sl}(n|n))$.
Thus
\[
\dim\mathfrak{g}_{0}^{e_{0}}=\dim\mathfrak{g}_{0}=\left(\sum_{i\in\mathbb{Z}}\dim\mathfrak{sl}(r_{i}|s_{i})\right)-1=\left(\sum_{i\in\mathbb{Z}:c_{i}>0}\left(c_{i}^{2}-1\right)\right)-1
\]
 by Proposition \ref{prop:dimg^e~c_i}. We have that $\dim\mathfrak{g}^{e}=\sum_{i\in\mathbb{Z}:c_{i}>0}c_{i}^{2}-2$.
Therefore, 
\begin{align*}
\dim\mathfrak{g}^{e}-\dim\mathfrak{g}_{0}^{e_{0}} & =\sum_{i\in\mathbb{Z}:c_{i}>0}c_{i}^{2}-2-\left(\sum_{i\in\mathbb{Z}:c_{i}>0}\left(c_{i}^{2}-1\right)\right)+1\\
 & =-1+\sum_{i\in\mathbb{Z}:c_{i}>0}1=\lambda_{1}-1=n_{2}(\Delta)
\end{align*}
 as there are in total $\lambda_{1}$ columns in $P$ with non-zero
boxes. Define $\mathfrak{I}=\{i:r_{i}=s_{i}\neq0\}$ and thus $\tau=|\mathfrak{I}|$.
Since $I_{r_{i}|s_{i}}\in\mathfrak{z}(\mathfrak{sl}(r_{i}|s_{i}))$
for all $i\in\mathfrak{I}$, we have that $\{I_{r_{i}|s_{i}}:i\in\mathfrak{I}\}$
spans $\mathfrak{z}(\mathfrak{g}_{0})$ and is a basis of $\mathfrak{z}(\mathfrak{g}_{0})$
if there exist some $j\notin\mathfrak{I}$. However, if $\mathfrak{I}=\{-\lambda_{1}+1,-\lambda_{1}+3,\dots,\lambda_{1}-3,\lambda_{1}-1\}$,
i.e. $\tau=\lambda_{1}$, we get a basis of $\mathfrak{z}(\mathfrak{g}_{0})$
with one fewer element as $\sum_{i}I_{r_{i}|s_{i}}\in\mathfrak{z}(\mathfrak{sl}(n|n))$
maps to zero. Hence, we have that 
\[
\dim\mathfrak{z}(\mathfrak{g}_{0}^{e_{0}})=\dim\mathfrak{z}(\mathfrak{g}_{0})=\begin{cases}
\tau-1 & \text{if }r_{i}=s_{i}\text{ for all }i\\
\tau & \text{otherwise.}
\end{cases}
\]

Therefore, we obtain that 
\begin{align*}
\dim\mathfrak{z}(\mathfrak{g}^{e})-\dim\mathfrak{z}(\mathfrak{g}_{0}^{e_{0}}) & =\dim\left(\mathfrak{z}(\mathfrak{g}^{e})\right)^{G^{e}}-\dim\left(\mathfrak{z}(\mathfrak{g}_{0}^{e_{0}})\right)^{G_{0}^{e_{0}}}\\
 & =\begin{cases}
\lambda_{1}-\tau=n_{2}(\Delta)-\tau+1=0 & \text{if }r_{i}=s_{i}\text{ for all }i\\
\lambda_{1}-1-\tau=n_{2}(\Delta)-\tau & \text{otherwise.}
\end{cases}
\end{align*}
 Note that if $r_{i}=s_{i}$ for all $i$, then $n_{2}(\Delta)-\tau+1=0$
as there are in total $\lambda_{1}$ columns with $r_{i}=s_{i}\neq0$.

\textbf{Case 2:} When $\Delta$ has some label equal to $1$. Recall
that $k>0$ is the minimal integer such that $c_{k}=0$. Note that
the Dynkin pyramid has $2\lambda_{1}+1$ columns with column number
from $-\lambda_{1}$ to $\lambda_{1}$. Based on \cite[Subsection 3.4]{Han-sl-osp},
we have $n_{2}(\Delta)=\lambda_{1}-k$. Next we look at $\mathfrak{g}_{0}^{e_{0}}$.
We know that 
\begin{align*}
\mathfrak{g}_{0} & \cong\biggl(\mathfrak{sl}\left(r_{-\lambda_{1}+1}|s_{-\lambda_{1}+1}\right)\oplus\dotsb\oplus\mathfrak{sl}\left(r_{-k-1}|s_{-k-1}\right)\oplus\mathfrak{sl}\left(\sum_{i=-k+1}^{k-1}r_{i}\bigg|\sum_{i=-k+1}^{k-1}s_{i}\right)\\
 & \oplus\mathfrak{sl}\left(r_{k+1}|s_{k+1}\right)\oplus\dotsb\oplus\mathfrak{sl}\left(r_{\lambda_{1}-1}|s_{\lambda_{1}-1}\right)\biggr)/\mathfrak{z}\left(\mathfrak{sl}\left(n|n\right)\right),
\end{align*}
 and the projection of $e_{0}$ in each $\mathfrak{sl}(r_{i}|s_{i})$
is zero for all $|i|>k$. By \cite[Subsection 3.4]{Han-sl-osp} we
have that 
\[
\dim\mathfrak{sl}\left(\sum_{i=-k+1}^{k-1}r_{i}\vert\sum_{i=-k+1}^{k-1}s_{i}\right)^{e_{0}}=\sum_{i=-k+1}^{k-1}c_{i}^{2}+\sum_{i=-k+1}^{k-1}c_{i}c_{i+1}-1
\]
 and $\dim\mathfrak{sl}(r_{i}|s_{i})=c_{i}^{2}-1$ for $|i|>k$. Since
the Dynkin pyramid is symmetric, we have
\[
\dim\mathfrak{g}_{0}^{e_{0}}=\sum_{i=-k+1}^{k-1}c_{i}^{2}+\sum_{i=-k+1}^{k-1}c_{i}c_{i+1}-1+2(c_{k+1}^{2}-1)+\dots+2(c_{\lambda_{1}-1}^{2}-1)-1.
\]
 Therefore, we obtain that 
\begin{align*}
\dim\mathfrak{g}^{e}-\dim\mathfrak{g}_{0}^{e_{0}} & =\sum_{i=-\lambda_{1}+1}^{\lambda_{1}-1}c_{i}^{2}+\sum_{i=-\lambda_{1}+1}^{\lambda_{1}-1}c_{i}c_{i+1}-2-\left(\sum_{i=-\lambda_{1}+1}^{\lambda_{1}-1}c_{i}^{2}-1-\sum_{i=1}^{\lambda_{1}-k}1-1\right)\\
 & =\lambda_{1}-k=n_{2}(\Delta).
\end{align*}

Next we investigate $\mathfrak{z}(\mathfrak{g}_{0}^{e_{0}})$. According
to \cite[Theorem 13]{Han-sl-osp}, 
\[
\dim\mathfrak{z}\left(\mathfrak{sl}\left(\sum_{i=-k+1}^{k-1}r_{i}\vert\sum_{i=-k+1}^{k-1}s_{i}\right)^{e_{0}}\right)=\begin{cases}
k-1 & \text{if }\sum_{|i|<k}r_{i}\neq\sum_{|i|<k}s_{i};\\
k & \text{if }\sum_{|i|<k}r_{i}=\sum_{|i|<k}s_{i}.
\end{cases}
\]
 Note that if $r_{i}\neq s_{i}$ for all $|i|>k$, then $\mathfrak{z}(\mathfrak{sl}(r_{i}|s_{i}))=0$
and thus $\dim\mathfrak{z}(\mathfrak{g}_{0}^{e_{0}})=k-1+\sigma$.
However, when there exist $r_{i}=s_{i}$ for some $|i|>k$, we have
that $I_{r_{i}|s_{i}}\in\mathfrak{z}(\mathfrak{sl}(r_{i}|s_{i}))$.
Moreover, if $r_{i}=s_{i}$ for all $i$ and $\sigma=1$, we have
that the $n\times n$ identity matrix maps to zero. Hence, 
\[
\dim\mathfrak{z}(\mathfrak{g}_{0}^{e_{0}})=\begin{cases}
k-1+\tau & \text{if }r_{i}=s_{i}\text{ for all i and }\sigma=1;\\
k-1+\sigma+\tau & \text{otherwise}.
\end{cases}
\]
Therefore, we obtain that 
\begin{align*}
\dim\mathfrak{z}(\mathfrak{g}^{e})-\dim\mathfrak{z}(\mathfrak{g}_{0}^{e_{0}}) & =\dim\left(\mathfrak{z}(\mathfrak{g}^{e})\right)^{G^{e}}-\dim\left(\mathfrak{z}(\mathfrak{g}_{0}^{e_{0}})\right)^{G_{0}^{e_{0}}}\\
 & =\begin{cases}
\lambda_{1}-k-\tau=n_{2}(\Delta)-\tau & \text{if }r_{i}=s_{i}\text{ for all i and }\sigma=1;\\
\lambda_{1}-k-\sigma-\tau=n_{2}(\Delta)-\sigma-\tau & \text{otherwise}.
\end{cases}
\end{align*}
\end{proof}

\section{Reachability and the Panyushev property in $\mathfrak{sl}(m|n)$
and $\mathfrak{psl}(n|n)$\label{sec:sl(m|n)}}

\noindent Let $V=V_{\bar{0}}\oplus V_{\bar{1}}$ be a finite-dimensional
vector superspace such that $\dim V_{\bar{0}}=m$, $\dim V_{1}=n$
and $m\neq n$. Let $\mathfrak{g}=\mathfrak{sl}\left(V\right)=\mathfrak{sl}\left(m|n\right)=\mathfrak{g}_{\bar{0}}\oplus\mathfrak{g}_{\bar{1}}$
and $e\in\mathfrak{g}_{\bar{0}}$ be nilpotent such that the partition
corresponding to $e$ is $\lambda$, which is defined in (\ref{eq:partition 1,2,3...}).
Fix an $\mathfrak{sl}(2)$-triple $\{e,h,f\}\subseteq\mathfrak{g}_{\bar{0}}$,
we recall that $\mathfrak{g}$ can be decomposed into the direct sum
of $\mathrm{ad}h$-eigenspaces $\mathfrak{g}=\bigoplus_{l\in\mathbb{Z}}\mathfrak{g}(l)$,
there also exists an $\mathrm{ad}h$-eigenspace decomposition for
$\mathfrak{g}^{e}$ such that $\mathfrak{g}^{e}=\bigoplus_{l\geq0}\mathfrak{g}^{e}(l)$.

In this section, we show that $e$ is reachable if and only if the
subalgebra $\mathfrak{g}^{e}(\ge1):=\bigoplus_{l\ge1}\mathfrak{g}^{e}(l)$
of $\mathfrak{g}^{e}$ is generated by $\mathfrak{g}^{e}(1)$, i.e.
$e$ satisfies the Panyushev property. We consider $\mathfrak{gl}(V)$
frequently in the following subsections, hence, let us fix the notation
$\bar{\mathfrak{g}}=\mathfrak{gl}\left(V\right)$.

\subsection{ad$h$-grading on $\bar{\mathfrak{g}}^{e}$}

We start by calculating the $\mathrm{ad}h$-eigenvalue of basis elements
in $\bar{\mathfrak{g}}^{e}$. Define a grading of the vector space
$V$ such that $V=\bigoplus_{t\in\mathbb{Z}}V(t)$ and $V(t)=\{v\in V:hv=tv\}$.
\begin{lem}
\label{lem:adh-eigenvalue of basis elements}Let $\{\xi_{i}^{j,k}:1\leq i,j\leq r+s,\max\{\lambda_{j}-\lambda_{i},0\}\leq k\leq\lambda_{j}-1\}$
be the basis of $\bar{\mathfrak{g}}^{e}$ as defined in (\ref{eq:basis gl}).
Then

(i) $x\in\bar{\mathfrak{g}}(l)$ if and only if \textup{$x\left(V\left(t\right)\right)\subseteq V(t+l)$
for all $t\in\mathbb{Z}$;}

(ii) $\xi_{i}^{j,k}\in\bar{\mathfrak{g}}^{e}(\lambda_{i}-\lambda_{j}+2k)$.
\end{lem}

\begin{proof}
\begin{singlespace}
(i) Let $x\in\bar{\mathfrak{g}}$ and $v\in V(t)$.

Suppose $[h,x]=lx$. Then 
\[
h\left(xv\right)=x\left(hv\right)+[h,x]v=txv+lxv=(t+l)xv.
\]
Hence, we know that $xv\in V(t+l)$ as required. 

Conversely, suppose $xv\in V(t+l)$. Then $h\left(xv\right)=(t+l)xv$.
Hence, we have that 
\[
h\left(xv\right)=x\left(hv\right)+[h,x]v=\left(t+l\right)xv,
\]
 which implies that $[h,x]v=lxv$ for all $v\in V(t)$. Therefore,
we have that $[h,x]=lx$. 
\end{singlespace}

(ii) Recall that a basis of $V$ is given in Subsection \ref{subsec:Centralizer-of-nilpoent},
i.e. $V=\mathrm{Span}\{e^{k}v_{i}:1\leq i\leq r+s,0\leq k\leq\lambda_{i}-1\}$.
Based on the way that the Dynkin pyramid P is constructed, we have
that $v_{i}\in V(-\lambda_{i}+1)$ and $\xi_{i}^{j,k}\left(v_{i}\right)=e^{k}v_{j}\in V\left(-\lambda_{j}+2k+1\right)$.
Hence, we obtain that $\xi_{i}^{j,k}\in\bar{\mathfrak{g}}^{e}(\lambda_{i}-\lambda_{j}+2k)$
by part (i).
\end{proof}

\subsection{Reachable nilpotent elements $e$ in $\mathfrak{sl}(m|n)_{\bar{0}}$}

In the following theorem we give a characterisation of reachable nilpotent
elements in $\mathfrak{g}_{\bar{0}}$ in terms of the partition $\lambda$.
\begin{thm}
\label{thm:reachable gl}If a nilpotent element $e\in\mathfrak{g}_{\bar{0}}$
is reachable, then the corresponding partition $\lambda$ satisfies
the condition $\lambda_{i}-\lambda_{i+1}\in\{0,1\}$ for all $1\leq i\leq r+s$
(by convention $\lambda_{i}=0$ for $i>r+s$).
\end{thm}

\begin{proof}
Note that $e=\sum_{i=1}^{r+s}\xi_{i}^{i,1}$ and $\xi_{i}^{i,1}\in\mathfrak{g}^{e}(2)$.
We know that $[\mathfrak{g}^{e}(i),\mathfrak{g}^{e}(j)]\subseteq\mathfrak{g}^{e}(i+j)$.
Thus $\xi_{i}^{i,1}$ only can be obtained from $[\mathfrak{g}^{e}(1),\mathfrak{g}^{e}(1)]+[\mathfrak{g}^{e}(0),\mathfrak{g}^{e}(2)]$.

By part (ii) in Lemma \ref{lem:adh-eigenvalue of basis elements},
we know that the ad$h$-grading of $\xi_{i}^{j,k}$ is $\lambda_{i}-\lambda_{j}+2k$.
Since the condition on $k$ is $\max\{\lambda_{j}-\lambda_{i},0\}\leq k\leq\lambda_{j}-1$,
we have that $\lambda_{i}-\lambda_{j}+2k\geq|\lambda_{i}-\lambda_{j}|$.
Notice that $\lambda_{i}-\lambda_{j}+2k=0$ if and only if $\lambda_{i}=\lambda_{j}$
and $k=0$. Hence, a basis for $\mathfrak{g}^{e}(0)$ is $\{\xi_{i}^{j,0}:\lambda_{i}=\lambda_{j},i\neq j\}\cup\{\lambda_{i+1}\xi_{i}^{i,0}-(-1)^{\bar{i}}\lambda_{i}\xi_{i+1}^{i+1,0}\}$.
Similarly, a basis for $\mathfrak{g}^{e}(1)$ consists of elements
$\{\xi_{i}^{j,0}:\lambda_{i}=\lambda_{j}+1\}\cup\{\xi_{j}^{i,1}:\lambda_{i}=\lambda_{j}+1\}$
and a basis for $\mathfrak{g}^{e}(2)$ consists of elements $\lbrace\xi_{i}^{j,0},\xi_{j}^{i,2}:\lambda_{i}=\lambda_{j}+2\rbrace\cup\lbrace\xi_{i}^{j,1}:\lambda_{i}=\lambda_{j}\rbrace.$

The only way to get $\xi_{i}^{i,1}$ as a commutator from $[\mathfrak{g}^{e}(1),\mathfrak{g}^{e}(1)]$
or $[\mathfrak{g}^{e}(0),\mathfrak{g}^{e}(2)]$ is when $\lambda_{i}=\lambda_{j}+1$
or $\lambda_{i}=\lambda_{j}$. If $\lambda_{i}=\lambda_{j}+1$. we
consider $\xi_{j}^{i,1},\xi_{i}^{j,0}\in\mathfrak{g}^{e}(1)$, we
have that $[\xi_{j}^{i,1},\xi_{i}^{j,0}]=\xi_{i}^{i,1}-(-1)^{\bar{a}\bar{b}}\xi_{j}^{j,1}$
where $\bar{a}$ (resp. $\bar{b}$) is the parity of $\xi_{i}^{j,0}$
(resp. $\xi_{j}^{i,1}$). If $\lambda_{i}=\lambda_{j}$ and $i\neq j$,
we consider $\xi_{j}^{i,1}\in\mathfrak{g}^{e}(2)$, $\xi_{i}^{j,0}\in\mathfrak{g}^{e}(0)$
and obtain that $[\xi_{j}^{i,1},\xi_{i}^{j,0}]=\xi_{i}^{i,1}-(-1)^{\bar{a}\bar{b}}\xi_{j}^{j,1}$. 

Since $e\in[\mathfrak{g}^{e},\mathfrak{g}^{e}]$, then $e$ can be
written as a linear combination of $\xi_{j}^{j,1}-(-1)^{\bar{a}\bar{b}}\xi_{i}^{i,1}$
for $1\leq i,j\leq r+s$, $\lambda_{i}=\lambda_{j}+1$ or $\lambda_{i}=\lambda_{j}$.
If there exists $i$ such that $\lambda_{i}\geq\lambda_{i+1}+2,$
then such a linear combination can not exist. Therefore, we deduce
that $\lambda_{i}-\lambda_{i+1}\in\{0,1\}$ for all $1\leq i\leq r+s$.
\end{proof}
\begin{cor}
\label{cor:reachable gl}Theorem \ref{thm:reachable gl} also holds
for $\bar{\mathfrak{g}}$. 
\end{cor}

\begin{proof}
For two elements $x',y'\in\bar{\mathfrak{g}}^{e}$, we have $x'=x+aI_{m|n}$
and $y'=y+aI_{m|n}$ for $x,y\in\mathfrak{g}^{e}$, then $[x',y']=[x+aI_{m|n},y+aI_{m|n}]=[x,y]$.
This shows that $[\bar{\mathfrak{g}}^{e},\bar{\mathfrak{g}}^{e}]\subseteq[\mathfrak{g}^{e},\mathfrak{g}^{e}]$.
Thus $e\in[\bar{\mathfrak{g}}^{e},\bar{\mathfrak{g}}^{e}]$ implies
$e\in[\mathfrak{g}^{e},\mathfrak{g}^{e}]$. Hence, the Corollary holds
by Theorem \ref{thm:reachable gl}.
\end{proof}
In the following theorem, we apply a similar argument as in \cite[Section 4]{Panyushev}
to determine when $\mathfrak{g}^{e}\left(\geq1\right)$ is generated
by $\mathfrak{g}^{e}\left(1\right)$ as a Lie superalgebra. 
\begin{thm}
\label{thm:Panyushev gl}Suppose $\lambda$ satisfies the condition
$\lambda_{i}-\lambda_{i+1}\in\{0,1\}$ for all $1\leq i\leq r+s$
(by convention $\lambda_{i}=0$ for $i>r+s$). Then $\mathfrak{g}^{e}\left(\geq1\right)$
is generated by $\mathfrak{g}^{e}\left(1\right)$ as a Lie superalgebra.
\end{thm}

\begin{proof}
We argue by induction on $\dim V$.

Based on the way that the Dynkin pyramid $P$ is constructed, we know
that $e$ sends $v_{i}$ to $v_{j}$ if the box labelled by $j$ is
the left neighbour of $i$ and sends $v_{i}$ to zero if the box labelled
by $v_{i}$ has no left neighbour. Thus we have that $\ker\left(e\right)=\mathrm{Span}\{e^{\lambda_{i}-1}v_{i}:i=1,\dots,r+s\}$.
Let $U=V/\ker\left(e\right)$ and $\hat{\mathfrak{g}}=\mathfrak{sl}\left(U\right)$.
Then we have that $\dim U=\dim V-(r+s)$ and $e$ induces the nilpotent
transformation $\hat{e}$ of $U$ with the corresponding partition
$\hat{\lambda}=(\lambda_{1}-1,\dots,\lambda_{r+s}-1)$. The mapping
from $V$ to $U$ induces the $\mathbb{Z}_{2}$-graded homomorphism
$\phi:\mathfrak{g}^{e}\rightarrow\hat{\mathfrak{g}}^{\hat{e}}$. Denote
$W=\mathrm{Span}\{v_{i}:i=1,\dots,r+s\}$. Note that $\mathrm{im}\left(e\right)=\mathrm{Span}\{e^{k}v_{i}:i=1,\dots,r+s,0<k\leq\lambda_{i}-1\}$.
Thus $V=W\oplus\mathrm{im}\left(e\right)$ and 
\[
\ker\left(\phi\right)=\{\xi_{i}^{j,k}\in\mathfrak{g}^{e}:\xi_{i}^{j,k}\left(W\right)\subseteq\ker\left(e\right)\}.
\]
 We also deduce that $\ker\left(\phi\right)\cong\mathrm{Hom}\left(W,\ker\left(e\right)\right)$
and $\dim\ker\left(\phi\right)=\left(r+s\right)^{2}$. Next we claim
that $\phi$ preserves the $\mathbb{Z}$-grading of both $\mathfrak{g}^{e}$
and $\hat{\mathfrak{g}}^{\hat{e}}$. According to Lemma \ref{lem:adh-eigenvalue of basis elements},
the $h$-eigenvalue of $e^{k}v_{i}$ is $-\lambda_{i}+2k+1$ for $0\leq k\leq\lambda_{i}-1$
and $1\leq i\leq r+s$. Moreover, we know that $\xi_{i}^{j,k}\in\mathfrak{g}^{e}(l)$
if and only if the difference of the $h$-eigenvalues of $\xi_{i}^{j,k}\left(v_{t}\right)$
and $v_{t}$ is equal to $l$ according to Lemma \ref{lem:adh-eigenvalue of basis elements}.
Let $\{\hat{e},\hat{h},\hat{f}\}$ be an $\mathfrak{sl}(2)$-triple
in $\hat{\mathfrak{g}}$ containing $\hat{e}$, then the $\hat{h}$-eigenvalue
of $e^{k}v_{i}$ is $-\lambda_{i}+2k+2$ for $\lambda_{i}\geq2$,
$0\leq k\leq\lambda_{i}-2$ and $1\leq i\leq r+s$. Note that this
shift does not affect the difference of the eigenvalues for elements
in $U$. Hence, we obtain that $\phi\left(\mathfrak{g}^{e}(l)\right)\subseteq\hat{\mathfrak{g}}^{\hat{e}}(l)$
which completes the claim.

Now we start the induction steps. It is clear that the partition $\hat{\lambda}$
also satisfies the condition of the theorem. By the inductive hypothesis,
we assume that $\hat{\mathfrak{g}}^{\hat{e}}(\geq1)$ is generated
by $\hat{\mathfrak{g}}^{\hat{e}}(1)$. Subsequently, it suffices to
prove that each element of $\ker\left(\phi\right)$ is contained in
the subspace generated by $\mathfrak{g}^{e}(1)$. Note that $\ker\left(\phi\right)$
has a basis $\{\xi_{i}^{j,\lambda_{j}-1}:1\leq i,j\leq r+s\}$. By
Lemma \ref{lem:adh-eigenvalue of basis elements}, we know that $\xi_{i}^{j,\lambda_{j}-1}\in\mathfrak{g}^{e}\left(\lambda_{i}+\lambda_{j}-2\right)$.
Observe that $\xi_{i}^{j,\lambda_{j}-1}\in\mathfrak{g}^{e}\left(0\right)$
if and only if $\lambda_{i}=\lambda_{j}=1$, which we exclude from
further considerations.

If $\lambda_{j}>1$, we can take $\lambda_{t}=\lambda_{j}-1$. Then
we consider elements $\xi_{t}^{j,1}\in\mathfrak{g}^{e}\left(1\right)$
and $\xi_{i}^{t,\lambda_{t}-1}\in\ker\left(\phi\right)$. Notice that
$[\xi_{t}^{j,1},\xi_{i}^{t,\lambda_{t}-1}]=\xi_{i}^{j,\lambda_{j}-1}\in\ker\left(\phi\right)$.
Since the $\mathrm{ad}h$-eigenvalue of $\xi_{i}^{t,\lambda_{t}-1}$
is one less than that of $\xi_{i}^{j,\lambda_{j}-1}$, we can show
that each element of $\ker\left(\phi\right)$ lies in the subspace
generated by $\mathfrak{g}^{e}(1)$ by an ascending induction on the
$\mathrm{ad}h$-eigenvalue of $\xi_{i}^{j,\lambda_{j}-1}$.

If $\lambda_{j}=1$, then $\lambda_{i}>1$ and we can take $\lambda_{t}=\lambda_{i}-1$.
Then we consider elements $\xi_{t}^{i,1}\in\mathfrak{g}^{e}\left(1\right)$
and $\xi_{j}^{t,\lambda_{t}-1}\in\ker\left(\phi\right)$. Notice that
$[\xi_{t}^{i,1},\xi_{j}^{t,\lambda_{t}-1}]=\xi_{j}^{i,\lambda_{i}-1}\in\ker\left(\phi\right)$.
Similarly we deduce that each element of $\ker\left(\phi\right)$
lies in the subspace generated by $\mathfrak{g}^{e}(1)$.

Therefore, we conclude that $\mathfrak{g}^{e}\left(\geq1\right)$
is generated by $\mathfrak{g}^{e}\left(1\right)$ as a Lie superalgebra.
\end{proof}
Now, we are ready to present the main result of this section.
\begin{thm}
\label{thm:three conditions are equivalent}Let $\mathfrak{g}=\mathfrak{sl}(m|n)=\mathfrak{g}_{\bar{0}}\oplus\mathfrak{g}_{\bar{1}}$
and $e\in\mathfrak{g}_{\bar{0}}$ be nilpotent. Let $\lambda$ be
the partition with respect to $e$ as defined in (\ref{eq:partition 1,2,3...}).
Then the following conditions are equivalent:

(1) $\lambda_{i}-\lambda_{i+1}\in\{0,1\}$ for all $1\leq i<r+s$
and $\lambda_{r+s}=1$;

(2) $\mathfrak{g}^{e}\left(\geq1\right)$ is generated by $\mathfrak{g}^{e}\left(1\right)$
as a Lie superalgebra;

(3) $e\in[\mathfrak{g}^{e}(1),\mathfrak{g}^{e}(1)]$;

(4) $e\in[\mathfrak{g}^{e},\mathfrak{g}^{e}]$.
\end{thm}

\begin{proof}
The implication from (1) to (2) is the content of Theorem \ref{thm:Panyushev gl}.

By definition, $\mathfrak{g}^{e}\left(\geq1\right)$ is generated
by $\mathfrak{g}^{e}\left(1\right)$ implies that $e\in[\mathfrak{g}^{e}(1),\mathfrak{g}^{e}(1)]$
as $e\in\mathfrak{g}^{e}(2)$. 

The implication from (3) to (4) is clear.

The implication from (4) to (1) is the content of Theorem \ref{thm:reachable gl}.
\end{proof}
\begin{rem}
Theorem \ref{thm:three conditions are equivalent} also holds for
$\mathfrak{psl}(n|n)$ since $\mathfrak{psl}(n|n)^{e}$ only has one
basis element less than $\mathfrak{g}^{e}$ and this element has ad$h$-grading
$0$, which does not affect the above content of proof. 

By Theorem \ref{thm:three conditions are equivalent}, we have the
following Corollary.
\end{rem}

\begin{cor}
Let $e\in\mathfrak{g}_{\bar{0}}$ be nilpotent. Then $e$ is reachable
if and only if $\lambda$ satisfies the condition $\lambda_{i}-\lambda_{i+1}\in\{0,1\}$
for all $1\leq i<r+s$ and $\lambda_{r+s}=1$.
\end{cor}

\section{Reachability and the Panyushev property for $\mathfrak{osp}(m|2n)$\label{sec:osp(m|2n)}}

\noindent Let $V=V_{\bar{0}}\oplus V_{\bar{1}}$ be a $\mathbb{Z}_{2}$-graded
vector space over $\mathbb{C}$ such that $\dim V_{\bar{0}}=m$ and
$\dim V_{\bar{1}}=2n$. Suppose there exists a non-degenerate supersymmetric
bilinear form on $V$. Let $\mathfrak{g}=\mathfrak{g}_{\bar{0}}\oplus\mathfrak{g}_{\bar{1}}=\mathfrak{osp}(V)$
and we write $\mathfrak{osp}(m|2n)$ for $\mathfrak{osp}(V)$. Let
$e\in\mathfrak{g}_{\bar{0}}$ be nilpotent with the corresponding
partition $\lambda$. Note that $\lambda=(\lambda_{1},\dots,\lambda_{r+s})$
is a partition of $(m|2n)$ defined in a similar way as in (\ref{eq:partition 1,2,3...})
such that $\sum_{\ensuremath{\bar{i}}=\bar{0}}\lambda_{i}=m$, $\sum_{\ensuremath{\bar{i}}=\bar{1}}\lambda_{i}=2n$.

In this section, we first recall a basis for the centralizer $\mathfrak{g}^{e}$
as given in \cite[Subsection 4.5]{Han-sl-osp}. Then we decompose
$\mathfrak{g}^{e}$ into subspaces, namely $\mathfrak{N}$, $N_{1}$
and $N_{2}$, which will be introduced in Subsection \ref{subsec:Decomposing-g^e}.
We further decompose $[\mathfrak{g}^{e},\mathfrak{g}^{e}]$ into a
finite collection of those subspaces of $\mathfrak{g}^{e}$ and deduce
that for a reachable element $e\in\mathfrak{g}_{\bar{0}}$, basis
elements in a complement to $[\mathfrak{g}^{e},\mathfrak{g}^{e}]$
in $\mathfrak{g}^{e}$ all lie in $\mathfrak{g}^{e}(0)$. Last, we
show that $e\in\mathfrak{g}_{\bar{0}}$ is reachable if and only if
$e\in\mathfrak{g}_{\bar{0}}$ satisfies the Panyushev property, i.e.
$\mathfrak{g}^{e}\left(\geq1\right)$ is generated by $\mathfrak{g}^{e}(1)$.

\subsection{Decomposing $\mathfrak{g}^{e}$ into subspaces\label{subsec:Decomposing-g^e}}

By \cite[Subsection 4.5]{Han-sl-osp}, there exists an involution
$i\mapsto i^{*}$ on the set $\{1,\dots,r+s\}$ such that $i^{*}\in\{i-1,i,i+1\}$
for all $1\leq i\leq r+s$, $\lambda_{i^{*}}=\lambda_{i}$ and $\bar{i^{*}}=\bar{i}$.
Let $\xi_{i}^{j,k}$ be the basis element in $\mathfrak{gl}(m|2n)^{e}$
defined as in (\ref{eq:basis gl}). We know that $\mathfrak{g}^{e}=\mathfrak{g}\cap\mathfrak{gl}(m|2n)^{e}$
and a basis for $\mathfrak{g}^{e}$ contains the following elements:

\begin{doublespace}
\noindent 
\[
\xi_{i}^{i^{*},\lambda_{i}-1-k}\text{ for }0\leq k\leq\lambda_{i}-1,\ k\text{ is odd if }\bar{i}=\bar{0}\text{ and }k\text{ is even if }\bar{i}=\bar{1};
\]

\noindent 
\[
\xi_{i}^{j,\lambda_{j}-1-k}+\varepsilon_{i,j,k}\xi_{j^{*}}^{i^{*},\lambda_{i}-1-k}\text{ for all }0\leq k\leq\min\{\lambda_{i},\lambda_{j}\}-1
\]
 such that $\varepsilon_{i,j,k}\in\{\pm1\}$ and 
\begin{equation}
\varepsilon_{i,j,k}=(-1)^{\lambda_{j}-k-\bar{x}\bar{i}}\theta_{j}\theta_{i}.\label{eq:epsilon i,j,k}
\end{equation}
 Note that $\bar{x}\in\{\bar{0},\bar{1}\}$ is the parity of $\xi_{i}^{j,\lambda_{j}-1-k}$
and $\theta_{j},\theta_{i}\in\{\pm1\}$ can be determined explicitly
by \cite[Section 3.2]{Jantzen}. An element $\xi_{i}^{j,\lambda_{j}-1-k}+\varepsilon_{i,j,k}\xi_{j^{*}}^{i^{*},\lambda_{i}-1-k}$
is in $\mathfrak{g}_{\bar{0}}^{e}$ if $\bar{i}=\bar{j}$ or in $\mathfrak{g}_{\bar{1}}^{e}$
if $\bar{i}\neq\bar{j}$. Moreover, for an element $\xi_{i}^{i^{*},\lambda_{i}-1-k}$,
we have $i=i^{*}$ if $\bar{i}=\bar{0}$, $(-1)^{\lambda_{i}}=-1$
or $\bar{i}=\bar{1},(-1)^{\lambda_{i}}=1$. 
\end{doublespace}

Next we decompose $\mathfrak{g}^{e}$ into subspaces in a similar
way as in \cite[Section 2]{Premet=000026Topley-derived} in preparation
of investigating the structure of $[\mathfrak{g}^{e},\mathfrak{g}^{e}]$
later in Subsection \ref{subsec:A-decomposition-for-=00005Bg^e,g^e=00005D}.
Let us denote 
\begin{align}
H & :=\left\{ \xi_{i}^{i,\lambda_{i}-1-k}:i=i^{*},0\leq k<\lambda_{i},\lambda_{i}-k\text{ even}\right\} \label{eq:decompose g^e}\\
 & \cup\left\{ \xi_{i}^{i,\lambda_{i}-1-k}+\varepsilon_{i,i,k}\xi_{i^{*}}^{i^{*},\lambda_{i}-1-k}:i\neq i^{*},0\leq k<\lambda_{i}\right\} ;\nonumber 
\end{align}
\[
N_{(1)}:=\left\{ \xi_{i}^{i^{*},\lambda_{i}-1-k}:i\neq i^{*},0\leq k<\lambda_{i},\lambda_{i}-k\text{ odd}\right\} ;
\]
\[
N_{(2)}=:\left\{ \xi_{i}^{j,\lambda_{j}-1-k}+\varepsilon_{i,j,k}\xi_{j^{*}}^{i^{*},\lambda_{i}-1-k}:i<j,j\neq i^{*},0\leq k<\min\{\lambda_{i},\lambda_{j}\}\right\} ,
\]
 and let $\mathfrak{N}=\mathrm{Span}\left(H\right),$ $N_{1}=\mathrm{Span}\left(N_{(1)}\right)$,
$N_{2}=\mathrm{Span}\left(N_{(2)}\right)$. Note that the set $H\cup N_{(1)}\cup N_{(2)}$
forms a basis for $\mathfrak{g}^{e}$. 

Let $\mathfrak{N}=\mathfrak{N}_{0}\oplus\mathfrak{N}_{1}\subseteq\mathfrak{g}_{\bar{0}}^{e}$
be a decomposition for $\mathfrak{N}$ such that 
\[
\mathfrak{N}_{0}=\mathrm{Span}\left\{ \left\{ \xi_{i}^{i,\lambda_{i}-1-k}:i\neq i^{*},\lambda_{i}-k\text{ even}\right\} \cup\left\{ \xi_{i}^{i,\lambda_{i}-1-k}+\xi_{i^{*}}^{i^{*},\lambda_{i}-1-k}:i\neq i^{*},\lambda_{i}-k\text{ even}\right\} \right\} 
\]
 and 
\[
\mathfrak{N}_{1}=\mathrm{Span}\left\{ \xi_{i}^{i,\lambda_{i}-1-k}-\xi_{i^{*}}^{i^{*},\lambda_{i}-1-k}:i\neq i^{*},\lambda_{i}-k\text{ odd}\right\} .
\]

\subsection{A decomposition for $[\mathfrak{g}^{e},\mathfrak{g}^{e}]$\label{subsec:A-decomposition-for-=00005Bg^e,g^e=00005D}}

According to our notation in (\ref{eq:decompose g^e}), we have that
$\mathfrak{g}^{e}=\mathfrak{N}\oplus N_{1}\oplus N_{2}$ and $\mathfrak{N}=\mathfrak{N}_{0}\oplus\mathfrak{N}_{1}$.
In this subsection, we obtain a decomposition for $[\mathfrak{g}^{e},\mathfrak{g}^{e}]$. 

We first describe $N_{1}\cap[\mathfrak{g}^{e},\mathfrak{g}^{e}]$
and $N_{2}\cap[\mathfrak{g}^{e},\mathfrak{g}^{e}]$ using a similar
approach as in \cite[Subsection 2.3]{Premet=000026Topley-derived}.
Let $N_{(2)}^{-}$ be the set consisting of basis elements 
\[
\left\{ \xi_{i}^{i+1,0}+\varepsilon_{i,i+1,\lambda_{i+1}-1}\xi_{(i+1)^{*}}^{i^{*},\lambda_{i}-\lambda_{i+1}}\in N_{(2)}:i=i^{*},i+1=(i+1)^{*},\lambda_{i-1}>\lambda_{i}\geq\lambda_{i+1}>\lambda_{i+2}\right\} .
\]
 Define $N_{(2)}^{+}$ to be the complement of basis elements to $N_{(2)}^{-}$
in $N_{2}$ and let $N_{2}^{-}=\mathrm{Span}\left(N_{(2)}^{-}\right)$,
$N_{2}^{+}=\mathrm{Span}\left(N_{(2)}^{+}\right)$. 
\begin{prop}
\label{prop:long proof} We have that $N_{1}\subseteq[\mathfrak{g}^{e},\mathfrak{g}^{e}]$
and $N_{2}\cap[\mathfrak{g}^{e},\mathfrak{g}^{e}]=N_{2}^{+}$.
\end{prop}

\begin{proof}
$N_{1}$ has a basis consisting of elements $\xi_{i}^{i^{*},\lambda_{i}-1-k}$
with $i\neq i^{*}$, $0\leq k<\lambda_{i}$ and $\lambda_{i}-k$ odd.
Suppose $\xi_{i}^{i^{*},\lambda_{i}-1-k}\in N_{1}$. We consider $\xi_{i}^{i,0}+\varepsilon_{i,i,\lambda_{i}-1}\xi_{i^{*}}^{i^{*},0}\in\mathfrak{N}$
and compute 
\[
[\xi_{i}^{i^{*},\lambda_{i}-1-k},\xi_{i}^{i,0}+\varepsilon_{i,i,\lambda_{i}-1}\xi_{i^{*}}^{i^{*},0}]=\xi_{i}^{i^{*},\lambda_{i}-1-k}-\varepsilon_{i,i,\lambda_{i}-1}\xi_{i}^{i^{*},\lambda_{i}-1-k}.
\]
 We have that $\varepsilon_{i,i,\lambda_{i}-1}=(-1)^{1}=-1$ by (\ref{eq:epsilon i,j,k}).
Thus the above commutator equals to $2\xi_{i}^{i^{*},\lambda_{i}-1-k}\in[\mathfrak{g}^{e},\mathfrak{g}^{e}]$.
Hence, we deduce that $N_{1}=[N_{1},\mathfrak{N}]\subseteq[\mathfrak{g}^{e},\mathfrak{g}^{e}]$.

The proof of the second part is divided into two steps. Note that
$N_{(2)}^{+}=N_{(2)}^{+,1}\cup N_{(2)}^{+,2}\cup N_{(2)}^{+,3}$ where
\[
N_{(2)}^{+,1}=\left\{ \xi_{i}^{j,\lambda_{j}-1-k}+\varepsilon_{i,j,k}\xi_{j^{*}}^{i^{*},\lambda_{i}-1-k}:j>i,j\neq i+1,0\leq k\leq\lambda_{j}-1\right\} ,
\]
\[
N_{(2)}^{+,2}=\left\{ \xi_{i}^{i+1,\lambda_{i+1}-1-k}+\varepsilon_{i,i+1,k}\xi_{(i+1)^{*}}^{i^{*},\lambda_{i}-1-k}:0\leq k\leq\lambda_{i+1}-2\right\} ,
\]
\begin{align*}
N_{(2)}^{+,3} & =\left\{ \xi_{i}^{i+1,0}+\varepsilon_{i,i+1,\lambda_{i+1}-1}\xi_{(i+1)^{*}}^{i^{*},\lambda_{i}-\lambda_{i+1}}:i\neq i^{*}\text{ or }i+1\neq(i+1)^{*}\right\} \\
 & \cup\left\{ \xi_{i}^{i+1,0}+\varepsilon_{i,i+1,\lambda_{i+1}-1}\xi_{(i+1)^{*}}^{i^{*},\lambda_{i}-\lambda_{i+1}}:i=i^{*},i+1=(i+1)^{*},\lambda_{i-1}=\lambda_{i}\text{ or }\lambda_{i+1}=\lambda_{i+2}\right\} .
\end{align*}

\textbf{Step 1}: We demonstrate that $N_{2}^{+}\subseteq[\mathfrak{g}^{e},\mathfrak{g}^{e}]$
by showing all elements in $N_{(2)}^{+}$ can be obtained from a commutator
of two elements in $\mathfrak{g}^{e}$. For the remainder of the proof
we fix $j>i\neq j^{*}$.

$\bullet$ We show that $\xi_{i}^{j,\lambda_{j}-1-k}+\varepsilon_{i,j,k}\xi_{j^{*}}^{i^{*},\lambda_{i}-1-k}\in[\mathfrak{g}^{e},\mathfrak{g}^{e}]$
for $j\neq i+1$ and $0\leq k\leq\lambda_{j}-1$.

(i) When $j\neq j^{*}$, we have 
\[
[\xi_{j}^{j,0}+\varepsilon_{j,j,\lambda_{j}-1}\xi_{j^{*}}^{j^{*},0},\xi_{i}^{j,\lambda_{j}-1-k}+\varepsilon_{i,j,k}\xi_{j^{*}}^{i^{*},\lambda_{i}-1-k}]=\xi_{i}^{j,\lambda_{j}-1-k}-\varepsilon_{j,j,\lambda_{j}-1}\varepsilon_{i,j,k}\xi_{j^{*}}^{i^{*},\lambda_{i}-1-k}.
\]
 Since $\varepsilon_{j,j,\lambda_{j}-1}=-1$ by equation (\ref{eq:epsilon i,j,k}),
the above commutator equals to $\xi_{i}^{j,\lambda_{l}-1-k}+\varepsilon_{i,j,k}\xi_{j^{*}}^{i^{*},\lambda_{i}-1-k}$
 and thus $\xi_{i}^{j,\lambda_{l}-1-k}+\varepsilon_{i,j,k}\xi_{j^{*}}^{i^{*},\lambda_{i}-1-k}$ lies in $[\mathfrak{g}^{e},\mathfrak{g}^{e}]$. Similarly, when
$i\neq i^{*}$, we have 
\[
[\xi_{i}^{i,0}+\varepsilon_{i,i,\lambda_{i}-1}\xi_{i^{*}}^{i^{*},0},\xi_{i}^{j,\lambda_{j}-1-k}+\varepsilon_{i,j,k}\xi_{j^{*}}^{i^{*},\lambda_{i}-1-k}]=\varepsilon_{i,i,\lambda_{i}-1}\varepsilon_{i,j,k}\xi_{j^{*}}^{i^{*},\lambda_{i}-1-k}-\xi_{i}^{j,\lambda_{j}-1-k}.
\]
 Again we know that $\varepsilon_{i,i,\lambda_{i}-1}=-1$, thus the
above commutator equals to $-(\xi_{i}^{j,\lambda_{j}-1-k}+\varepsilon_{i,j,k}\xi_{j^{*}}^{i^{*},\lambda_{i}-1-k})$
 and thus $\xi_{i}^{j,\lambda_{j}-1-k}+\varepsilon_{i,j,k}\xi_{j^{*}}^{i^{*},\lambda_{i}-1-k}$ lies in $[\mathfrak{g}^{e},\mathfrak{g}^{e}]$. 

(ii) When $j=j^{*}$ and $i=i^{*}$, we know that there exists $l$
such that $j>l>i$ as $j>i$ and $j\neq i+1.$ We compute 
\[
[\xi_{l}^{j,\lambda_{j}-1-k}+\varepsilon_{l,j,k}\xi_{j}^{l^{*},\lambda_{l}-1-k},\xi_{i}^{l,0}+\varepsilon_{i,l,\lambda_{l}-1}\xi_{l^{*}}^{i,\lambda_{i}-\lambda_{l}}]=\xi_{i}^{j,\lambda_{j}-1-k}-(-1)^{\overline{a_{1}}\overline{b_{1}}}\varepsilon_{l,j,k}\varepsilon_{i,l,\lambda_{l}-1}\xi_{j}^{i,\lambda_{i}-1-k}
\]
 where $\overline{a_{1}}$ (resp. $\overline{b_{1}}$) is the parity
of $\xi_{l^{*}}^{i,\lambda_{i}-\lambda_{l}}$ (resp. $\xi_{j}^{l^{*},\lambda_{l}-1-k}$).
Using (\ref{eq:epsilon i,j,k}) we obtain that 
\[
-(-1)^{\overline{a_{1}}\overline{b_{1}}}\varepsilon_{l,j,k}\varepsilon_{i,l,\lambda_{l}-1}=\varepsilon_{i,j,k}.
\]
 Hence, the above commutator equals to $\xi_{i}^{j,\lambda_{j}-1-k}+\varepsilon_{i,j,k}\xi_{j}^{i,\lambda_{i}-1-k}$
and we deduce that $\xi_{i}^{j,\lambda_{j}-1-k}+\varepsilon_{i,j,k}\xi_{j}^{i,\lambda_{i}-1-k}\in[\mathfrak{g}^{e},\mathfrak{g}^{e}]$.

$\bullet$ We show that $\xi_{i}^{i+1,\lambda_{i+1}-1-k}+\varepsilon_{i,i+1,k}\xi_{(i+1)^{*}}^{i^{*},\lambda_{i}-1-k}\in[\mathfrak{g}^{e},\mathfrak{g}^{e}]$
for $0\leq k\leq\lambda_{i+1}-2$.

(iii) When $i+1\neq(i+1)^{*}$ or $i\neq i^{*}$, we deduce that $\xi_{i}^{i+1,\lambda_{i+1}-1-k}+\varepsilon_{i,i+1,k}\xi_{(i+1)^{*}}^{i^{*},\lambda_{i}-1-k}\in[\mathfrak{g}^{e},\mathfrak{g}^{e}]$
for $0\leq k\leq\lambda_{i+1}-2$ from part (i) by letting $j=i+1$
in part (i).

(iv) When $i+1=(i+1)^{*}$ and $i=i^{*}$, we have 
\[
[\xi_{i+1}^{i+1,1},\xi_{i}^{i+1,\lambda_{i+1}-1-k}+\varepsilon_{i,i+1,k}\xi_{i+1}^{i,\lambda_{i}-1-k}]=\xi_{i}^{i+1,\lambda_{i+1}-k}-\varepsilon_{i,i+1,k}\xi_{i+1}^{i,\lambda_{i}-k}.
\]
Take $k'=k-1$, then the above expression is equal to $\xi_{i}^{i+1,\lambda_{i+1}-1-k'}+\varepsilon_{i,i+1,k'}\xi_{i+1}^{i,\lambda_{i}-1-k'}\in[\mathfrak{g}^{e},\mathfrak{g}^{e}]$
for all $0\leq k'\leq\lambda_{j}-2$ because $\varepsilon_{i,i+1,k'}=-\varepsilon_{i,i+1,k}$
by (\ref{eq:epsilon i,j,k}). 

$\bullet$ We show that $\xi_{i}^{i+1,0}+\varepsilon_{i,i+1,\lambda_{i+1}-1}\xi_{(i+1)^{*}}^{i^{*},\lambda_{i}-\lambda_{i+1}}\in[\mathfrak{g}^{e},\mathfrak{g}^{e}]$
for $i+1\neq(i+1)^{*}$ or $i\neq i^{*}$ or $i+1=(i+1)^{*}$, $i=i^{*}$
and $\lambda_{i-1}=\lambda_{i}\text{ or }\lambda_{i+1}=\lambda_{i+2}$.

(v) When $i+1\neq(i+1)^{*}$ or $i\neq i^{*}$, we have that $\xi_{i}^{i+1,0}+\varepsilon_{i,i+1,\lambda_{i+1}-1}\xi_{(i+1)^{*}}^{i^{*},\lambda_{i}-\lambda_{i+1}}\in[\mathfrak{g}^{e},\mathfrak{g}^{e}]$
from part $(i)$ by letting $j=i+1$ and $k=\lambda_{i+1}-1$ in part
(i).

(vi) When $i+1=(i+1)^{*}$, $i=i^{*}$ and $\lambda_{i}=\lambda_{i-1}$,
we have 
\begin{align*}
[\xi_{i-1}^{i+1,0}+\varepsilon_{i-1,i+1,\lambda_{i+1}-1}\xi_{i+1}^{(i-1)^{*},\lambda_{i}-\lambda_{i+1}},\xi_{i}^{i-1,0}+\varepsilon_{i,i-1,\lambda_{i}-1}\xi_{(i-1)^{*}}^{i,0}]\\
=\xi_{i}^{i+1,0}-(-1)^{\overline{a_{2}}\overline{b_{2}}}\varepsilon_{i-1,i+1,\lambda_{i+1}-1}\varepsilon_{i,i-1,\lambda_{i}-1}\xi_{i+1}^{i,\lambda_{i}-\lambda_{i+1}}
\end{align*}
 where $\overline{a_{2}}$ (resp. $\overline{b_{2}}$) is the parity
of $\xi_{(i-1)^{*}}^{i,0}$ (resp. $\xi_{i+1}^{(i-1)^{*},\lambda_{i}-\lambda_{i+1}}$).
Using (\ref{eq:epsilon i,j,k}) we obtain that $-(-1)^{\overline{a_{2}}\overline{b_{2}}}\varepsilon_{i-1,i+1,\lambda_{i+1}-1}\varepsilon_{i,i-1,\lambda_{i}-1}=\varepsilon_{i,i+1,\lambda_{i+1}-1}$.
Hence, the above commutator equals to $\xi_{i}^{i+1,0}+\varepsilon_{i,i+1,\lambda_{i+1}-1}\xi_{i+1}^{i,\lambda_{i}-\lambda_{i+1}}$
and we deduce that $\xi_{i}^{i+1,0}+\varepsilon_{i,i+1,\lambda_{i+1}-1}\xi_{i+1}^{i,\lambda_{i}-\lambda_{i+1}}\in[\mathfrak{g}^{e},\mathfrak{g}^{e}]$.
Similarly, when $i+1=(i+1)^{*}$, $i=i^{*}$ and $\lambda_{i+1}=\lambda_{i+2}$,
we have 
\[
[\xi_{i+2}^{i+1,0}+\varepsilon_{i+2,i+1,\lambda_{i+1}-1}\xi_{i+1}^{(i+2)^{*},0},\xi_{i}^{i+2,0}+\varepsilon_{i,i+2,\lambda_{i+2}-1}\xi_{(i+2)^{*}}^{i,\lambda_{i}-\lambda_{i+2}}]=\xi_{i}^{i+1,0}+\varepsilon_{i,i+1,\lambda_{i+1}}\xi_{i+1}^{i,\lambda_{i}-\lambda_{i+1}}.
\]
 This implies that $\xi_{i}^{i+1,0}+\varepsilon_{i,i+1,\lambda_{i+1}}\xi_{i+1}^{i,\lambda_{i}-\lambda_{i+1}}\in[\mathfrak{g}^{e},\mathfrak{g}^{e}]$.

In parts (i)--(vi), we checked all elements in $N_{(2)}^{+}$ and
we obtain that $N_{(2)}^{+}\subseteq[\mathfrak{g}^{e},\mathfrak{g}^{e}]$.
This implies that $N_{2}^{+}\subseteq[\mathfrak{g}^{e},\mathfrak{g}^{e}]$. 

\textbf{Step 2}: We demonstrate that $N_{2}\cap[\mathfrak{g}^{e},\mathfrak{g}^{e}]=N_{2}^{+}$.
Choose an element $\xi_{i}^{j,\lambda_{j}-1-k_{1}}+\varepsilon_{i,j,k_{1}}\xi_{j^{*}}^{i^{*},\lambda_{i}-1-k_{1}}\in N_{2}$.

$\bullet$ We start with showing that $[\mathfrak{N},N_{2}]\subseteq N_{2}^{+}$. 

(i). For $\xi_{l}^{l,\lambda_{l}-1-k_{2}}\in\mathfrak{N}$ with $l=l^{*},0\leq k_{2}<\lambda_{l}$,
$\lambda_{l}-k_{2}$ even, the commutator $[\xi_{i}^{j,\lambda_{j}-1-k_{1}}+\varepsilon_{i,j,k}\xi_{j^{*}}^{i^{*},\lambda_{i}-1-k_{1}},\xi_{l}^{l,\lambda_{l}-1-k_{2}}]$
is nonzero if and only if $i=l$ or $j=l$. Assume $i=l$ first. Then
\[
[\xi_{i}^{j,\lambda_{j}-1-k_{1}}+\varepsilon_{i,j,k_{1}}\xi_{j^{*}}^{i^{*},\lambda_{i}-1-k_{1}},\xi_{l}^{l,\lambda_{l}-1-k_{2}}]=\xi_{i}^{j,\lambda_{j}+\lambda_{l}-2-k_{1}-k_{2}}-\varepsilon_{i,j,k_{1}}\xi_{j^{*}}^{i^{*},\lambda_{i}+\lambda_{l}-2-k_{1}-k_{2}}.
\]
Let $k=k_{1}+k_{2}-(\lambda_{l}-1)$. The above commutator equals
to $\xi_{i}^{j,\lambda_{j}-1-k}+\varepsilon_{i,j,k}\xi_{j^{*}}^{i^{*},\lambda_{i}-1-k}$
because $\varepsilon_{i,j,k}=-\varepsilon_{i,j,k_{1}}$ by (\ref{eq:epsilon i,j,k}).
This element $\xi_{i}^{j,\lambda_{j}-1-k}+\varepsilon_{i,j,k}\xi_{j^{*}}^{i^{*},\lambda_{i}-1-k}$
lies either in $N_{2}^{-}$ or $N_{2}^{+}$. Since $\lambda_{l}-k_{2}$
is even, we have that $k_{2}\leq\lambda_{l}-2$ and $k\leq k_{1}-1<\lambda_{j}-1$.
This forces $\xi_{i}^{j,\lambda_{j}-1-k}+\varepsilon_{i,j,k}\xi_{j^{*}}^{i^{*},\lambda_{i}-1-k}$
to lie in $N_{2}^{+}$. Applying a similar approach to the case $j=l$,
we also obtain that $[\xi_{i}^{j,\lambda_{j}-1-k_{1}}+\varepsilon_{i,j,k_{1}}\xi_{j^{*}}^{i^{*},\lambda_{i}-1-k_{1}},\xi_{l}^{l,\lambda_{l}-1-k_{2}}]\in N_{2}^{+}$. 

(ii). For $\xi_{l}^{l,\lambda_{l}-1-k_{2}}+\varepsilon_{l,l,k_{2}}\xi_{l^{*}}^{l^{*},\lambda_{l}-1-k_{2}}\in\mathfrak{N}$
with $l\neq l^{*}$, $0\leq k_{2}<\lambda_{l}$, similarly the commutator
$[\xi_{i}^{j,\lambda_{j}-1-k_{1}}+\varepsilon_{i,j,k_{1}}\xi_{j^{*}}^{i^{*},\lambda_{i}-1-k_{1}},\xi_{l}^{l,\lambda_{l}-1-k_{2}}+\varepsilon_{l,l,k_{2}}\xi_{l^{*}}^{l^{*},\lambda_{l}-1-k_{2}}]$
is nonzero if and only if $i=l$ or $j=l$. Assume $i=l$. Then 
\begin{align*}
[\xi_{i}^{j,\lambda_{j}-1-k_{1}}+\varepsilon_{i,j,k_{1}}\xi_{j^{*}}^{i^{*},\lambda_{i}-1-k_{1}},\xi_{l}^{l,\lambda_{l}-1-k_{2}}+\varepsilon_{l,l,k_{2}}\xi_{l^{*}}^{l^{*},\lambda_{l}-1-k_{2}}]\\
=\xi_{i}^{j,\lambda_{j}+\lambda_{l}-2-k_{1}-k_{2}}-\varepsilon_{i,j,k_{1}}\varepsilon_{l,l,k_{2}}\xi_{j^{*}}^{i^{*},\lambda_{i}+\lambda_{l}-2-k_{1}-k_{2}}.
\end{align*}
 Again let $k=k_{1}+k_{2}-(\lambda_{l}-1)$. The above commutator
equals to $\xi_{i}^{j,\lambda_{j}-1-k}+\varepsilon_{i,j,k}\xi_{j^{*}}^{i^{*},\lambda_{i}-1-k}$
because $\varepsilon_{i,j,k}=-\varepsilon_{i,j,k_{1}}\varepsilon_{l,l,k_{2}}$
by (\ref{eq:epsilon i,j,k}). Since $i=l$, $l\neq l^{*}$, we have
that $i\neq i^{*}$. This forces $\xi_{i}^{j,\lambda_{j}-1-k}+\varepsilon_{i,j,k}\xi_{j^{*}}^{i^{*},\lambda_{i}-1-k}$
to lie in $N_{2}^{+}$. Applying a similar approach to the case $j=l$,
we also obtain that $[\xi_{i}^{j,\lambda_{j}-1-k_{1}}+\varepsilon_{i,j,k_{1}}\xi_{j^{*}}^{i^{*},\lambda_{i}-1-k_{1}},\xi_{l}^{l,\lambda_{l}-1-k_{2}}+\varepsilon_{l,l,k_{2}}\xi_{l^{*}}^{l^{*},\lambda_{l}-1-k_{2}}]\in N_{2}^{+}$. 

Hence, we deduce that $[\mathfrak{N},N_{2}]\subseteq N_{2}^{+}$.

$\bullet$ Next we show that $[N_{1},N_{2}]\subseteq N_{2}^{+}$.

For an element $\xi_{l}^{l^{*},\lambda_{l}-1-k_{2}}\in N_{1}$ with
$l\neq l^{*}$, $0\leq k_{2}<\lambda_{l}$, $\lambda_{l}-k_{2}$ is
odd. The commutator $[\xi_{i}^{j,\lambda_{j}-1-k_{1}}+\varepsilon_{i,j,k_{1}}\xi_{j^{*}}^{i^{*},\lambda_{i}-1-k_{1}},\xi_{l}^{l^{*},\lambda_{l}-1-k_{2}}]$
is nonzero if and only if $i=l^{*}$ or $j=l$. Let $k=k_{1}+k_{2}-(\lambda_{l}-1)$.

Assume $i=l^{*}$, then 
\[
[\xi_{i}^{j,\lambda_{j}-1-k_{1}}+\varepsilon_{i,j,k_{1}}\xi_{j^{*}}^{i^{*},\lambda_{i}-1-k_{1}},\xi_{l}^{l^{*},\lambda_{l}-1-k_{2}}]=\xi_{l}^{j,\lambda_{j}-1-k}-\varepsilon_{i,j,k_{1}}\xi_{j^{*}}^{l^{*},\lambda_{i}-1-k}.
\]
 We compute that $\varepsilon_{l,j,k}=-\varepsilon_{i,j,k_{1}}$.
Thus the above commutator equals to $\xi_{i^{*}}^{j,\lambda_{j}-1-k}+\varepsilon_{i^{*},j,k}\xi_{j^{*}}^{i,\lambda_{i}-1-k}\in N_{2}$.
We further observe that $[\xi_{i}^{j,\lambda_{j}-1-k_{1}}+\varepsilon_{i,j,k_{1}}\xi_{j^{*}}^{i^{*},\lambda_{i}-1-k_{1}},\xi_{l}^{l^{*},\lambda_{l}-1-k_{2}}]\in N_{2}^{+}$
because $l\neq l^{*}$ based on our assumption and thus $i\neq i^{*}$.

Next assume $j=l$, then 
\[
[\xi_{i}^{j,\lambda_{j}-1-k_{1}}+\varepsilon_{i,j,k_{1}}\xi_{j^{*}}^{i^{*},\lambda_{i}-1-k_{1}},\xi_{l}^{l^{*},\lambda_{l}-1-k_{2}}]=\varepsilon_{i,j,k_{1}}\xi_{l}^{i^{*},\lambda_{i}-1-k}-\xi_{i}^{l^{*},\lambda_{j}-1-k}.
\]
 We compute that $\varepsilon_{i,j^{*},k}=-\varepsilon_{i,j,k_{1}}$.
Thus the above commutator equals to 
\[
-(\xi_{i}^{j^{*},\lambda_{j}-1-k}+\varepsilon_{i,j^{*},k}\xi_{j}^{i^{*},\lambda_{i}-1-k})\in N_{2}.
\]
 Note that $j=l$ and $l\neq l^{*}$, thus $j\neq j^{*}$ which implies
that 
\[
[\xi_{i}^{j,\lambda_{j}-1-k_{1}}+\varepsilon_{i,j,k_{1}}\xi_{j^{*}}^{i^{*},\lambda_{i}-1-k_{1}},\xi_{l}^{l^{*},\lambda_{l}-1-k_{2}}]\in N_{2}^{+}.
\]

$\bullet$ Now we demonstrate that $N_{2}\cap[N_{2},N_{2}]\subseteq N_{2}^{+}$
by showing for any two basis elements 
\[
\xi_{i}^{j,\lambda_{j}-1-k_{1}}+\varepsilon_{i,j,k_{1}}\xi_{j^{*}}^{i^{*},\lambda_{i}-1-k_{1}},\xi_{t}^{l,\lambda_{l}-1-k_{2}}+\varepsilon_{t,l,k_{2}}\xi_{l^{*}}^{t^{*},\lambda_{t}-1-k_{2}}\in N_{2},
\]
 the commutator between them lies in either $\mathfrak{N}$, $N_{1}$
or $N_{2}^{+}$.

We may assume that $j>i$ and $l>k$. Note that 
\begin{equation}
[\xi_{i}^{j,\lambda_{j}-1-k_{1}}+\varepsilon_{i,j,k_{1}}\xi_{j^{*}}^{i^{*},\lambda_{i}-1-k_{1}},\xi_{t}^{l,\lambda_{l}-1-k_{2}}+\varepsilon_{t,l,k_{2}}\xi_{l^{*}}^{t^{*},\lambda_{t}-1-k_{2}}]\label{eq:pro}
\end{equation}
 is nonzero when at least one of the following equalities hold: $i=l$,
$j=t$, $i^{*}=t$, $j^{*}=l$. 

We first consider the case when $i=l$. Since $i^{*}\neq j>i$ and
$l>t\neq l^{*}$, we have that $i\neq t^{*}$, $j\neq t$ and $j^{*}\neq l$.
Thus the commutator in (\ref{eq:pro}) is equal to $\xi_{t}^{j,\lambda_{j}-1-k}-(-1)^{\overline{a_{3}}\overline{b_{3}}}\varepsilon_{i,j,k_{1}}\varepsilon_{t,l,k_{2}}\xi_{j^{*}}^{t^{*},\lambda_{t}-1-k}$
where $\overline{a_{3}}$ (resp. $\overline{b_{3}}$) is the parity
of $\xi_{l^{*}}^{t^{*},\lambda_{t}-1-k_{2}}$ (resp. $\xi_{j^{*}}^{i^{*},\lambda_{i}-1-k_{1}}$).
By (\ref{eq:epsilon i,j,k}) we have that $-(-1)^{\overline{a_{3}}\overline{b_{3}}}\varepsilon_{i,j,k_{1}}\varepsilon_{t,l,k_{2}}=\varepsilon_{t,j,k}$
and the commutator in (\ref{eq:pro}) further equals to $\xi_{t}^{j,\lambda_{j}-1-k}+\varepsilon_{t,j,k}\xi_{j^{*}}^{t^{*},\lambda_{t}-1-k}\in N_{2}$.
Observe that this element does not lie in $N_{2}^{-}$ because for
$j>i=l>t$, it is impossible for $j=t+1$. Hence, the commutator in
equation (\ref{eq:pro}) lies in $N_{2}^{+}$. Applying a similar
argument to the case $j=t$ we also have that the commutator in equation
(\ref{eq:pro}) lies in $N_{2}^{+}$. 

Now suppose $i^{*}=t$, then $i\neq l$ and $j\neq t$. If $j^{*}=l$
then (\ref{eq:pro}) is equal to $\varepsilon_{t,l,k_{2}}((\xi_{j}^{j,\lambda_{j}+\lambda_{i}-2-k_{1}-k_{2}}+\varepsilon_{j,j,k_{1}+k_{2}-\lambda_{i}+1}\xi_{j^{*}}^{j^{*},\lambda_{j}+\lambda_{i}-2-k_{1}-k_{2}})-(\xi_{i}^{i,\lambda_{j}+\lambda_{i}-2-k_{1}-k_{2}}+\varepsilon_{i,i,k_{1}+k_{2}-\lambda_{j}+1}\xi_{i^{*}}^{i^{*},\lambda_{j}+\lambda_{i}-2-k_{1}-k_{2}}))\in\mathfrak{N}$.
If $j^{*}\neq l$ then (\ref{eq:pro}) is equal to 
\begin{equation}
\varepsilon_{t,l,k_{2}}(\xi_{l^{*}}^{j,\lambda_{j}+\lambda_{i}-2-k_{1}-k_{2}}+\varepsilon_{l^{*},j,k_{1}+k_{2}-\lambda_{i}+1}\xi_{j^{*}}^{l,\lambda_{l}+\lambda_{i}-2-k_{1}-k_{2}}).\label{eq:pro 1}
\end{equation}
 If $j=l$, then the expression in (\ref{eq:pro 1}) lies in $N_{1}$.
If $j\neq l$, then the expression in (\ref{eq:pro 1}) lies in $N_{2}$.
If this expression lies in $N_{2}^{-}$ we must have $j=j^{*}$, $l=l^{*}$,
$j=l^{*}+1$ and $\lambda_{l^{*}-1}\neq\lambda_{l^{*}}$. This implies
that $i<l$ and $\lambda_{l}>\lambda_{i}$. Since $0\leq k_{1}<\lambda_{j}$
and $0\leq k_{2}<\lambda_{l}$ by definition of $N_{2}$, we have
that $k_{1}+k_{2}-(\lambda_{i}-1)<\lambda_{j}-1$ which contrary to
the definition of $N_{2}^{-}$. Hence, the expression in (\ref{eq:pro 1})
lies in $N_{2}^{+}$. Note that the proof for the case $j^{*}=l$
is identical to the case $i^{*}=t$, hence is omitted.

Therefore, we obtain that $N_{2}\cap[N_{2},N_{2}]\subseteq N_{2}^{+}$.

$\bullet$ Now we are ready to show that $N_{2}\cap[\mathfrak{g}^{e},\mathfrak{g}^{e}]=N_{2}^{+}$.
Note that $[\mathfrak{N},\mathfrak{N}]=0$ by computing commutators
between basis elements in $\mathfrak{N}$. We also claim that $[N_{1},N_{1}]\subseteq\mathfrak{N}$.
For $\xi_{i}^{i^{*},\lambda_{i}-1-k_{1}},\xi_{j}^{j^{*},\lambda_{j}-1-k_{2}}\in N_{1}$
with $i\neq i^{*}$, $j\neq j^{*}$ and $\lambda_{i}-k_{1}$ and $\lambda_{j}-k_{2}$
are both odd. The commutator $[\xi_{i}^{i^{*},\lambda_{i}-1-k_{1}},\xi_{j}^{j^{*},\lambda_{j}-1-k_{2}}]$
is nonzero only when $i=j^{*}$. In which case we have that $[\xi_{i}^{i^{*},\lambda_{i}-1-k_{1}},\xi_{j}^{j^{*},\lambda_{j}-1-k_{2}}]=\xi_{j}^{j,\lambda_{j}+\lambda_{i}-2-k_{1}-k_{2}}-\xi_{j^{*}}^{j^{*},\lambda_{j}+\lambda_{i}-2-k_{1}-k_{2}}\in\mathfrak{N}$.

We have already showed that $N_{2}^{+}\subseteq N_{2}\cap[\mathfrak{g}^{e},\mathfrak{g}^{e}]$
in Step 1. Since $\mathfrak{g}^{e}=\mathfrak{N}\oplus N_{1}\oplus N_{2}$,
we further deduce that $N_{2}\cap[\mathfrak{g}^{e},\mathfrak{g}^{e}]\subseteq N_{2}^{+}$
based on the previous arguments. Therefore, we conclude that $N_{2}^{+}=N_{2}\cap[\mathfrak{g}^{e},\mathfrak{g}^{e}]$.
\end{proof}
Using a similar argument as in \cite[Theorem 6]{Premet=000026Topley-derived},
we give a decomposition for $[\mathfrak{g}^{e},\mathfrak{g}^{e}]$
in the following theorem.
\begin{thm}
\label{thm:=00005Bg^e,g^e=00005D}Let $\mathfrak{g}=\mathfrak{g}_{\bar{0}}\oplus\mathfrak{g}_{\bar{1}}=\mathfrak{osp}(m|2n)$
and $e\in\mathfrak{g}_{\bar{0}}$ be nilpotent. Then $[\mathfrak{g}^{e},\mathfrak{g}^{e}]=N_{1}\oplus N_{2}^{+}\oplus\mathfrak{N}_{1}\oplus(\mathfrak{N}_{0}\cap[N_{2},N_{2}])$.
\end{thm}

\begin{proof}
As shown in the proof of Proposition \ref{prop:long proof}, we have
that 
\[
[\mathfrak{g}^{e},\mathfrak{g}^{e}]=\left(N_{1}\cap[\mathfrak{g}^{e},\mathfrak{g}^{e}]\right)+\left(N_{2}\cap[\mathfrak{g}^{e},\mathfrak{g}^{e}]\right)+\left(\mathfrak{N}\cap[\mathfrak{g}^{e},\mathfrak{g}^{e}]\right).
\]
 Recall that $\mathfrak{N}=\mathfrak{N}_{0}\oplus\mathfrak{N}_{1}$.
Since $\mathfrak{N}_{1}\subseteq\mathfrak{g}_{\bar{0}}^{e}$, we have
that $\mathfrak{N}_{1}\subseteq[\mathfrak{g}^{e},\mathfrak{g}^{e}]$
by \cite[Proposition 3]{Premet=000026Topley-derived} and thus $\mathfrak{N}_{1}\cap[\mathfrak{g}^{e},\mathfrak{g}^{e}]=\mathfrak{N}_{1}$.
Based on the proof of Proposition \ref{prop:long proof}, we deduce
that $\mathfrak{N}_{0}\cap[\mathfrak{g}^{e},\mathfrak{g}^{e}]=\left(\mathfrak{N}_{0}\cap[N_{1},N_{1}]\right)+\left(\mathfrak{N}_{0}\cap[N_{2},N_{2}]\right)$.
Since $N_{1}\subseteq\mathfrak{g}_{\bar{0}}^{e}$, then $[N_{1},N_{1}]\subseteq\mathfrak{N}_{1}$
by \cite[Proposition 3]{Premet=000026Topley-derived}. Thus we deduce
that $\mathfrak{N}_{0}\cap[\mathfrak{g}^{e},\mathfrak{g}^{e}]=\mathfrak{N}_{0}\cap[N_{2},N_{2}].$
Hence, we have $\mathfrak{N}\cap[\mathfrak{g}^{e},\mathfrak{g}^{e}]=\mathfrak{N}_{1}+\left(\mathfrak{N}_{0}\cap[N_{2},N_{2}]\right)$
using the fact $\mathfrak{N}=\mathfrak{N}_{0}\oplus\mathfrak{N}_{1}$.
Note that $\left(N_{1}\cap[\mathfrak{g}^{e},\mathfrak{g}^{e}]\right)+\left(N_{2}\cap[\mathfrak{g}^{e},\mathfrak{g}^{e}]\right)=N_{1}+N_{2}^{+}$
by Proposition \ref{prop:long proof}. Therefore, we obtain that $[\mathfrak{g}^{e},\mathfrak{g}^{e}]=N_{1}\oplus N_{2}^{+}\oplus\mathfrak{N}_{1}\oplus(\mathfrak{N}_{0}\cap[N_{2},N_{2}])$.
\end{proof}

\subsection{The centralizer of a reachable nilpotent $e$ in $\mathfrak{osp}\left(m|2n\right)_{\bar{0}}$}

It is clear that if $e\in\mathfrak{g}_{\bar{0}}$ satisfies the Panyushev
property, then $e\in\mathfrak{g}_{\bar{0}}$ is reachable. Now suppose
$e\in\mathfrak{g}_{\bar{0}}$ is reachable. Since $\mathfrak{g}^{e}=\mathfrak{g}\cap\mathfrak{gl}(m|2n)^{e}$,
we have that $e\in[\mathfrak{gl}(m|2n)^{e},\mathfrak{gl}(m|2n)^{e}]$.
This implies that the partition $\lambda$ with respect to $e$ satisfies
the condition $\lambda_{i}-\lambda_{i+1}\in\{0,1\}$ for all $1\leq i\leq r+s$
by Theorem \ref{thm:reachable gl} and Corollary \ref{cor:reachable gl}.
In the following theorem, we apply a similar approach as in \cite[Prposition 3.3]{Modular}
in order to show that if $\lambda$ satisfies the above condition,
then the Panyushev property holds for $e$. Once this is done, we
have that $e\in\mathfrak{g}_{\bar{0}}$ is reachable if and only if
$e\in\mathfrak{g}_{\bar{0}}$ satisfies the Panyushev property.
\begin{thm}
Suppose $\lambda$ satisfies the condition $\lambda_{i}-\lambda_{i+1}\in\{0,1\}$
for all $1\leq i\leq r+s$ (by convention $\lambda_{i}=0$ for $i>r+s$).
Then $\mathfrak{g}^{e}(\geq1)$ is generated by $\mathfrak{g}^{e}(1)$
as a Lie superalgebra. 
\end{thm}

\begin{proof}
We prove the theorem by showing that $\mathfrak{g}^{e}$ is generated
by $\mathfrak{g}^{e}(0)$ and $\mathfrak{g}^{e}(1)$.

Suppose $\lambda$ satisfies the condition in the theorem. In this
case $\mathfrak{N}_{0}\cap[N_{2},N_{2}]$ in Theorem \ref{thm:=00005Bg^e,g^e=00005D}
is equal to $\mathfrak{N}_{0}$ by \cite[Section 2.4]{Premet=000026Topley-derived}.
Combining with the decomposition for $[\mathfrak{g}^{e},\mathfrak{g}^{e}]$
in Theorem \ref{thm:=00005Bg^e,g^e=00005D}, we deduce that a complement
to $[\mathfrak{g}^{e},\mathfrak{g}^{e}]$ in $\mathfrak{g}^{e}$ contains
basis elements 
\begin{equation}
\lbrace\xi_{i}^{i+1,0}+\varepsilon_{i,i+1,\lambda_{i+1}-1}\xi_{(i+1)^{*}}^{i^{*},\lambda_{i}-\lambda_{i+1}}\in N_{(2)}:i=i^{*},i+1=(i+1)^{*},\lambda_{i-1}>\lambda_{i}\geq\lambda_{i+1}>\lambda_{i+2}\rbrace.\label{eq:=00005Bg^e=00005D/g^e}
\end{equation}
 Since $\lambda_{i}-\lambda_{i+1}\in\{0,1\}$ for all $1\leq i\leq r+s$,
we have that basis elements in (\ref{eq:=00005Bg^e=00005D/g^e}) are
in $\mathfrak{g}^{e}(0)$ if $\lambda_{i}-\lambda_{i+1}=0$ or in
$\mathfrak{g}^{e}(1)$ if $\lambda_{i}-\lambda_{i+1}=1$ by Lemma
\ref{lem:adh-eigenvalue of basis elements}.(ii). Hence, in order
to prove the theorem, it suffices to show that $[\mathfrak{g}^{e}(1),\mathfrak{g}^{e}(c-1)]=\mathfrak{g}^{e}(c)$
for all $c>1$.

For $1\leq i,j\leq r+s$ and $|\lambda_{i}-\lambda_{j}|\geq2$, we
know that $\xi_{i}^{j,\lambda_{j}-1-k}+\varepsilon_{i,j,k}\xi_{j^{*}}^{i^{*},\lambda_{i}-1-k}\in\mathfrak{g}^{e}(\lambda_{i}+\lambda_{j}-2-2k)$.
We may assume that $i<j$. Let $\lambda_{i}+\lambda_{j}-2-2k=c$.
Since $\lambda$ satisfies the condition in the theorem, there exists
$l$ such that $i<l<j$ and $\lambda_{i}=\lambda_{l}+1>\lambda_{j}$.
Then $\xi_{l}^{j,\lambda_{j}-1-k}+\varepsilon_{l,j,k}\xi_{j^{*}}^{l^{*},\lambda_{l}-1-k}\in\mathfrak{g}^{e}(c-1)$,
$\xi_{i}^{l,0}+\varepsilon_{i,l,\lambda_{l}-1}\xi_{l^{*}}^{i^{*},1}\in\mathfrak{g}^{e}(1)$
and 
\[
[\xi_{l}^{j,\lambda_{j}-1-k}+\varepsilon_{l,j,k}\xi_{j^{*}}^{l^{*},\lambda_{l}-1-k},\xi_{i}^{l,0}+\varepsilon_{i,l,\lambda_{l}-1}\xi_{l^{*}}^{i^{*},1}]=\xi_{i}^{j,\lambda_{j}-1-k}-(-1)^{\overline{a_{4}}\overline{b_{4}}}\varepsilon_{l,j,k}\varepsilon_{i,l,\lambda_{l}-1}\xi_{j^{*}}^{i^{*},\lambda_{l}-k}
\]
 where $\overline{a_{4}}$ (resp. $\overline{b_{4}}$) is the parity
of $\xi_{l^{*}}^{i^{*},1}$ (resp. $\xi_{j^{*}}^{l^{*},\lambda_{l}-1-k}$).
By (\ref{eq:epsilon i,j,k}) we calculate that 
\[
-(-1)^{\overline{a_{4}}\overline{b_{4}}}\varepsilon_{l,j,k}\varepsilon_{i,l,\lambda_{l}-1}=\varepsilon_{i,j,k}
\]
 and thus 
\[
\xi_{i}^{j,\lambda_{j}-1-k}+\varepsilon_{i,j,k}\xi_{j^{*}}^{i^{*},\lambda_{i}-1-k}=[\xi_{l}^{j,\lambda_{j}-1-k}+\varepsilon_{l,j,k}\xi_{j^{*}}^{l^{*},\lambda_{l}-1-k},\xi_{i}^{l,0}+\varepsilon_{i,l,\lambda_{l}-1}\xi_{l^{*}}^{i^{*},1}]\in[\mathfrak{g}^{e}(c-1),\mathfrak{g}^{e}(1)].
\]

For $|\lambda_{i}-\lambda_{j}|=1$, similarly we may assume that $i<j$
and let $\xi_{i}^{j,\lambda_{j}-1-k}+\varepsilon_{i,j,k}\xi_{j^{*}}^{i^{*},\lambda_{i}-1-k}\in\mathfrak{g}^{e}(c)$.
We first suppose $\lambda_{j}-k$ is odd. Then $\xi_{j^{*}}^{j,\lambda_{j}-1-k}\in\mathfrak{g}^{e}(c-1)$
with $j\neq j^{*}$ and $\xi_{i}^{j^{*},0}+\varepsilon_{i,j,\lambda_{j}-1}\xi_{j}^{i^{*},1}\in\mathfrak{g}^{e}(1)$,
we have 
\[
[\xi_{j^{*}}^{j,\lambda_{j}-1-k},\xi_{i}^{j^{*},0}+\varepsilon_{i,j^{*},\lambda_{j}-1}\xi_{j}^{i^{*},1}]=\xi_{i}^{j,\lambda_{j}-1-k}-\varepsilon_{i,j^{*},\lambda_{j}-1}\xi_{j^{*}}^{i^{*},\lambda_{i}-1-k}.
\]
 By (\ref{eq:epsilon i,j,k}) we know that $\varepsilon_{i,j,k}=-\varepsilon_{i,j^{*},\lambda_{j}-1}$
as $\lambda_{j}-k$ is odd. Hence, 
\[
\xi_{i}^{j,\lambda_{j}-1-k}+\varepsilon_{i,j,k}\xi_{j^{*}}^{i^{*},\lambda_{i}-1-k}=[\xi_{j^{*}}^{j,\lambda_{j}-1-k},\xi_{i}^{j^{*},0}+\varepsilon_{i,j^{*},\lambda_{j}-1}\xi_{j}^{i^{*},1}]\in[\mathfrak{g}^{e}(c-1),\mathfrak{g}^{e}(1)].
\]
 Now suppose $\lambda_{j}-k$ is even, then $\xi_{j}^{j,\lambda_{j}-1-k}\in\mathfrak{g}^{e}(c-1)$
if $j=j^{*}$ or $\xi_{j}^{j,\lambda_{j}-1-k}+\xi_{j^{*}}^{j^{*},\lambda_{j}-1-k}\in\mathfrak{g}^{e}(c-1)$
if $j\neq j^{*}$. The argument for the first of these possibilities
is identical to the case when $\lambda_{j}-k$ is odd. We next consider
the second possibility. Note that 
\[
[\xi_{j}^{j,\lambda_{j}-1-k}+\xi_{j^{*}}^{j^{*},\lambda_{j}-1-k},\xi_{i}^{j,0}+\varepsilon_{i,j,\lambda_{j}-1}\xi_{j^{*}}^{i^{*},1}]=\xi_{i}^{j,\lambda_{j}-1-k}-\varepsilon_{i,j,\lambda_{j}-1}\xi_{j^{*}}^{i^{*},\lambda_{i}-1-k}.
\]
 By (\ref{eq:epsilon i,j,k}) we know that $\varepsilon_{i,j,k}=-\varepsilon_{i,j,\lambda_{j}-1}$
as $\lambda_{j}-k$ is even. Hence, 
\[
\xi_{i}^{j,\lambda_{j}-1-k}+\varepsilon_{i,j,k}\xi_{j^{*}}^{i^{*},\lambda_{i}-1-k}=[\xi_{j}^{j,\lambda_{j}-1-k}+\xi_{j^{*}}^{j^{*},\lambda_{j}-1-k},\xi_{i}^{j,0}+\varepsilon_{i,j,\lambda_{j}-1}\xi_{j^{*}}^{i^{*},1}]\in[\mathfrak{g}^{e}(c-1),\mathfrak{g}^{e}(1)].
\]

For $\lambda_{i}=\lambda_{j}$, let $2(\lambda_{i}-1-k)=c>1$. We
use induction on $\lambda_{i}$ to show that $\xi_{i}^{j,\lambda_{j}-1-k}+\varepsilon_{i,j,k}\xi_{j^{*}}^{i^{*},\lambda_{i}-1-k}$,
$\xi_{i}^{i^{*},\lambda_{i}-1-k}\in[\mathfrak{g}^{e}(c-1),\mathfrak{g}^{e}(1)]$.
Fix $i$ and assume for all $l,w,t$ and $\lambda_{l}<\lambda_{i}$,
$\xi_{l}^{w,\lambda_{w}-1-t}+\varepsilon_{l,w,t}\xi_{w^{*}}^{l^{*},\lambda_{l}-1-t}\in\mathfrak{g}^{e}(c)$
implies it belongs to $[\mathfrak{g}^{e}(c-1),\mathfrak{g}^{e}(1)]$.
Since $\lambda$ satisfies the condition in the theorem and $2(\lambda_{i}-1-k)=c>1$,
there exists $l$ such that $\lambda_{l}=1$. Note that the induction
starts with the case where $\lambda_{i}=\lambda_{j}=2$ and $k=0$,
then the induction base is dealt by calculating 
\[
\xi_{i}^{j,1}+\varepsilon_{i,j,0}\xi_{j^{*}}^{i^{*},1}=[\xi_{l}^{j,1}+\varepsilon_{l,j,0}\xi_{j^{*}}^{l^{*}.0},\xi_{i}^{l,1}+\varepsilon_{i,l,0}\xi_{l^{*}}^{i^{*},1}]\in[\mathfrak{g}^{e}(c-1),\mathfrak{g}^{e}(1)].
\]
 The induction step is as follows. Fix $i,j,k$ such that $2(\lambda_{i}-1-k)=c>1$.
Choose $l$ with $\lambda_{l}=\lambda_{i}-1$. Then for $\xi_{l}^{j,1}+\varepsilon_{i,j,\lambda_{l}-1}\xi_{j^{*}}^{l^{*},0}\in\mathfrak{g}^{e}(1)$,
$\xi_{i}^{l,\lambda_{l}-1-k}+\varepsilon_{i,l,k}\xi_{l^{*}}^{i^{*},\lambda_{i}-1-k}\in\mathfrak{g}^{e}(c-1)$,
we compute 
\begin{equation}
[\xi_{l}^{j,1}+\varepsilon_{i,j,\lambda_{l}-1}\xi_{j^{*}}^{l^{*},0},\xi_{i}^{l,\lambda_{l}-1-k}+\varepsilon_{i,l,k}\xi_{l^{*}}^{i^{*},\lambda_{i}-1-k}].\label{eq:proof-g^e(0)-g^e(1)}
\end{equation}

When $i\neq j$, the commutator in (\ref{eq:proof-g^e(0)-g^e(1)})
equals to $\xi_{i}^{j,\lambda_{j}-1-k}+\varepsilon_{i,j,k}\xi_{j^{*}}^{i^{*},\lambda_{i}-1-k}\in\mathfrak{g}^{e}(c)$.
If $i=j^{*}$, we have 
\[
\xi_{i}^{j,\lambda_{j}-1-k}+\varepsilon_{i,j,k}\xi_{j^{*}}^{i^{*},\lambda_{i}-1-k}=\begin{cases}
0 & \lambda_{i}-k\text{ even;}\\
2\xi_{i}^{i^{*},\lambda_{i}-1-k}\in\mathfrak{g}^{e}(c) & \lambda_{i}-k\text{ odd}.
\end{cases}
\]

Hence, we obtain that $\xi_{i}^{i^{*},\lambda_{i}-1-k}$ and $\xi_{i}^{j,\lambda_{j}-1-k}+\varepsilon_{i,j,k}\xi_{j^{*}}^{i^{*},\lambda_{i}-1-k}$
with $i\neq j,i\neq j^{*}$ lies in $[\mathfrak{g}^{e}(c-1),\mathfrak{g}^{e}(1)]$.

When $i=j$ and $i\neq i^{*}$, the commutator in (\ref{eq:proof-g^e(0)-g^e(1)})
equals to $(\xi_{i}^{i,\lambda_{i}-1-k}+\varepsilon_{i,i,k}\xi_{i^{*}}^{i^{*},\lambda_{i}-1-k})-(\xi_{l}^{l,\lambda_{i}-1-k}+\varepsilon_{l,l,k-1}\xi_{l^{*}}^{l^{*},\lambda_{i}-1-k})$.
If $i=i^{*}$, then
\[
\xi_{i}^{i,\lambda_{i}-1-k}+\varepsilon_{i,i,k}\xi_{i^{*}}^{i^{*},\lambda_{i}-1-k}=\begin{cases}
2\xi_{i}^{i,\lambda_{i}-1-k}\in\mathfrak{g}^{e}(c) & \lambda_{i}-k\text{ even};\\
0 & \lambda_{i}-k\text{ odd.}
\end{cases}
\]
By our inductive hypothesis we have $\xi_{l}^{l,\lambda_{i}-1-k}+\varepsilon_{l,l,k-1}\xi_{l^{*}}^{l^{*},\lambda_{i}-1-k}\in[\mathfrak{g}^{e}(c-1),\mathfrak{g}^{e}(1)]$.
Hence, we deduce that $\xi_{i}^{i,\lambda_{i}-1-k}+\varepsilon_{i,i,k}\xi_{i^{*}}^{i^{*},\lambda_{i}-1-k}$
with $i\neq i^{*}$ lies in $[\mathfrak{g}^{e}(c-1),\mathfrak{g}^{e}(1)]$.
The above argument goes through all basis elements in $\mathfrak{g}^{e}$.
Hence we obtain that $\mathfrak{g}^{e}$ is generated by $\mathfrak{g}^{e}(0)$
and $\mathfrak{g}^{e}(1)$, which is equivalent to $\mathfrak{g}^{e}(\geq1)$
is generated by $\mathfrak{g}^{e}(1)$ as a Lie superalgebra. 
\end{proof}

\section{On reachable nilpotent elements in exceptional Lie superalgebras\label{sec:exceptional}}

\noindent Let $\mathfrak{g}=\mathfrak{g}_{\bar{0}}\oplus\mathfrak{g}_{\bar{1}}=D(2,1;\alpha)$,
$G(3)$ or $F(4)$ and let $e\in\mathfrak{g}_{\bar{0}}$ be nilpotent.
Let $\{e,h,f\}\subseteq\mathfrak{g}_{\bar{0}}$ be an $\mathfrak{sl}(2)$-triple
containing $e$. In this section, we give the classification of $e\in\mathfrak{g}_{\bar{0}}$
which are reachable, strongly reachable or satisfy the Panyushev property.
We start each subsection by recalling the notation of basis elements
for $\mathfrak{g}$ and the commutators of basis elements as given
in \cite[Sections 4--6]{han-exc}. In this paper we adopt the notation
in \cite{han-exc}. Throughout the following subsections, we consider
the Lie algebra $\mathfrak{sl}(2)$ frequently. Hence, let us fix
the notation $\mathfrak{sl}(2)=\left\langle E,H,F\right\rangle $
where 
\[
E=\begin{pmatrix}0 & 1\\
0 & 0
\end{pmatrix},H=\begin{pmatrix}1 & 0\\
0 & -1
\end{pmatrix},F=\begin{pmatrix}0 & 0\\
1 & 0
\end{pmatrix}.
\]

\subsection{Reachability, strong reachability and the Panyushev property for
$D(2,1;\alpha)$}

\noindent Let $\mathfrak{g}=\mathfrak{g}_{\bar{0}}\oplus\mathfrak{g}_{\bar{1}}=D(2,1;\alpha)$.
Note that $\mathfrak{g}$ sometimes is denoted by $\Gamma(\sigma_{1},\sigma_{2},\sigma_{3})$
where $\sigma_{1},\sigma_{2},\sigma_{3}$ are complex numbers such
that $\sigma_{1}+\sigma_{2}+\sigma_{3}=0$. For any $\alpha\in\mathbb{C}\setminus\{0,-1\}$,
we have $D(2,1;\alpha)\cong\Gamma(1+\alpha,-1,-\alpha)\cong\Gamma(\frac{1+\alpha}{\alpha},-1,-\frac{1}{\alpha})\cong\Gamma(-\alpha,-1,1+\alpha)$
and the supercommutator $\mathfrak{g}_{\bar{0}}\times\mathfrak{g}_{\bar{1}}\rightarrow\mathfrak{g}_{\bar{1}}$
depends on $\sigma_{i}$. The even part of $\mathfrak{g}$ is $\mathfrak{g}_{\bar{0}}=\mathfrak{sl}(2)\oplus\mathfrak{sl}(2)\oplus\mathfrak{sl}(2)$
and the odd part of $\mathfrak{g}$ is $\mathfrak{g}_{\bar{1}}=V_{1}\otimes V_{2}\otimes V_{3}$.
Note that each $V_{i}$ is a copy of a two-dimensional vector space
$V$ with basis elements $v_{1}=(1,0)^{t}$ and $v_{-1}=(0,1)^{t}$.
In order to give a basis for $\mathfrak{g}$, we denote $E_{1}=(E,0,0)$,
$E_{2}=(0,E,0)$ and $E_{3}=(0,0,E)$. Similarly, we denote $F_{1}=(F,0,0)$,
$F_{2}=(0,F,0)$, $F_{3}=(0,0,F)$, $H_{1}=(H,0,0)$, $H_{2}=(0,H,0)$
and $H_{3}=(0,0,H)$. By \cite[Section 4.1]{han-exc}, a basis for
$\mathfrak{g}_{\bar{0}}$ is $\{E_{i},H_{i},F_{i}:i=1,2,3\}$ and
a basis for $\mathfrak{g}_{\bar{1}}$ is $\{v_{i}\otimes v_{j}\otimes v_{k}:i,j,k=\pm1\}$.

Representatives $e$ of nilpotent orbits in $\mathfrak{g}_{\bar{0}}$
are given in \cite[Section 4.3]{han-exc}, these are $e=0$, $E_{1}$,
$E_{2}$, $E_{3}$, $E_{1}+E_{2}$, $E_{1}+E_{3}$, $E_{2}+E_{3}$,
$E_{1}+E_{2}+E_{3}$. In Table \ref{tab:Reachable D(2,1;)}, we classify
which ones are reachable, strongly reachable or satisfy the Panyushev
property based on the structure of $\mathfrak{g}^{e}$ in \cite[Table 4.3]{han-exc}.

\begin{table}[H]
\noindent \begin{centering}
\begin{tabular}{|>{\centering}p{5cm}|c|c|>{\centering}p{4cm}|}
\hline 
Representatives of nilpotent orbits $e\in\mathfrak{g}_{\bar{0}}$ & Reachable & Strongly reachable & Satisfying the Panyushev property\tabularnewline
\hline 
\hline 
$0$ & $\checkmark$ & $\checkmark$ & $\checkmark$\tabularnewline
\hline 
$E_{1},E_{2},E_{3}$ & $\checkmark$ & $\checkmark$ & $\checkmark$\tabularnewline
\hline 
$E_{1}+E_{2},E_{1}+E_{3},E_{2}+E_{3}$ &  &  & \tabularnewline
\hline 
$E_{1}+E_{2}+E_{3}$ & $\checkmark$ &  & $\checkmark$\tabularnewline
\hline 
\end{tabular}
\par\end{centering}
\caption{\label{tab:Reachable D(2,1;)}Reachable, strongly reachable and Panyushev
elements in $D(2,1;\alpha)$}
\end{table}

In the remaining part of this subsection, we explain our calculations
explicitly for $e=0$, $E_{1}$, $E_{1}+E_{2}$, $E_{1}+E_{2}+E_{3}$.
Note that analysis for cases $e=E_{2}$, $E_{3}$ are similar to $e=E_{1}$
and analysis for the cases $e=E_{1}+E_{3}$, $E_{2}+E_{3}$ are similar
to $e=E_{1}+E_{2}$, for which we omitted in this paper.

\begin{onehalfspace}
\noindent (1) $e=0$
\end{onehalfspace}

\noindent It is clear that $e$ is reachable and satisfies the Panyushev
property as $\mathfrak{g}^{e}=\mathfrak{g}$ and $\mathfrak{g}^{e}(\geq1)=\mathfrak{g}^{e}(1)=0$.

\begin{onehalfspace}
\noindent (2) $e=E_{1}$
\end{onehalfspace}

\noindent The semisimple element $h=H_{1}$. By \cite[Table 4.3]{han-exc},
$\mathfrak{g}^{e}=\mathfrak{g}^{e}(0)\oplus\mathfrak{g}^{e}(1)\oplus\mathfrak{g}^{e}(2)$
where $\mathfrak{g}^{e}(0)=\left\langle E_{2},H_{2},F_{2},E_{3},H_{3},F_{3}\right\rangle $,
$\mathfrak{g}^{e}(1)=\left\langle v_{1}\otimes v_{j}\otimes v_{k}:j,k=\pm1\right\rangle $
and $\mathfrak{g}^{e}(2)=\left\langle E_{1}\right\rangle $. Note
that $[v_{1}\otimes v_{1}\otimes v_{1},v_{1}\otimes v_{-1}\otimes v_{-1}]=[v_{1}\otimes v_{1}\otimes v_{-1},v_{1}\otimes v_{-1}\otimes v_{1}]=2\sigma_{1}E_{1}\in\mathfrak{g}^{e}(2)$,
hence $\mathfrak{g}^{e}(\geq1)$ is generated by $\mathfrak{g}^{e}(1)$.
This also implies that $e$ is reachable. In order to prove the strong
reachability, it suffices to show that $\mathfrak{g}^{e}(0)\oplus\mathfrak{g}^{e}(1)\subseteq[\mathfrak{g}^{e},\mathfrak{g}^{e}]$.
It is clear that $\mathfrak{g}^{e}(0)\subseteq[\mathfrak{g}^{e},\mathfrak{g}^{e}]$
as $\{E_{2},H_{2},F_{2}\}$ and $\{E_{3},H_{3},F_{3}\}$ are $\mathfrak{sl}(2)$-triples.
Observe that $[H_{2},v_{1}\otimes v_{j}\otimes v_{k}]=jv_{1}\otimes v_{j}\otimes v_{k}\in\mathfrak{g}^{e}(1)$,
thus $\mathfrak{g}^{e}(1)\subseteq[\mathfrak{g}^{e},\mathfrak{g}^{e}]$.
Therefore, we have that $\mathfrak{g}^{e}=[\mathfrak{g}^{e},\mathfrak{g}^{e}]$.

\begin{onehalfspace}
\noindent (3) $e=E_{1}+E_{2}$
\end{onehalfspace}

\noindent The semisimple element $h=H_{1}+H_{2}$. By \cite[Table 4.3]{han-exc},
$\mathfrak{g}^{e}=\mathfrak{g}^{e}(0)\oplus\mathfrak{g}^{e}(2)$ where
$\mathfrak{g}^{e}(0)=\langle E_{3},H_{3},F_{3},v_{1}\otimes v_{-1}\otimes v_{1}-v_{-1}\otimes v_{1}\otimes v_{1},v_{1}\otimes v_{-1}\otimes v_{-1}-v_{-1}\otimes v_{1}\otimes v_{-1}\rangle$
and $\mathfrak{g}^{e}(2)=\langle E_{1},E_{2},v_{1}\otimes v_{1}\otimes v_{1},v_{1}\otimes v_{1}\otimes v_{-1}\rangle$.
It is clear that $e$ does not satisfy the Panyushev property as $\mathfrak{g}^{e}(1)=0$.
In this case $e\in[\mathfrak{g}^{e},\mathfrak{g}^{e}]$ implies $e\in[\mathfrak{g}^{e}(0),\mathfrak{g}^{e}(2)]$.
However, $[v_{1}\otimes v_{1}\otimes v_{-1},v_{1}\otimes v_{-1}\otimes v_{1}-v_{-1}\otimes v_{1}\otimes v_{1}]=-[v_{1}\otimes v_{1}\otimes v_{1},v_{1}\otimes v_{-1}\otimes v_{-1}-v_{-1}\otimes v_{1}\otimes v_{-1}]=-2\sigma_{1}E_{1}+2\sigma_{2}E_{2}$
and commutators between the other pairs of $[x,y]$ such that $x\in\mathfrak{g}^{e}(0),y\in\mathfrak{g}^{e}(2)]$
are all equal to zero. Hence, we deduce that $e\notin[\mathfrak{g}^{e},\mathfrak{g}^{e}]$.

\noindent (4) $e=E_{1}+E_{2}+E_{3}$

\noindent The semisimple element $h=H_{1}+H_{2}+H_{3}$. By \cite[Table 4.3]{han-exc},
$\mathfrak{g}^{e}=\mathfrak{g}^{e}(1)\oplus\mathfrak{g}^{e}(2)\oplus\mathfrak{g}^{e}(3)$
where $\mathfrak{g}^{e}(1)=\langle v_{1}\otimes v_{1}\otimes v_{-1}-v_{-1}\otimes v_{1}\otimes v_{1},v_{1}\otimes v_{-1}\otimes v_{1}-v_{-1}\otimes v_{1}\otimes v_{1}\rangle,$
$\mathfrak{g}^{e}(2)=\left\langle E_{1},E_{2},E_{3}\right\rangle $
and $\mathfrak{g}^{e}(3)=\left\langle v_{1}\otimes v_{1}\otimes v_{1}\right\rangle $.
Denote by $x=v_{1}\otimes v_{1}\otimes v_{-1}-v_{-1}\otimes v_{1}\otimes v_{1}$
and $y=v_{1}\otimes v_{-1}\otimes v_{1}-v_{-1}\otimes v_{1}\otimes v_{1}$.
Then $[x,x]=4\sigma_{2}E_{2}$, $[y,y]=4\sigma_{3}E_{3}$, $[x,y]=-2\sigma_{1}E_{1}+2\sigma_{2}E_{2}+2\sigma_{3}E_{3}$
and $[E_{2},y]=v_{1}\otimes v_{1}\otimes v_{1}$. Hence, we conclude
that $\mathfrak{g}^{e}(\geq1)$ is generated by $\mathfrak{g}^{e}(1)$.
We say that $e$ is reachable because $e$ can be written as a linear
combination of $[x,x]$, $[y,y]$ and $[x,y]$. However, $e$ is not
strongly reachable as $\mathfrak{g}^{e}(1)\nsubseteq[\mathfrak{g}^{e},\mathfrak{g}^{e}]$.

\subsection{Reachability, strong reachability and the Panyushev property for
$G(3)$}

Let $\mathfrak{g}=\mathfrak{g}_{\bar{0}}\oplus\mathfrak{g}_{\bar{1}}=G(3)$.
According to \cite[Section 5]{han-exc}, $\mathfrak{g}_{\bar{0}}=\mathfrak{sl}(2)\oplus G_{2}$
and $\mathfrak{g}_{\bar{1}}=V_{2}\otimes V_{7}$ where $V_{2}$ is
a two-dimensional vector space with basis elements $v_{1}=(1,0)^{t}$
and $v_{-1}=(0,1)^{t}$ and $V_{7}$ is a seven-dimensional simple
representation of the Lie algebra $G_{2}$. A basis of $V_{7}$ is
denoted by $\{e_{3},e_{2},e_{1},e_{0},e_{-1},e_{-2},e_{-3}\}$, see
\cite[Chapter 4]{Musson2012}. Note that $G_{2}$ can be viewed as
a Lie subalgebra of $\mathfrak{gl}(V_{7})$, thus we write the elements
of $G_{2}$ with respect to the basis of $V_{7}$ such that $G_{2}$
has a basis $\{h_{1},h_{2},x_{i},y_{i}:i=1,\dots,6\}$ where $h_{1}=\text{diag}(1,-1,2,0,-2,1,-1)$,
$h_{2}=\text{diag}(0,1,-1,0,1,-1,0)$, 
\[
x_{1}=\begin{pmatrix}0 & -1 & 0 & 0 & 0 & 0 & 0\\
0 & 0 & 0 & 0 & 0 & 0 & 0\\
0 & 0 & 0 & 1 & 0 & 0 & 0\\
0 & 0 & 0 & 0 & -2 & 0 & 0\\
0 & 0 & 0 & 0 & 0 & 0 & 0\\
0 & 0 & 0 & 0 & 0 & 0 & 1\\
0 & 0 & 0 & 0 & 0 & 0 & 0
\end{pmatrix},\ x_{2}=\begin{pmatrix}0 & 0 & 0 & 0 & 0 & 0 & 0\\
0 & 0 & 1 & 0 & 0 & 0 & 0\\
0 & 0 & 0 & 0 & 0 & 0 & 0\\
0 & 0 & 0 & 0 & 0 & 0 & 0\\
0 & 0 & 0 & 0 & 0 & -1 & 0\\
0 & 0 & 0 & 0 & 0 & 0 & 0\\
0 & 0 & 0 & 0 & 0 & 0 & 0
\end{pmatrix},
\]
\[
y_{1}=\begin{pmatrix}0 & 0 & 0 & 0 & 0 & 0 & 0\\
-1 & 0 & 0 & 0 & 0 & 0 & 0\\
0 & 0 & 0 & 0 & 0 & 0 & 0\\
0 & 0 & 2 & 0 & 0 & 0 & 0\\
0 & 0 & 0 & -1 & 0 & 0 & 0\\
0 & 0 & 0 & 0 & 0 & 0 & 0\\
0 & 0 & 0 & 0 & 0 & 1 & 0
\end{pmatrix},\ y_{2}=\begin{pmatrix}0 & 0 & 0 & 0 & 0 & 0 & 0\\
0 & 0 & 0 & 0 & 0 & 0 & 0\\
0 & 1 & 0 & 0 & 0 & 0 & 0\\
0 & 0 & 0 & 0 & 0 & 0 & 0\\
0 & 0 & 0 & 0 & 0 & 0 & 0\\
0 & 0 & 0 & 0 & -1 & 0 & 0\\
0 & 0 & 0 & 0 & 0 & 0 & 0
\end{pmatrix},
\]
 $x_{3}=[x_{1},x_{2}]$, $x_{4}=[x_{1},x_{3}]$, $x_{5}=[x_{1},x_{4}]$
and $x_{6}=[x_{5},x_{2}]$. The remaining $y_{i}$ for $i=3,\dots,6$
can be generated by $y_{1}$ and $y_{2}$ in a similar way. A basis
for $\mathfrak{g}_{\bar{1}}$ is $\{v_{i}\otimes e_{j}:i=\pm1,j=0,\pm1,\pm2,\pm3\}$.

Representatives $e$ of nilpotent orbits in $\mathfrak{g}_{\bar{0}}$
are given in \cite[Section 5.3]{han-exc} which we list in Table \ref{tab:Reachable G(3)}.
In addition, we determine which $e$ are reachable, strongly reachable
or satisfy the Panyushev property in Table \ref{tab:Reachable G(3)}.

\begin{table}[H]
\noindent \begin{centering}
\begin{tabular}{|>{\centering}m{4cm}|>{\centering}p{2cm}|>{\centering}p{2cm}|>{\centering}m{4cm}|}
\hline 
Representatives of nilpotent orbits $e\in\mathfrak{g}_{\bar{0}}$ & Reachable & Strongly reachable & Satisfying the Panyushev property\tabularnewline
\hline 
\hline 
$E+(x_{1}+x_{2})$ &  &  & \tabularnewline
\hline 
$E+x_{2}$ & $\checkmark$ & $\checkmark$ & $\checkmark$\tabularnewline
\hline 
$E+x_{1}$ & $\checkmark$ &  & \tabularnewline
\hline 
$E+(x_{2}+x_{5})$ & $\checkmark$ &  & $\checkmark$\tabularnewline
\hline 
$E$ & $\checkmark$ & $\checkmark$ & $\checkmark$\tabularnewline
\hline 
$x_{1}+x_{2}$ &  &  & \tabularnewline
\hline 
$x_{2}$ & $\checkmark$ & $\checkmark$ & $\checkmark$\tabularnewline
\hline 
$x_{1}$ & $\checkmark$ & $\checkmark$ & \tabularnewline
\hline 
$x_{2}+x_{5}$ &  &  & \tabularnewline
\hline 
$0$ & $\checkmark$ & $\checkmark$ & $\checkmark$\tabularnewline
\hline 
\end{tabular}
\par\end{centering}
\caption{\label{tab:Reachable G(3)}Reachable, strongly reachable and Panyushev
elements in $G(3)$}

\end{table}

In the remaining part of this subsection, we take $e=E+x_{2}$ and
$e=x_{1}$ as examples to show our explicit calculations. Other cases
are being dealt with using a similar approach.

\begin{onehalfspace}
\noindent (1) $e=E+x_{2}$
\end{onehalfspace}

\noindent The semisimple element $h=H+h_{2}$. Note that $\mathfrak{g}^{e}=\mathfrak{g}^{e}(0)\oplus\mathfrak{g}^{e}(1)\oplus\mathfrak{g}^{e}(2)$
where basis elements are shown in the following table.

\noindent 
\begin{table}[H]
\noindent \centering{}%
\begin{tabular}{|c|>{\centering}m{10cm}|}
\hline 
$\mathfrak{g}^{e}(0)$ & $\left\langle x_{4},y_{4},2h_{1}+3h_{2},v_{1}\otimes e_{1}-v_{-1}\otimes e_{2},v_{1}\otimes e_{-2}+v_{-1}\otimes e_{-1}\right\rangle $\tabularnewline
\hline 
$\mathfrak{g}^{e}(1)$ & $\left\langle y_{1},x_{3},x_{6},y_{5},v_{1}\otimes e_{3},v_{1}\otimes e_{0},v_{1}\otimes e_{-3}\right\rangle $\tabularnewline
\hline 
$\mathfrak{g}^{e}(2)$ & $\left\langle E,x_{2},v_{1}\otimes e_{2},v_{1}\otimes e_{-1}\right\rangle $\tabularnewline
\hline 
\end{tabular}\caption{$\mathrm{ad}h$-eigenspaces for $\mathfrak{g}^{e}$ when $e=E+x_{2}$}
\end{table}
 By calculating commutators of elements in $\mathfrak{g}^{e}(1)$,
we have that $[y_{1},v_{1}\otimes e_{3}]=-v_{1}\otimes e_{2}\in\mathfrak{g}^{e}(2)$,
$[y_{1},v_{1}\otimes e_{0}]=-v_{1}\otimes e_{-1}\in\mathfrak{g}^{e}(2)$,
$[y_{1},x_{3}]=3x_{2}\in\mathfrak{g}^{e}(2)$ and $[v_{1}\otimes e_{3},v_{1}\otimes e_{-3}]=16E\in\mathfrak{g}^{e}(2)$.
Hence, we have that $\mathfrak{g}^{e}(\geq1)$ is generated by $\mathfrak{g}^{e}(1)$.
Since $e$ can be written as a linear combination of $[v_{1}\otimes e_{3},v_{1}\otimes e_{-3}]$
and $[y_{1},x_{3}]$, we say that $e$ is reachable. 

To show that $e$ is strongly reachable, it remains to show that $\mathfrak{g}^{e}(0)\oplus\mathfrak{g}^{e}(1)\subseteq[\mathfrak{g}^{e},\mathfrak{g}^{e}]$.
Denote by $x=v_{1}\otimes e_{1}-v_{-1}\otimes e_{2}$ and $y=v_{1}\otimes e_{-2}+v_{-1}\otimes e_{-1}$.
Note that $\mathfrak{g}^{e}(0)\cong\mathfrak{sl}(2)\oplus\left\langle x,y\right\rangle $
and $[2h_{1}+3h_{2},x]=x$, $[2h_{1}+3h_{2},y]=-y$. This implies
that $\mathfrak{g}^{e}(0)\subseteq[\mathfrak{g}^{e},\mathfrak{g}^{e}]$.
Next we look at $\mathfrak{g}^{e}(1)$. We know that $[x,v_{1}\otimes e_{-3}]=-4y_{1}$,
$[x,v_{1}\otimes e_{0}]=-4x_{3}$, $[x,v_{1}\otimes e_{3}]=2x_{6}$,
$[y,v_{1}\otimes e_{-3}]=2y_{5}$, $[x_{4},v_{1}\otimes e_{0}]=2v_{1}\otimes e_{3}$,
$[x_{4},v_{1}\otimes e_{-3}]=-4v_{1}\otimes e_{0}$ and $[y_{4},v_{1}\otimes e_{0}]=-2v_{1}\otimes e_{-3}$.
This provides us every basis elements of $\mathfrak{g}^{e}(1)$, thus
$\mathfrak{g}^{e}(1)\subseteq[\mathfrak{g}^{e},\mathfrak{g}^{e}]$.
Therefore, we deduce that $\mathfrak{g}^{e}=[\mathfrak{g}^{e},\mathfrak{g}^{e}]$.

\begin{onehalfspace}
\noindent (2) $e=x_{1}$
\end{onehalfspace}

\noindent The semisimple element $h=h_{1}$. The centralizer $\mathfrak{g}^{e}=\mathfrak{g}^{e}(0)\oplus\mathfrak{g}^{e}(1)\oplus\mathfrak{g}^{e}(2)\oplus\mathfrak{g}^{e}(3)$
where basis elements are shown in the following table.

\noindent 
\begin{table}[H]
\noindent \centering{}%
\begin{tabular}{|c|>{\centering}m{10cm}|}
\hline 
$\mathfrak{g}^{e}(0)$ & $\left\langle E,H,F,x_{6},y_{6},h_{1}+2h_{2}\right\rangle $\tabularnewline
\hline 
$\mathfrak{g}^{e}(1)$ & $\left\langle v_{i}\otimes e_{j}:i=\pm1,j=3,-2\right\rangle $\tabularnewline
\hline 
$\mathfrak{g}^{e}(2)$ & $\left\langle x_{1},v_{i}\otimes e_{1}:i=\pm1\right\rangle $\tabularnewline
\hline 
$\mathfrak{g}^{e}(3)$ & \noindent $\left\langle x_{5},y_{2}\right\rangle $\tabularnewline
\hline 
\end{tabular}\caption{$\mathrm{ad}h$-eigenspaces for $\mathfrak{g}^{e}$ when $e=x_{1}$}
\end{table}

\noindent We have that $[v_{1}\otimes e_{3},v_{-1}\otimes e_{-2}]=-[v_{-1}\otimes e_{3},v_{1}\otimes e_{-2}]=-4x_{1}$
and other commutators between basis elements for $\mathfrak{g}^{e}(1)$
are zero. This implies that $e$ is reachable but $\mathfrak{g}^{e}(\geq1)$
is not generated by $\mathfrak{g}^{e}(1)$ as we cannot obtain other
basis elements of $\mathfrak{g}^{e}(2)$ by computing commutators
between basis elements of $\mathfrak{g}^{e}(1)$. 

We next show that $e$ is indeed strongly reachable. It is clear that
$\mathfrak{g}^{e}(0)\subseteq[\mathfrak{g}^{e},\mathfrak{g}^{e}]$
as $\mathfrak{g}^{e}(0)\cong\mathfrak{sl}(2)\oplus\mathfrak{sl}(2)$.
Note that $[H,v_{i}\otimes e_{j}]=iv_{i}\otimes e_{j}$, then $\mathfrak{g}^{e}(1)\oplus\mathfrak{g}^{e}(2)\subseteq[\mathfrak{g}^{e},\mathfrak{g}^{e}]$.
Moreover, $\mathfrak{g}^{e}(3)\subseteq[\mathfrak{g}^{e},\mathfrak{g}^{e}]$
since $[v_{1}\otimes e_{3},v_{-1}\otimes e_{1}]=2x_{5}$ and $[v_{1}\otimes e_{-2},v_{-1}\otimes e_{1}]=-12y_{2}$.
Therefore, we obtain that $e$ is strongly reachable.

\subsection{Reachability, strong reachability and the Panyushev property for
$F(4)$}

Let $\mathfrak{g}=\mathfrak{g}_{\bar{0}}\oplus\mathfrak{g}_{\bar{1}}=F(4)$
such that $\mathfrak{g}_{\bar{0}}=\mathfrak{sl}(2)\oplus\mathfrak{so}(7)$
and $\mathfrak{g}_{\bar{1}}=V_{2}\otimes V_{8}$ where $V_{2}$ is
a two-dimensional vector space with basis elements $v_{1}=(1,0)^{t}$
and $v_{-1}=(0,1)^{t}$ and $V_{8}$ is the spin representation for
$\mathfrak{so}(7)$. Let $V_{\mathfrak{so}}$ be a vector space such
that $\mathfrak{so}(V_{\mathfrak{so}})=\mathfrak{so}(7)$, we use
the notation from \cite[Section 6.1]{han-exc} for a basis for $V_{\mathfrak{so}}$
where $V_{\mathfrak{so}}=\langle e_{3},e_{2},e_{1},e_{0},e_{-1},e_{-2},e_{-3}\rangle$.
Then a basis for $\mathfrak{so}(7)$ consists of the following elements
\[
R_{e_{i},e_{-j}}=e_{i,j}-e_{-j,-i}\text{ and }R_{e_{i},e_{0}}=2e_{i,0}-e_{0,-i}\text{ for }i,j\in\{\pm1,\pm2,\pm3\}
\]
 such that $e_{i,j}$ is the elementary transformation which sends
$e_{i}$ to $e_{j}$ and the other basis vectors to $0$. A basis
for $V_{8}$ is denoted by $\{s,e_{1}s,e_{2}s,e_{3}s,e_{1}e_{2}s,e_{1}e_{3}s,e_{2}e_{3}s,e_{1}e_{2}e_{3}s\}$
where $e_{-i}s=0$ for $i=1,2,3$ and $e_{0}s=s$.

Representatives $e$ of nilpotent orbits in $\mathfrak{g}_{\bar{0}}$
are given in \cite[Section 6.5]{han-exc}. We list these representatives
in Table \ref{tab:Reachable F(4)} and state which $e$ are reachable,
strongly reachable or satisfy the Panyushev property.

\begin{table}[H]
\noindent \begin{centering}
\begin{tabular}{|>{\centering}p{5cm}|>{\centering}m{2cm}|>{\centering}m{2cm}|>{\centering}m{4cm}|}
\hline 
Representatives of nilpotent orbits $e\in\mathfrak{g}_{\bar{0}}$ & Reachable & Strongly reachable & Satisfying the Panyushev property\tabularnewline
\hline 
\hline 
$E+(R_{e_{1},e_{-2}}+R_{e_{2},e_{-3}}+R_{e_{3},e_{0}})$ &  &  & \tabularnewline
\hline 
$E+(R_{e_{1},e_{-2}}+R_{e_{2},e_{0}})$ &  &  & \tabularnewline
\hline 
$E+(R_{e_{1},e_{-3}}+R_{e_{2},e_{3}})$ & $\checkmark$ &  & $\checkmark$\tabularnewline
\hline 
$E+(R_{e_{1},e_{0}}+R_{e_{2},e_{3}})$ & $\checkmark$ &  & $\checkmark$\tabularnewline
\hline 
$E+R_{e_{1},e_{0}}$ &  &  & \tabularnewline
\hline 
$E+R_{e_{1},e_{2}}$ & $\checkmark$ & $\checkmark$ & $\checkmark$\tabularnewline
\hline 
$E$ & $\checkmark$ & $\checkmark$ & $\checkmark$\tabularnewline
\hline 
$R_{e_{1},e_{-2}}+R_{e_{2},e_{-3}}+R_{e_{3},e_{0}}$ &  &  & \tabularnewline
\hline 
$R_{e_{1},e_{-2}}+R_{e_{2},e_{0}}$ &  &  & \tabularnewline
\hline 
$R_{e_{1},e_{-3}}+R_{e_{2},e_{3}}$ &  &  & \tabularnewline
\hline 
$R_{e_{1},e_{0}}+R_{e_{2},e_{3}}$ & $\checkmark$ & $\checkmark$ & $\checkmark$\tabularnewline
\hline 
$R_{e_{1},e_{0}}$ & $\checkmark$ & $\checkmark$ & $\checkmark$\tabularnewline
\hline 
$R_{e_{1},e_{2}}$ & $\checkmark$ & $\checkmark$ & $\checkmark$\tabularnewline
\hline 
$0$ & $\checkmark$ & $\checkmark$ & $\checkmark$\tabularnewline
\hline 
\end{tabular}
\par\end{centering}
\caption{\label{tab:Reachable F(4)}Reachable, strongly reachable and Panyushev
elements in $F(4)$}

\end{table}

In the remaining part of this subsection, we demonstrate details of
our calculation for cases $e=E+(R_{e_{1},e_{-3}}+R_{e_{2},e_{3}})$
and $e=E+R_{e_{1},e_{2}}$ . The method for general cases is similar
to these two cases. 

\noindent (1) $e=E+(R_{e_{1},e_{-3}}+R_{e_{2},e_{3}})$

\noindent The semisimple element $h=H+(2R_{e_{1},e_{-1}}+2R_{e_{2},e_{-2}})$.
The centralizer $\mathfrak{g}^{e}=\mathfrak{g}^{e}(0)\oplus\mathfrak{g}^{e}(1)\oplus\mathfrak{g}^{e}(2)\oplus\mathfrak{g}^{e}(3)\oplus\mathfrak{g}^{e}(4)$
where basis elements are shown in the following table. 

\noindent 
\begin{table}[H]
\noindent \begin{centering}
\begin{tabular}{|c|>{\centering}m{11cm}|}
\hline 
$\mathfrak{g}^{e}(0)$ & $\left\langle R_{e_{1},e_{-1}}-R_{e_{2},e_{-2}}+R_{e_{3},3}\right\rangle $\tabularnewline
\hline 
$\mathfrak{g}^{e}(1)$ & $\left\langle v_{1}\otimes e_{1}s-v_{-1}\otimes e_{1}e_{2}e_{3}s,v_{1}\otimes e_{2}s,v_{-1}\otimes e_{1}e_{2}s+v_{1}\otimes e_{2}e_{3}s,v_{1}\otimes e_{1}e_{3}s\right\rangle $\tabularnewline
\hline 
$\mathfrak{g}^{e}(2)$ & $\left\langle E,R_{e_{1},e_{-3}}+R_{e_{2},e_{3}},R_{e_{1},e_{-3}},R_{e_{2},e_{-3}},R_{e_{2},e_{0}},R_{e_{1},e_{0}},R_{e_{1},e_{3}}\right\rangle $\tabularnewline
\hline 
$\mathfrak{g}^{e}(3)$ & $\left\langle v_{1}\otimes e_{1}e_{2}e_{3}s,v_{1}\otimes e_{1}e_{2}s\right\rangle $\tabularnewline
\hline 
$\mathfrak{g}^{e}(4)$ & $\left\langle R_{e_{1},e_{2}}\right\rangle $\tabularnewline
\hline 
\end{tabular}
\par\end{centering}
\caption{$\mathrm{ad}h$-eigenspaces for $\mathfrak{g}^{e}$ when $e=E+(R_{e_{1},e_{-3}}+R_{e_{2},e_{3}})$}

\end{table}
 Let $x=v_{1}\otimes e_{1}s-v_{-1}\otimes e_{1}e_{2}e_{3}s$, $y=v_{-1}\otimes e_{1}e_{2}s+v_{1}\otimes e_{2}e_{3}s\in\mathfrak{g}^{e}(1)$.
Then we compute $[x,x]=R_{e_{1},e_{0}}$, $[x,v_{1}\otimes e_{2}s]=\frac{1}{2}R_{e_{2},e_{0}}$,
$[x,v_{1}\otimes e_{1}e_{3}s]=R_{e_{1},e_{3}}$, $[x,y]=R_{e_{1},e_{-3}}+R_{e_{2},e_{3}}-6E$,
$[v_{1}\otimes e_{2}s,y]=R_{e_{2},e_{-3}}$ and $[v_{1}\otimes e_{2}s,v_{1}\otimes e_{1}e_{3}s]=6E$.
The above commutators give us all basis elements for $\mathfrak{g}^{e}(2)$.
Similarly, we have that $[R_{e_{1},e_{0}},v_{1}\otimes e_{2}s]=-v_{1}\otimes e_{1}e_{2}s\in\mathfrak{g}^{e}(3)$,
$[R_{e_{1},e_{0}},y]=v_{1}\otimes e_{1}e_{2}e_{3}s\in\mathfrak{g}^{e}(3)$
and $[v_{1}\otimes e_{1}e_{2}e_{3}s,y]=R_{e_{1},e_{2}}\in\mathfrak{g}^{e}(4)$.
Hence, we conclude that $\mathfrak{g}^{e}(\geq1)$ is generated by
$\mathfrak{g}^{e}(1)$. 

The above argument also implies that $e\in[\mathfrak{g}^{e},\mathfrak{g}^{e}]$
as $e$ can be written as a linear combination of $[x,y]$ and $[v_{1}\otimes e_{2}s,v_{1}\otimes e_{1}e_{3}s]$.
Therefore, we have that $e$ is reachable. Note that $e$ is not strongly
reachable because basis elements in $\mathfrak{g}^{e}(0)$ cannot
be written as a sum of commutators between basis elements in $\mathfrak{g}^{e}$.

\noindent (2) $e=E+R_{e_{1},e_{2}}$

\noindent The semisimple element $h=R_{e_{1},e_{-1}}+R_{e_{2},e_{-2}}$.
The centralizer $\mathfrak{g}^{e}=\mathfrak{g}^{e}(0)\oplus\mathfrak{g}^{e}(1)\oplus\mathfrak{g}^{e}(2)$
where basis elements are shown in the following table. 

\begin{table}[H]
\noindent \begin{centering}
\begin{tabular}{|c|>{\centering}m{14cm}|}
\hline 
$\mathfrak{g}^{e}(0)$ & \noindent $\langle R_{e_{1},e_{-2}},R_{e_{1},e_{-1}}-R_{e_{2},e_{-2}},R_{e_{3},e_{-3}},R_{e_{2},e_{-1}},R_{e_{-3},e_{0}},R_{e_{3},e_{0}},v_{1}\otimes s-v_{-1}\otimes e_{1}e_{2}s,v_{1}\otimes e_{3}s-v_{-1}\otimes e_{1}e_{2}e_{3}s\rangle$\tabularnewline
\hline 
$\mathfrak{g}^{e}(1)$ & $\left\langle R_{e_{1},e_{3}},R_{e_{1},e_{0}},R_{e_{2},e_{-3}},R_{e_{2},e_{3}},R_{e_{1},e_{-3}},R_{e_{2},e_{0}},v_{1}\otimes e_{1}s,v_{1}\otimes e_{2}s,v_{1}\otimes e_{1}e_{3}s,v_{1}\otimes e_{2}e_{3}s\right\rangle $\tabularnewline
\hline 
$\mathfrak{g}^{e}(2)$ & $\left\langle E,R_{e_{1},e_{2}},v_{1}\otimes e_{1}e_{2}s,v_{1}\otimes e_{1}e_{2}e_{3}s\right\rangle $\tabularnewline
\hline 
\end{tabular}\caption{$\mathrm{ad}h$-eigenspaces for $\mathfrak{g}^{e}$ when $e=E+R_{e_{1},e_{2}}$}
\par\end{centering}
\end{table}
 By computing $[R_{e_{1},e_{0}},v_{1}\otimes e_{2}e_{3}s]=v_{1}\otimes e_{1}e_{2}e_{3}s$,
$[R_{e_{2},e_{0}},v_{1}\otimes e_{1}s]=v_{1}\otimes e_{1}e_{2}s$,
$[v_{1}\otimes e_{1}s,v_{1}\otimes e_{2}e_{3}s]=-6E$ and $[R_{e_{1},e_{3}},R_{e_{2},e_{-3}}]=-R_{e_{1},e_{2}}$,
we obtain all basis elements for $\mathfrak{g}^{e}(2)$. Hence, we
have that $\mathfrak{g}^{e}(\geq1)$ is generated by $\mathfrak{g}^{e}(1)$.

The above argument also implies that $\mathfrak{g}^{e}(2)\subseteq[\mathfrak{g}^{e},\mathfrak{g}^{e}]$
and $e\in[\mathfrak{g}^{e},\mathfrak{g}^{e}]$. Next we show that
$e$ is strongly reachable. Note that $[R_{e_{-3},e_{0}},R_{e_{3},e_{0}}]=R_{e_{3},e_{-3}}$
and $[R_{e_{1},e_{-2}},R_{e_{2},e_{-1}}]=R_{e_{1},e_{-1}}-R_{e_{2},e_{-2}}$.
The other basis elements for $\mathfrak{g}^{e}(0)$ and $\mathfrak{g}^{e}(1)$
can be obtained by applying $\mathrm{ad}(R_{e_{1},e_{-1}}-R_{e_{2},e_{-2}})$
or $\mathrm{ad}(R_{e_{3},e_{-3}})$ to basis elements for $\mathfrak{g}^{e}(0)\oplus\mathfrak{g}^{e}(1)$.
Hence, we deduce that $\mathfrak{g}^{e}(0)\oplus\mathfrak{g}^{e}(1)\subseteq[\mathfrak{g}^{e},\mathfrak{g}^{e}]$
and thus $\mathfrak{g}^{e}=[\mathfrak{g}^{e},\mathfrak{g}^{e}]$.

School of Mathematics, University of Birmingham, Edgbaston, Birmingham,
B15 2TT, UK
\end{document}